\documentclass[a4paper]{amsart}

\usepackage[]{graphicx}
\usepackage{psfrag}

\newcommand{\ncmd}{\newcommand}
\ncmd{\rencmd}{\renewcommand}
\ncmd{\dspst}{\displaystyle }

\ncmd{\preu}{\noindent \mbox{\bf Proof  }}
\ncmd{\sketch}{\noindent \mbox{\bf Sketch of proof  }}
\ncmd{\preud}[1]
{\noindent \mbox{\bf Proof of \ref{#1}   }}

\ncmd{\ex}{\noindent \mbox{\bf exemple  }}
\ncmd{\exs}{\noindent \mbox{\bf exemples  }}
\ncmd{\fin}{\hspace*{\fill} 
\quad\hbox{\hskip 1pt\vrule width 4pt height 6pt
          depth 1.5pt\hskip 1pt} \medskip }
\ncmd{\lemv}{\noindent \mbox{\bf Lemma  }}
\ncmd{\thmv}{\noindent \mbox{\bf Theorem  }}

\ncmd{\ab}{\mbox{$\mbox{AdS}_{n}\/$ }}
\ncmd{\ei}{\mbox{${\mbox{Ein}}_{n}\/$ }}
\ncmd{\hei}{\mbox{$\widehat{\mbox{Ein}}_{n}\/$ }}
\ncmd{\tei}{\mbox{$\widetilde{\mbox{Ein}}_{n}\/$ }}

\ncmd{\mr}{\mbox{$\widetilde{M}\/$ }}

\ncmd{\slrr}{\mbox{$\widetilde{SL}(2,{\mathbb R}) \times \widetilde{SL}(2,{\mathbb R})\/$ }}

\newtheorem{prop}{Proposition}[section]
\newtheorem{thm}[prop]{Theorem}
\newtheorem{lem}[prop]{Lemma}
\newtheorem{cor}[prop]{Corollary}
\newtheorem{defin}[prop]{Definition}
\newtheorem{rac}[prop]{Remark}
\newtheorem{ric}[prop]{Example}

\ncmd{\exe}{\begin{ric} \em }
\ncmd{\eexe}{\em \end{ric}}
\ncmd{\rque}{\begin{rac} \em}
\ncmd{\erque}{\em \end{rac}}

\rencmd{\thesubsubsection}{\thesubsection-\alph{subsubsection}}
\newcounter{inc}
\rencmd{\theequation}{\arabic{inc}}

\makeindex

\begin{document}

\title[Causal actions and limit sets]{Causal properties of AdS-isometry groups I: \\
Causal actions and limit sets}
\author{Thierry Barbot}
\email{Thierry.Barbot@umpa.ens-lyon.fr}
\thanks
{Work supported 
by CNRS, ACI ``Structures g\'eom\'etriques et Trous Noirs''.}
\address{CNRS, UMR 5669\\
Ecole Normale Sup\'erieure de Lyon\\ 
46 all\'ee d'Italie, 69364 Lyon}
\date{\today}

\begin{abstract}
We study the causality relation in the $3$-dimensional 
anti-de Sitter space AdS and its conformal boundary $\mbox{Ein}_2$. 
To any closed achronal subset $\Lambda$ in $\mbox{Ein}_2$ we associate 
the invisible domain $E(\Lambda)$ from $\Lambda$ in AdS. We show that if
$\Gamma$ is a torsion-free discrete group of isometries of AdS preserving 
$\Lambda$ and is non-elementary (for example, not abelian)
then the action of $\Gamma$ on $E(\Lambda)$ is free, properly 
discontinuous and strongly causal. If $\Lambda$ is a topological
circle then the quotient space 
$M_\Lambda(\Gamma) = \Gamma\backslash{E}(\Lambda)$ is a maximal globally
hyperbolic AdS-spacetime admitting a Cauchy surface $S$ such that the induced
metric on $S$ is complete.

In a forthcoming paper \cite{ba2} we study the case where $\Gamma$ is elementary
and use the results of the
present paper to define a large family of AdS-spacetimes including
all the previously known examples of BTZ multi-black holes.
\end{abstract}

\maketitle

\section{Introduction}
The anti-de Sitter space AdS is the complete Lorentzian manifold with constant 
sectional curvature $-1$ (\S~\ref{sub.ads}). This is the lorentzian version of the hyperbolic 
space ${\mathbb H}^n$. Here we only consider the $3$-dimensional case. In this introduction 
Isom(AdS) denotes the group of isometries of AdS preserving the orientation and the 
time-orientation (see definition~\ref{def.timeori}).

We intend to reach two goals:

\begin{itemize}
\item For every discrete group $\Gamma$ of isometries on AdS, study $\Gamma$-invariant domains of AdS 
on which $\Gamma$ acts properly discontinuously,

\item Provide a general geometrical framework including the notion of BTZ black-hole and multi-black hole. 

\end{itemize}

BTZ black-holes were defined in \cite{BTZ, BTZ2} and studied in many papers,
serving as toy models in the attempt to put together the black-hole notion 
(which arises from general relativity) and quantum physics. We will
discuss these spacetimes with detail in the second part of this work.

The present paper is devoted to the first problem.
We essentially mimic the classical study of groups of isometry of the hyperbolic 
space. Let's recall few basic facts of this well-known theory: let $\Gamma$ be a 
discrete group of isometry of ${\mathbb H}^n$. The action of $\Gamma$ on the entire 
${\mathbb H}^n$ is properly discontinuous and if $\Gamma$ is torsion-free, 
this action is free. Moreover, the action of $\Gamma$ on the boundary at infinity 
$\partial{\mathbb H}^n \approx {\mathbb S}^{n-1}$ admits a unique minimal 
invariant closed subset $\Lambda(\Gamma)$ and the action of $\Gamma$ on 
$\partial{\mathbb H}^n$ is properly discontinuous.

These features do not apply directly in the AdS context, mainly because 
the action of Isom(AdS) on AdS is not proper. 
In this paper we adopt the following point of view: when dealing 
with this kind of questions it is pertinent to take into account 
related \emph{causality notions,\/} which, in the riemannian context, 
remain hidden since automatically fulfilled. Moreover, this causal 
aspect lies in the very foundation of the notion of BTZ black hole. 

\subsection{Causality notions}
A lorentzian manifold  - in short, a spacetime - is causal if no 
point can be joined to itself by a non-trivial causal curve, i.e., a $C^1$ 
immersed curve for which the tangent vectors have non-positive norm for 
the ambient lorentzian metric. Observe that this notion remains meaningful 
for any pseudo-riemannian metric. In the riemannian case this property is 
always true since all tangent vectors have positive norms.

This notion extends to isometry groups in the following way: a group 
$\Gamma$ of isometries on a pseudo-riemannian manifold $\Omega$ 
is causal if for every $x$ in $M$ and every $\gamma$ in $\Gamma$ 
there is no causal nontrivial curve joining $x$ to $\gamma{x}$. 
In the riemannian context any action is causal.

In \S~\ref{sec.defcause} we collect general definitions and facts about causality 
notions in general spacetimes.

\subsection{Einstein universe}
Actually we spend a great amount of time to the 
detailed discussion of the subtle causality notion in AdS-spacetimes: 
see \S~\ref{sub.acausalads}. One important aspect is that in AdS the 
causality relation is trivial: every pair of points in AdS is causally 
related! But this triviality vanishes in the universal 
covering of AdS and we can define \emph{achronal\/} subsets (see next {\S}) 
in AdS as the projections in AdS of achronal subsets 
of the universal covering (see \S~\ref{sec.causaladsein}).

It appears extremely useful for the causal study of AdS to use the conformal embedding of 
the AdS spacetime in the universal conformally flat Lorentzian manifold: the \emph{Einstein universe} 
$\mbox{Ein}_3$ (\S~\ref{sec.conform}). The causality notion persists in the conformal framework 
and the understanding of the causality relation in AdS follows easily from the study in $\mbox{Ein}_3$ 
(\S~\ref{sub.causalein}). This ingredient is particularly useful for the study of spacelike 
surfaces (see \S~\ref{SPA}) in AdS-spacetimes: compare our proof of proposition~\ref{facile!} 
with the similar statement (Lemma $7$)  in the pioneering paper \cite{mess}.

Einstein universe plays another important role: the two-dimensional Einstein universe $\mbox{Ein}_2$ is 
the natural conformal boundary of AdS (remark~\ref{noms}). This feature is completely similar to the 
fact that the natural conformal boundary of ${\mathbb H}^3$ is the universal conformally flat riemannian 
surface; i.e. the round sphere ${\mathbb S}^2$ (observe also 
that ${\mathbb H}^3$ admits a conformal embedding in ${\mathbb S}^3$). The causality notion extends to 
the conformal completion $\mbox{AdS} \cup \partial\mbox{AdS}$ (\S~\ref{adsdads}) simply by restriction 
of the causality relation in $\mbox{Ein}_3$.

\subsection{Invisible domains}
A subset $\Lambda$ of the conformal boundary $\mbox{Ein}_2$ is \emph{achronal\/} if pairs of points in 
$\Lambda$ are not causally related. The \emph{invisible domain from $\Lambda$\/}, denoted $E(\Lambda)$, 
is the set of points in AdS which are not causally related to any point in $\Lambda$ 
(\S~\ref{sub.invisible}). When $\Lambda$ is \emph{generic\/,} i.e. not \emph{pure lightlike\/} 
(which is a particularly exceptional case, see definitions~\ref{defgeneric}, \ref{defgeneric2}), 
then $E(\Lambda)$ is a non-empty convex open domain containing $\Lambda$ in its closure 
which is \emph{geodesically convex,\/} i.e. any timelike geodesic segment in AdS joining 
two points in $U$ is contained in $E(\Lambda)$. 

We study extensively this notion in \S~\ref{sec.cauchy} in the ``non-elementary''
case, i.e. when $\Lambda$ is not contained in one 
or two lightlike geodesics (\S~\ref{sub.elementary}). The elementary case
will be studied in \cite{ba2}. We also consider in detail the case where 
$\Lambda$ is a topological circle in 
$\mbox{Ein}_2$: $E(\Lambda)$ is then \emph{globally hyperbolic,\/} 
with \emph{regular cosmological time\/} (see definitions~\ref{def.gh}, \ref{def.cosmo}, 
\ref{def.cosmoregular} and 
Propositions~\ref{pro.Ecosmologique}, \ref{remplir1}, \ref{pro.E=T}).

The main issue of \S~\ref{sec.cauchy} is Proposition~\ref{decompropre}: in the general case,
invisible domains decompose in (non-disjoint) two globally hyperbolic domains, called globally 
hyperbolic cores, and \emph{closed ends.\/} This notion will be more detailed in the next paper
when we will consider BTZ black-holes. 
We just mention here that they are simply intersections between 
the Klein model ${\mathbb A}{\mathbb D}{\mathbb S}$ and tetraedra in the projective space ${\mathbb R}P^3$.
They are also simple pieces of elementary invisible domains. 

\subsection{Proper and causal actions}
Let $\Gamma$ be a discrete group of isometries of AdS. We say 
that $\Gamma$ is \emph{admissible} if it preserves a generic 
non-elementary achronal subset of $\mbox{Ein}_2$. This definition 
extends to the elementary cases, but requires then a more detailed 
discussion that we postpone to \cite{ba2}. Actually, it is easy to see that 
$\Lambda$ is necessarily non-elementary if $\Gamma$ is not abelian.

However, one of the main result of the present paper is Theorem~\ref{noneleOK} 
that we reproduce here:

\begin{thm}
\label{thm.main1}
Let ${\Lambda}$ be a nonelementary generic achronal subset preserved by a discrete group 
$\Gamma \subset \mbox{Isom(AdS)}$. Then the action of
$\Gamma$ on $E({\Lambda})$ is properly discontinuous.
\end{thm}

Theorem~\ref{noneleOK} states also another result: the quotient spacetime $M_{{\Lambda}}(\Gamma)$ is 
\emph{strongly causal\/} (see \S~\ref{sub.strong}).

\rque
\label{rk.torus}
The theorem above is still true when $\Gamma$ is a 
non-cyclic abelian group:
we then obtain as quotient spaces the (AdS)-\emph{Torus universes\/} 
described in \cite{martinec, carlip}. 
They correspond, through the AdS-rescaling (\cite{BenBon, BenBon2}) 
to the flat Torus Universes described in 
\cite{BBZ}, \cite{carlip}, \cite{BenGua}.
\erque

Furthermore, we have a description of admissible groups (Theorem~\ref{thm.admissible}). 
In the non-abelian case the formulation is a follows: Isom(AdS) is isomorphic to the 
quotient of $\mbox{SL}(2,{\mathbb R}) \times \mbox{SL}(2, {\mathbb R})$ by a subgroup 
of order two (see \S~\ref{sub.psl}). Under this identification every admissible group 
is the projection of the image of a representation 
$\rho = (\rho_L, \rho_R): \Gamma \rightarrow \mbox{SL}(2,{\mathbb R}) \times \mbox{SL}(2, {\mathbb R})$ 
where $\rho_{L, R}$ are fuchsian (i.e. faithfull and discrete) representations into $\mbox{SL}(2, {\mathbb R})$ 
one semi-conjugate to the other, i.e. defining the same bounded Euler cohomology class in 
$H_b^2(\Gamma, {\mathbb Z})$ (see \cite{ghyseuler1, ghyseuler2}, remark~\ref{rk.euler}).

\subsection{Limit sets}
We still assume that $\Gamma$ is admissible and non-abelian. Then, there is
an unique minimal closed generic achronal $\Gamma$-invariant subset $\Lambda(\Gamma)$ 
which is contained in every closed achronal $\Gamma$-invariant subset (Theorem~\ref{thm.minimal}, 
corollary~\ref{cor.cor}).
Hence, following the classical terminology used for isometry groups of ${\mathbb H}^n$, it is natural 
to call $\Lambda(\Gamma)$ the \emph{limit set} of $\Gamma$. 

This analogy with the riemannian case can be pursued further: 
consider the Klein model of ${\mathbb H}^n$ 
as a (convex) ellipsoid in ${\mathbb R}P^n$. Let $\Lambda$ be the limit set 
of a discrete group of isometries $\Gamma$ (a Kleinian group). For any pair $(x,y)$ of points in 
$\Lambda$ let $E_{xy}$ be the unique connected component of ${\mathbb R}P^n \setminus (T_x \cup T_y)$ 
containing ${\mathbb H}^n$ where $T_x$, $T_y$ are the projective hyperplanes tangent to ${\mathbb H}^n$ 
at $x$, $y$. Then the intersection of all the $E_{xy}$ with $(x,y)$ describing $\Lambda \times \Lambda$ 
is a convex domain $E(\Lambda)$ on which $\Gamma$ acts properly: this statement is the true riemannian counterpart of 
Theorem~\ref{thm.main1}.

The complement in ${\mathbb R}P^n$ of the closure of ${\mathbb H}^n$ admits a natural 
$\mbox{SO}(1,n)$-invariant lorentzian metric: this is the (Klein model of) the de Sitter space 
$\mbox{dS}^n$. Then the action of $\Gamma$ on $\mbox{dS}^n \cap E(\Lambda)$ is causal.
When $\Gamma$ is torsion-free the quotient space is strongly causal - even more, it is globally hyperbolic. 
This is a particular case of the content of Scannell thesis (\cite{scannel2}) where maximal globally 
hyperbolic dS-spacetimes are classified.
For more details in this direction we refer to the survey \cite{Bar.Zeg}.

\rque
However there is an important difference between the de Sitter case and the anti de Sitter one: whereas 
Scannell proved that the quotient spaces described above are all globally hyperbolic, the quotient space 
$M_\Lambda(\Gamma) = \Gamma\backslash{E}(\Lambda)$ is \emph{not} globally  hyperbolic, except if $\Lambda$ 
is a topological circle. 
\erque

\subsection{Maximal globally hyperbolic AdS-spacetimes}
We mention here another result of the present paper (corollary~\ref{cor.maxigh}): \emph{the quotient 
spacetimes $M_\Lambda(\Gamma)$, where $\Lambda$ is an achronal topological circle, admits Cauchy 
surfaces such that the ambient AdS-metric restricts as a Cauchy-complete riemannian metric\/} (for 
the definition of Cauchy surfaces, see \S~\ref{s.defgh}). This theorem is also proved in \cite{BenBon2} 
Proposition $6.4.19$. Our proof relies on an general construction interesting by itself
and which presumably works for the the higher dimensional case, associating to 
any embedded spacelike surfaces in AdS 
a surface embedded in ${\mathbb H}^2 \times {\mathbb H}^2$ (\S~\ref{sub.timebundle}, \ref{sub.isomh2}). 
This result completes nicely the classification of maximal globally hyperbolic AdS-spacetimes admitting 
Cauchy-complete Cauchy surfaces (Proposition~\ref{pro.maximalgh}, see also the proposition $6.5.7$ in 
\cite{BenBon2}). Thus the AdS-rescaling defined in \cite{BenBon, BenBon2} establishes a natural bijection 
between these maximal globally hyperbolic AdS-spacetimes and non complete maximal 
Cauchy-complete globally hyperbolic 
\emph{flat} spacetimes which have been classified in \cite{barflat}.

The spacetimes $M_\Lambda(\Gamma)$ defined here, even if $\Lambda$ 
is not a topological sphere, are always strongly causal. In particular 
they are never compact. This is an important difference between our work 
with usual studies of discrete subgroups of Isom(AdS) where a special 
focus is put on the cocompact case. In our framework the spacetimes 
enjoying a ``compact'' character are the globally hyperbolic AdS-spacetimes 
admitting a closed Cauchy surface.

\subsection{Admissible groups}
The main drawback of our approach is that many discrete subgroups are not admissible, i.e. do not preserve 
generic achronal subsets of $\mbox{Ein}_2$ (for example the lattices of Isom(AdS)). In 
\S~\ref{subsec.benoistconvexe} we give a characterization of (non-abelian) admissible groups: 
they are the subgroups of $\mbox{SO}_0(2,2)$ preserving some proper convex domain of ${\mathbb R}P^3$ 
which, in the terminology of \cite{benoistconvex}, are \emph{positively proximal.\/} In some way, the 
present paper provides a geometrical illustration of some cases considered in \cite{benoistconvex}. 
Actually we need to be more precise. The claim above is not exactly correct: it is true that admissible 
groups are positively proximal, but positively proximal subgroups are not always admissible. But our claim 
is not so far to be correct: the Klein model $\overline{{\mathbb A}{\mathbb D}{\mathbb S}}$ in 
${\mathbb R}P^3$ is one connected component of the complement of a quadric.
The other connected component is another copy $\overline{{\mathbb A}{\mathbb D}{\mathbb S}}'$ of AdS 
(see remark~\ref{autreads}, \ref{autretemps}). Observe that $\mbox{SO}_0(2,2)$ is also the isometry group of 
$\overline{{\mathbb A}{\mathbb D}{\mathbb S}}$. Then the correct statement is proposition~\ref{pro.-}: 
any positively proximal subgroup of $\mbox{SO}_0(2,2)$ is either admissible, or admissible as considered 
as a group of isometry of $\overline{{\mathbb A}{\mathbb D}{\mathbb S}}'$.

\subsection{Higher dimension}
Many results in the present paper extend in higher dimensions, particularly if we restrict ourselves 
to strongly irreducible subgroups. But the number of ``elementary'' cases increases with the dimension 
and a systematic treatment requires a non-elementary case-by-case study of these ``elementary'' cases. 
We prefer to postpone such a study to another circumstance. Dropping here the elementary case is 
not conceivable since it corresponds to our second main goal: the systematic description of BTZ 
multi-black holes, in particular, of single BTZ black holes.

\section*{Acknowledgements}
This work would not exist without many conversations with
F. B\'eguin, C. Frances and A. Zeghib.
I'm glad to thank A. Fathi, who clarified for me some aspects of 
the problem of extensions of Lipschitz functions (remark~\ref{deff}, 
proof of proposition~\ref{remplir1}). 
I also benefited from several conversations with F. Bonsante, K. Melnick.
The figures have been computed by P. Saad\'e.

\subsection*{Additionnal comments}
Besides all the references specific to BTZ black-holes themselves, this work is based on many ideas present 
in the unpublished preprint \cite{mess}. 
The elaboration of this paper has been announced in \cite{Bar.Zeg} with the title ``Limit sets of discrete 
Lorentzian groups''.

\section{General notions}
\label{sec.defcause}
A \emph{spacetime} $M$ is a manifold equipped with a lorentzian metric - actually 
we will soon restrict to the constant curvature case.
In our convention a lorentzian metric has signature $(-, +, \ldots, +)$; an 
orthonormal frame
is a frame $(e_1, e_2, \ldots , e_n)$ where $e_1$ has norm $-1$, every 
$e_i \;\; (i \geq 2)$ has norm $+1$ and every scalar product $\langle e_i \mid e_j \rangle \;\; (i \neq j)$ is $0$.
A tangent vector is \emph{spacelike} if its norm is positive; \emph{timelike} if 
its norm is negative; \emph{lightlike\/} if its norm is $0$.  
We also define \emph{causal} vectors as tangent vectors which are timelike 
or lightlike. An immersed surface $S$ is spacelike if all vectors tangent to 
$S$ are spacelike; it is nontimelike if tangent vectors are all spacelike or lightlike.

A causal (resp. timelike) curve is an immersion $c: I \subset {\mathbb R} \rightarrow M$ 
such that for every $t$ in $I$ the derivative $c'(t)$ is causal (resp. timelike). 
This notion extends naturally to non-differentiable curves (see \cite{beem}).
Such a curve is \emph{extendible\/} if there is another causal curve 
$\hat{c}: J \rightarrow M$ and a homeomorphism $\varphi: I \rightarrow K \subset J$ 
such that $K \neq J$ and $c$ coincide with $\hat{c} \circ \varphi$. The causal curve 
$c$ is \emph{inextendible\/} if it is not extendible.

\subsection{Time oriention}
We always assume that the lorentzian manifold is oriented. On spacetimes we have another orientability notion:

\begin{defin}
A spacetime $M$ is time-orientable (or chronologically orientable) if it admits a
a continuous field of timelike vectors.
\end{defin}

\begin{defin}
\label{def.timeori}
A time-orientation on $M$ is an equivalence class of continuous timelike vector fields,
for the following equivalence relation: two timelike vector fields $X$, $Y$ are equivalent
if for any $x$ in $M$ the scalar product $\langle X(x) \mid Y(x) \rangle$ is negative.

$M$ is time oriented when a time-orientation on $M$ has been selected.
\end{defin}

It is easy to show that any spacetime admits a continuous field of timelike lines. Hence,
every spacetime is doubly covered by a time-orientable
lorentzian manifold. We will always assume that
$M$ is time-oriented. Once selected the time-orientation $X$, the set of causal tangent vectors
splits into the union of two bundles of convex cones: the cone of future-oriented vectors 
$\{ v \in T_xM \mid \langle v \mid X(p) \rangle  < 0 \}$ and the cone of past-oriented vectors 
$\{ v \in T_xM \mid \langle v \mid X(p) \rangle  > 0 \}$.

Time-orientation provides naturally an orientation on every causal curve: a causal curve is either 
future-oriented or past-oriented.

\subsection{Causality notions}

Two points in $M$ are causally related if there exists a causal curve joining
them; they are {\em strictly\/} causally related if moreover this joining
curve can be chosen timelike. 

More generally: let $E$ a subset of $M$ and $U$ an open neighborhood of 
$E$ in $M$.
A subset $E$ of $M$ is said \emph{achronal\/} in $U$ if there is no 
timelike curve contained in $U$ joining two points of the subset. 
It is \emph{acausal,\/} or \emph{strictly achronal\/} in $U$ if there is no 
causal curve contained in $U$ joining two points of $E$. We say simply 
that $E$ is (strictly) achronal if it is (strictly) achronal in $U = M$. Finally, 
we say that $E$ is locally (strictly) achronal if every point $x$ 
in $E$ admits a neighborhood $U$ in $M$ such that $E \cap U$ is (strictly)
achronal in $U$.

\rque
\label{spacelocalachronal}
Spacelike hypersurfaces are locally acausal; nontimelike hypersurfaces are locally achronal.
\erque

\subsection{Past, future}
\label{sub.strong}

\begin{defin}
The future of a subset $A$ of $M$ is the set of final points of future oriented timelike curves not 
reduced to one point and starting from a point of $A$. The causal future of $A$ is the set of final 
points of future oriented causal curves, maybe reduced to one point and starting from a point of $S$ 
(hence $A$ itself belongs to its causal future).
The (causal) past of $A$ is the (causal) future of $A$ when the time-orientation of $M$ is reversed.
\end{defin}

\begin{defin}
Let $x$, $y$ be two points in $M$ with $y$ in the future of $x$. The common past-future region 
$U(x,y)$ is the intersection between the past of $y$ and the future of $x$.
\end{defin}

The domains $U(x,y)$ form the basis for some topology on $M$, the so-called \emph{Alexandrov topology\/} 
(see \cite{beem}). Observe that every $U(x,y)$ is open
for the manifold topology. The converse in general is false:

\begin{defin}
If the Alexandrov topology coincide with the manifold topology, $M$ is strongly causal.
\end{defin}

\rque
\label{strongstable}
If $M$ is strongly causal, every open domain $U \subset M$ equipped with the 
restriction of the ambient lorentzian metric is strongly causal.
\erque

\begin{prop}[proposition $3.11$ of \cite{beem}]
\label{pro.strongconvex}
The lorentzian manifold $M$ is strongly causal if and only if it satisfies the
following property: for every point $x$ in $M$ every neighborhood of $x$ 
contains an open neighborhood $U$ (for the usual manifold topology) which 
is \emph{causally convex,\/} i.e. such that any causal curve in $M$ 
joining two points in $U$ is actually contained in $U$.\fin
\end{prop}

\subsection{Global hyperbolicity}
\label{s.defgh}

\begin{defin}
\label{def.gh}
$M$ is globally hyperbolic if:

-- it is strongly causal,

-- for any $x$, $y$ in $M$ the intersection between the causal future of $x$ and 
the causal past of $y$ is compact or empty.
\end{defin}

From now we assume that $M$ is strongly causal.

The notion of global hyperbolicity is closely related to the notion of Cauchy 
surfaces: let $S$ be a spacelike surface embedded in $M$.

\begin{defin}
The past development $P(S)$ (resp. the future development $F(S)$) is the set of points $x$ in $M$ such that every inextendible causal path containing $x$ meets $S$ in its future (resp. in its past). The Cauchy development ${\mathcal C}(S)$ is the union $P(S) \cup F({S})$.
\end{defin}

\begin{defin}
If ${\mathcal C}(S)$ is the entire $M$, $S$ is a Cauchy surface.
\end{defin}

\begin{thm}[\cite{geroch}]
\label{gerochthm1}
A strongly causal lorentzian manifold $M$ is globally hyperbolic if and only if it admits a Cauchy surface. \fin
\end{thm}

\begin{thm}[\cite{geroch}, Proposition 6.6.8 of \cite{hawking}]
\label{thm.geroch}
If $M$ is globally hyperbolic and $S$ a Cauchy surface of $M$, there is a diffeomorphism
$f: M \rightarrow S \times {\mathbb R}$ such that every $f^{-1}(S \times \{ \ast \})$ is a 
Cauchy surface in $M$.\fin
\end{thm}

\rque
\label{rk.smooth}
There has been some imprecision in the litterature concerning the proof 
the smoothness of the splitting of globally hyperbolic spacetimes. 
See \cite{sanchez1, sanchez2, sanchez3} for a survey on this question and a complete proof 
of the smoothness of the splitting $M \approx  S \times {\mathbb R}$.
\erque

\rque
All the notions discussed in this section above only depends on the conformal class of the metric. 
Hence they are well-defined in every conformally lorentzian manifolds, in particular, in 
$\mbox{Ein}_{n}$ (see \S\ref{sec.conform}).
\erque

\begin{prop}
\label{pro.quotientgh}
Let $M$ be a globally hyperbolic spacetime and let $\Gamma$ be a group of isometries of $M$
acting freely and properly discontinuously on $M$ and preserving the chronological orientation. 
Assume that $\Gamma$ preserves a 
Cauchy surface $S$ in $M$. Then, the quotient spacetime $\Gamma\backslash{M}$ is globally hyperbolic 
and the projection of $S$ in this quotient is a Cauchy surface.
\end{prop}

\sketch
Let $M'$ be the quotient $\Gamma\backslash{M}$ and $S'$ be the projection of $S$ in $M'$.
Since $\Gamma$ preserves the chronological orientation, the future (resp. past) of $S'$ in $M'$ is 
the projection of the future (resp. past) of $S$ in $M$.
Every inextendible causal curve in $M'$ lifts in $M$ as a inextendible causal curve in $M$.
The proposition follows.\fin

\subsection{Maximal globally hyperbolic spacetimes}

In this \S, we assume that $M$ has constant curvature\footnote{We could actually only suppose that $(M,g)$ is a solution of the Einstein equation $\mbox{Ricci}_g - \frac{R}{2}g = \Lambda g$.}.

\begin{defin}
An isometric embedding $f: M \rightarrow N$ is a Cauchy embedding if the image by $f$ of any Cauchy surface in $M$ is a Cauchy surface of $N$.
\end{defin}

\begin{defin}
\label{def.maxigh}
A globally hyperbolic spacetime $M$ is maximal if every Cauchy embedding $f: M \rightarrow N$ in a spacetime with constant curvature is surjective.
\end{defin}

\begin{thm}[see Choquet-Bruhat-Geroch \cite{choquet.geroch}]
Let $M$ be a globally hyperbolic spacetime with constant curvature. Then, there is a Cauchy embedding $f: M \rightarrow N$ in a maximal globally hyperbolic spacetime $N$ with constant curvature. Moreover, this maximal globally hyperbolic extension is unique up to right composition by an isometry.\fin
\end{thm}

\subsection{Lorentzian distance}
Let $M$ be a time-oriented spacetime. 

\begin{defin}
The length-time $L(c)$ of a causal curve $c: I \rightarrow M$ is the integral over $I$ of
the square root of $-\langle c'(t) \mid c'(t) \rangle$.
\end{defin}

Observe that it is well-defined, since causal curves are always Lipschitz.

\begin{defin}
The lorentzian distance $d_{lor}(x,y)$ between two points $x$, $y$ in $M$ 
is $\mbox{Sup}\{ L(c) / c \in C(x,y) \}$ where $C(x,y)$ 
is the set of causal curves with extremities $x$, $y$. 
By convention, if $x$, $y$ are not causally related, $d_{lor}(x,y) = 0$.
\end{defin}

\begin{thm}[Corollary $4.7$ and Theorem $6.1$ of \cite{beem}]
\label{th.cosmogh}
If $M$ is globally hyperbolic, then $d: M \times M \rightarrow [0, +\infty]$ 
is continuous and admits only finite values. Moreover, if $y$ is in 
the causal future of $x$, then there exist a geodesic $c$ with extremities $x$, $y$ such that 
$L(c) = d(x,y)$.\fin
\end{thm}

\subsection{Cosmological time}
In any spacetime, we can define the \emph{cosmological time\/} (see \cite{galloway}):

\begin{defin}
\label{def.cosmo}
For any $x$ in $M$, the cosmological time $\tau(x)$ is the $\mbox{Sup}\{ L(c) / c \in {\mathcal R}^-(x) \}$, where ${\mathcal R}^-(x)$ is the set of past-oriented causal curves starting at $x$,
\end{defin}

This function could have in general a bad behavior: for example, in Minkowski space, the cosmological time is everywhere infinite.

\begin{defin}
\label{def.cosmoregular}
$M$ has \emph{regular cosmological time\/} if:

-- $M$ has \emph{finite existence time,\/} i.e. $\tau(x) < \infty$ for every $x$ in $M$,

-- for every past-oriented inextendible curve $c: [0, +\infty[ \rightarrow M$, $\lim_{t \to \infty} \tau(c(t)) = 0$.

\end{defin}

Theorem $1.2$ in \cite{galloway} expresses many nice properties of spacetimes with regular cosmological time function. We need only the following statement:

\begin{thm}
\label{lem.cosmogood}
If $M$ has regular cosmological time, then the cosmological time is Lipschitz 
regular and $M$ is globally hyperbolic. \fin
\end{thm}

Obviously:

\begin{lem}
\label{lem.cosmopreserve}
On any spacetime $M$, every level set of the cosmological time is preserved by the isometry group of $M$.\fin
\end{lem}

\begin{prop}
\label{pro.quotientcosmook}
If $M$ has regular cosmological time and if $\Gamma$ is a group of isometries acting freely and properly discontinuously on $M$, then the quotient space $\Gamma\backslash{M}$ has regular cosmological time.
\end{prop}

\sketch
It follows from the fact that inextendible causal curves in the quotient are projections of  inextendible causal curves in $M$.\fin

\section{Anti-de Sitter and Einstein spaces}
\label{s.adsein}

\subsection{Anti-de Sitter space}
\label{sub.ads}
Let $E$ denote the vector space ${\mathbb R}^4$ equipped with the quadratic form 
$Q = -u^{2}-v^{2}+x_{1}^{2}+ x_{2}^{2}$ (we also denote $E = {\mathbb R}^{2,2}$). 
The associated bilinear form
is denoted by $Q(.,.)$ or $\langle . \mid . \rangle$ as well.
The $3$-dimensional anti-de Sitter space $\mbox{AdS}_{3}$ is the set $\{ Q = -1 \}$, equipped
with the lorentzian metric obtained by restriction of $Q$. We will most of the time drop the index $3$, since here we only consider the three-dimensional case. This space has constant negative curvature. Its isometry group is naturally $O(2,2)$, acting freely transitively on the bundle of orthonormal frames of AdS. Let $SO_{0}(2,2)$ the neutral component of $O(2,2)$ (the so-called orthochronal
component). Fix an orientation of $SO(2) \subset SO_{0}(2,2)$ (the subgroup
preserving the $x_{1}, x_{2}$-coordinates): it defines a time-orientation on AdS. The neutral component $SO_{0}(2,2)$ is precisely the group of isometries of AdS preserving the orientation and the time orientation.

\subsection{Klein models}
\label{sub.klein}
Let $S(E)$ be the half-projectivization of $E$, i.e. the space of rays. It is the double covering of
the projective space $P(E) \approx {\mathbb RP}^3$. We lift to $S(E)$ all the usual notions 
in $P(E)$: for example, a projective line in $S(E)$ is the radial projection of a $2$-plane in $E$.

$S(E)$ is homeomorphic to the sphere ${\mathbb S}^3$. Let ${\mathbb A}{\mathbb D}{\mathbb S}$ be the radial projection of AdS in $S(E)$: we call it the Klein model of AdS. Let $\overline{{\mathbb A}{\mathbb D}{\mathbb S}}$ be the radial projection of AdS in $P(E)$.
Observe that $\overline{{\mathbb A}{\mathbb D}{\mathbb S}}$ is still time-oriented.

The boundary of ${\mathbb A}{\mathbb D}{\mathbb S}$ (resp. $\overline{{\mathbb A}{\mathbb D}{\mathbb S}}$) in $S(E)$ (resp. $P(E)$) 
is the quadric ${\mathcal Q}$
(resp. $\overline{\mathcal Q}$), projection of the zero set
of $Q$. We call it the \emph{Klein boundary\/} of AdS.

The main interesting feature of the Klein model is that geodesics there are connected components of intersections with projective lines. More generally, the totally geodesic subspaces 
in ${\mathbb A}{\mathbb D}{\mathbb S}$ are the traces in ${\mathbb A}{\mathbb D}{\mathbb S}$ of projective subspaces. Projective lines avoiding ${\mathcal Q}$ are timelike
geodesics, projective lines intersecting transversely $\mathcal Q$ (resp. tangent to ${\mathcal Q}$) induce spacelike (resp. lightlike) geodesics.

\rque
\label{autreads}
The complement in $S(E)$ of the closure of ${\mathbb A}{\mathbb D}{\mathbb S}$ can be considered as another copy of anti-de Sitter space: it is the projection of $\{ Q=+1 \}$; this locus, equipped with the restriction of $-Q$, is isometric to AdS. Hence, $S(E)$ is the union of two copies of AdS and of their common boundary.
\erque

\subsection{The $\mbox{SL}(2, {\mathbb R})$-model}
\label{sub.psl}
Consider the linear space $\mbox{gl}(2, {\mathbb R})$ of $2$ by $2$ matrices, equipped 
with the quadratic form $-det$. It is obviously isometric to $(E, Q)$. Hence,
AdS is canonically identified with the group $\mbox{SL}(2,{\mathbb R})$ of $2$ by $2$ matrices with determinant $1$. The actions of $\mbox{SL}(2,{\mathbb R})$ on itself by left and right translations are both isometric actions: we obtain a morphism $\mbox{SL}(2,{\mathbb R}) \times \mbox{SL}(2,{\mathbb R}) \rightarrow SO_0(2,2)$. A dimension argument proves easily that this morphism is surjective. Its kernel is the pair $\{ (id, id), (-id, -id) \}$.

The Klein model $\overline{{\mathbb A}{\mathbb D}{\mathbb S}}$ is canonically identified with $\mbox{PSL}(2,{\mathbb R})$.
The group of orientation and time orientation preserving isometries is $\mbox{PSL}(2,{\mathbb R}) \times \mbox{PSL}(2,{\mathbb R})$, acting by left and right translations.

This Lie group structure on $\overline{{\mathbb A}{\mathbb D}{\mathbb S}}$ provides a natural parallelism on the tangent bundle: if $\mathcal G$ denotes the Lie algebra of $G = \mbox{PSL}(2,{\mathbb R})$, i.e. the algebra of $2$ by $2$ matrices with zero trace, the differential of left translations identify $G \times {\mathcal G}$ with $TG = T\overline{{\mathbb A}{\mathbb D}{\mathbb S}}$. Then, the AdS-norm of a pair $(g, Y)$ is simply
$-det(Y)$.

\subsection{The universal anti-de Sitter space}
\label{sub.universel}
The anti-de Sitter space AdS is homeomorphic to ${\mathbb R}^{2} \times {\mathbb S}^{1}$:
in particular, it is not simply connected. Let $p: \widetilde{\mbox{AdS}} \rightarrow \mbox{AdS}$ be the universal covering. 
The composition of $p$ with the radial projection $\mbox{AdS} \rightarrow \overline{{\mathbb A}{\mathbb D}{\mathbb S}}$ is denoted by $\bar{p}$: this is a universal covering of $\overline{{\mathbb A}{\mathbb D}{\mathbb S}}$. 

Let $\delta$ be a generator of the Galois group of $\bar{p}$, i.e. the group of covering automorphisms: then, $\delta^2$ generates the Galois group of $p$.

\subsection{Affine domains}
\label{sub.affine}

\begin{defin}
Let $x$ be an element of AdS. The affine domain $A(x)$ is the subset of AdS formed by elements $y$ satisfying $\langle x \mid y \rangle < 0$.
\end{defin}

The restrictions of the radial projections of AdS over ${\mathbb A}{\mathbb D}{\mathbb S}$ or $\overline{{\mathbb A}{\mathbb D}{\mathbb S}}$ to affine domains are injective. The images of these projections are also called affine domains.

For any point $x$ of ${\mathbb A}{\mathbb D}{\mathbb S}$, let $x^{\ast}$ be the projection in ${\mathbb A}{\mathbb D}{\mathbb S}$ of the $Q$-orthogonal
hyperplane in $E$ of the direction defined by $x$: we call it the (totally geodesic) hypersurface
dual to $x$. Observe that $x^{\ast}$ has two connected components. 
Every connected component is a spacelike totally geodesic disc in ${\mathbb A}{\mathbb D}{\mathbb S}$, isometric to (the Klein model of) the hyperbolic disc ${\mathbb H}^2$. The boundary 
of $x^{\ast}$ in $S(E)$ is the set of tangency 
between ${\mathcal Q}$ and lightlike geodesics containing ${\mathbb A}{\mathbb D}{\mathbb S}$. Moreover, $x^{\ast}$ is orthogonal to
every timelike geodesic containing $x$. All these geodesics also contain $-x$.

Denote also by $A(x)$ the projection in ${\mathbb A}{\mathbb D}{\mathbb S}$ of 
the affine domain $A(y)$ in AdS, where $x$ is the projection of $y$ in 
${\mathbb A}{\mathbb D}{\mathbb S}$. Observe that $A(x)$ is 
the connected component of ${\mathbb A}{\mathbb D}{\mathbb S} \setminus x^\ast$ 
containing $x$. It is also the intersection between ${\mathbb A}{\mathbb D}{\mathbb S}$ 
and the affine patch
$V(x) = S(E) \cap S(\{ y / \langle y \mid x \rangle < 0 \})$. $V(x)$ admits 
a natural affine structure and is affinely isomorphic to ${\mathbb R}^3$.

\begin{defin}
Let $\tilde{x}$ be an element of $\widetilde{\mbox{AdS}}$. The affine domain $A(\tilde{x})$ is 
the connected component of $p^{-1}(A(p(\tilde{x})))$ containing $\tilde{x}$.
\end{defin}

Affine domains are simple blocks quite easy to visualize, from which 
$\widetilde{\mbox{AdS}}$ can be 
nicely figured out:

-- every affine domain is naturally identified with the interior in ${\mathbb R}^3$ of the 
one-sheet hyperboloid: $\{ (x,y,z) / x^2 + y^2 < 1+z^2 \}$.

-- for any $\tilde{x}$ in $\widetilde{\mbox{AdS}}$, let $A_i$ be the affine domain 
$A(\delta^i \tilde{x})$ (recall that $\delta$ generates the Galois group of $\bar{p}$). 
Then, the affine domains $A_i$ are disjoint $2$ by $2$, $\widetilde{\mbox{AdS}}$ is the union of the closures $\overline{A}_i$ and two such closures $\overline{A}_i$, $\overline{A}_j$ are disjoint, except if $j = i \pm 1$ (keeping away the trivial case $i=j$), in which case $\overline{A}_j \cap \overline{A}_i$ is a totally geodesic surface isometric to a connected component of $p(\tilde{x})^\ast$.

In other words, the universal anti-de Sitter space can be obtained by 
adding up a bi-infinite sequence of affine domains, every affine domain 
being attached to the next one along a copy of the hyperbolic plane.

\subsection{The projectivized timelike tangent bundle}
\label{sub.timebundle}

We will also consider the \emph{projectivized timelike tangent bundle:\/} 
it is the bundle $PT_{-1}\mbox{AdS}$ over AdS admitting as fibers over a point $x$ of AdS 
the set of timelike rays in $T_x\mbox{AdS}$. It can also be defined as the subset ${\mathcal T}$ 
of $\mbox{AdS} \times \mbox{AdS}$ formed by pairs $(x,y)$ such that $\langle x \mid y \rangle = 0$. 
Indeed, for such a pair, the tangent direction at $t=0$ of the curve $t \mapsto {\cos{t}}x+{\sin{t}}y$ 
for $t \geq 0$ defines a timelike ray in $T_x\mbox{AdS}$ and every timelike ray can be obtained 
in this manner in an unique way.

It will be useful to consider $\mathcal T$ as a quadric in $E \times E$. 
For convenience, we write explicitly the definition:

\begin{defin}
\label{def.timebundle}
The projectivized timelike tangent bundle $\mathcal T$ is the set of pairs  $(x,y)$ in $E \times E$ such that:

-- $\mid x \mid = \mid y \mid = -1$,

-- $\langle x \mid y \rangle = 0$.

\end{defin}

Using the canonical parallelism of the vector space $E \times E$, we 
can identify the tangent space to $\mathcal T$ over $(x,y)$ with the vector space of 
vectors $(u,v) \in E \times E$ such that:

-- $\langle x \mid u \rangle = \langle y \mid v \rangle = 0$,

-- $\langle x \mid v \rangle + \langle u \mid y \rangle = 0$.

Define $\mid (u,v) \mid^2$ as the sum $\frac{1}{4}(Q(u) + Q(v))$. It endows $\mathcal T$ with a 
pseudo-riemannian metric. 

The diagonal action of $\mbox{O}(2,2)$ on $E \times E$ preserves $\mathcal T$ and the restriction 
of this action on $\mathcal T$ is isometric for the pseudo-riemannian metric we have just defined. 
We claim that this metric is lorentzian. Indeed, by transitivity of the $\mbox{O}(2,2)$-action, it 
suffices to check at the special point $(x, y) = ((1,0,0,0), (0,1,0,0))$.
Tangent vectors at this point correspond to pairs $(u,v)$, with $u = (0, \alpha, \eta, \nu)$ and 
$v  = (-\alpha, 0, \eta', \nu')$. The pseudo-riemannian norm is therefore 
$-\frac{1}{2}\alpha^2 + \frac{1}{4}(\eta^2+\nu^2+\eta'^2+\nu'^2)$. The claim follows.

Observe that the identification of $\mathcal T$ with $PT_{-1}\mbox{AdS}$ defined above 
is $\mbox{O}(2,2)$-equivariant, where the $\mbox{O}(2,2)$-action on $PT_{-1}\mbox{AdS}$ 
to be considered is the action induced by the differential of its isometric action on AdS.

The space $\mathcal T$ admits two connected components: one of them, called ${\mathcal T}^+$, 
corresponds to future-oriented timelike tangent vectors to AdS and the other, called ${\mathcal T}^-$, 
corresponds to past-oriented timelike tangent vectors.

\section{Conformal embedding of AdS in the Einstein Universe}
\label{sec.conform}

Sometimes (for example, when the causality notion is involved, see next \S), 
it is worth considering the natural embedding of AdS in the so-called 
\emph{Einstein Universe\/} (see \cite{francesthese}).

Let $Q_n$ be the quadratic form $-u^{2}-v^{2}+x_{1}^{2}+ \ldots x_{n}^{2}$ on ${\mathbb R}^{2,n}$ (we only consider here the cases $n=2$ or $n=3$). Let ${\mathcal Q}_n$ be the projection of $\{ Q_n = 0 \}$ is the sphere
$S({\mathbb R}^{2,n})$ of half-directions. As we have seen above, ${\mathcal Q}_2$ can be naturally thought as the (Klein) boundary of AdS.

The $n$-dimensional Einstein universe, denoted by $\mbox{Ein}_{n}$, is the quadric 
${\mathcal Q}_{n}$, equip\-ped with a conformally lorentzian structure 
as follows: let $\pi: {\mathbb R}^{2,n} \setminus \{ 0 \} \rightarrow 
S({\mathbb R}^{2,n})$ be
the radial projection. For any open domain $U$ in ${\mathcal Q}_{n}$ and 
any section $\sigma: U \rightarrow {\mathbb R}^{2,n}$, 
we can define the norm $Q_\sigma(v)$ of any tangent vector $v$ as the $Q_n$ norm 
of $d\sigma(v)$. We obtain by this procedure a lorentzian metric on $U$. 
This lorentzian metric depends on the selected section $\sigma$, but if 
$\sigma' = f\sigma$ is another section, then ${Q}_{\sigma'} = f^2 Q_\sigma$. 
Hence, the {\em conformal class\/} of ${Q}_{\sigma}$ does not
depend on $\sigma$. Moreover, the choice of a representant (of class $C^{r}$) of 
this conformal class is equivalent to the choice of a section (of class $C^{r}$) of $\pi$ 
over $U$.
Eventually, the group of conformal transformations of $\mbox{Ein}_{n}$ is $O(2,n)$.
 
As a first application of this remark, we obtain that $\mbox{Ein}_{n}$ is conformally isomorphic to
${\mathbb S}^{n-1} \times {\mathbb S}^{1}$ equipped with the
metric $ds^{2}-dt^{2}$, where $ds^{2}$ is the usual metric on the unit sphere ${\mathbb S}^{n-1}$,
and $dt^{2}$ the usual metric on ${\mathbb S}^{1} \approx {\mathbb R}/{2\pi\mathbb Z}$. This lorentzian metric
appears when we select the global section $\sigma$ with image contained in the
sphere $u^{2}+v^{2}+x_{1}^{2}+\ldots+x_{n}^{2} = 2$ of ${\mathbb R}^{2,n}$. We will denote by $p: \widehat{\mbox{Ein}}_{n} \approx {\mathbb S}^{n-1} \times {\mathbb R} \rightarrow \mbox{Ein}_{n}$ the cyclic covering 
(it is the universal covering when $n \geq 3$) (it is coherent with the convention in \S~\ref{sub.universel} in view of the natural embedding $\widetilde{\mbox{AdS}} \subset \widehat{\mbox{Ein}}_3$, see remark~\ref{opposite} below).
Observe that ${\mbox{Ein}}_{n}$ is time-orientable.

Throughout this paper, we denote by $d$ the spherical distance on ${\mathbb S}^{n-1}$. 
Keeping in mind the identification $\widehat{\mbox{Ein}}_{n} \approx {\mathbb S}^{n-1}
\times {\mathbb R}$, timelike curves in $\widehat{\mbox{Ein}}_n$ correspond to curves
$t \mapsto (\varphi(t),t)$ where $t$ describes some segment $I \subset {\mathbb R}$ and 
$\varphi: I \rightarrow {\mathbb S}^{n-1}$ is a contracting map
(i.e. the spherical distance $d(\varphi(t),\varphi(t'))$ is
strictly less than $\mid t-t' \mid$). When $\varphi$ is just $1$-lipschitz,
the curve $t \mapsto (\varphi(t),t)$ is only causal.

It is well-known that the notion of lightlike geodesic is still meaningful in the conformally
lorentzian context, but they are not naturally parametrized. More precisely, if we forget
their parametrizations, lightlike geodesics does not depend on the lorentzian metric
in a given conformal class. Under the identification $\widehat{\mbox{Ein}}_{n} \approx {\mathbb S}^{n-1}
\times {\mathbb R}$, inextendible lightlike geodesics are curves 
$t \mapsto (\varphi(t),t)$ where $\varphi: {\mathbb R} \rightarrow {\mathbb S}^{n-1}$
is a geodesic on the sphere.

\begin{thm}(\cite{francesliouville})
The Einstein space $\widehat{\mbox{Ein}}_{n}$ is universal in the category of locally conformally
flat lorentzian spaces; i.e. every simply connected lorentzian manifold of dimension $n$ which is conformally
flat can be conformally immersed in $\widehat{\mbox{Ein}}_{n}$.\fin
\end{thm}

We will not give a proof of this theorem here, but will exhibit 
the natural embedding of anti-de Sitter space AdS in $\mbox{Ein}_3$:
let $v$ be any spacelike vector in ${\mathbb R}^{2,3}$; and $v^{\perp}$ its $Q_3$-orthogonal 
hyperplane in ${\mathbb R}^{2,3}$. Let ${\mathcal A}(v)$ be the projection in 
$S({\mathbb R}^{2,3})$ of the intersection
between $\{ Q_3 = 0 \}$ and one connected component 
of the complementary part of $v^{\perp}$ and 
let $\partial {\mathcal A}(v)$
be the projection of $\{ Q_{3} = 0 \} \cap v^{\perp}$. 
The notations $\mathcal A$ and $\partial {\mathcal A}$
will be reserved to the special case $v = v_{0} = (0, \ldots, 0,1)$.
There is a natural section $\sigma$ over ${\mathcal A}(v)$: take $\sigma(x)$ such that the $Q_{3}$-scalar
product between $\sigma(x)$ and $v$ is equal to $\pm 1$.
A straightforward computation shows that $({\mathcal A}(v), \bar{Q}_\sigma)$ is isometric to AdS.

\begin{defin}
An anti-de Sitter domain in ${\mbox{Ein}}_3$ is an open domain ${\mathcal A}(v)$ for some spacelike vector $v$ in ${\mathbb R}^{2,3}$.
\end{defin}

\rque
\label{2copies}
The notation is a little misleading since every $v$ 
defines actually {\em two\/} domains, since
${\mathbb R}^{2,3} \setminus v^{\perp}$ has two connected components. We can withdraw 
this undeterminancy by defining more precisely ${\mathcal A}(v)$ as the radial 
projection of $\{ x / \langle x \vert v \rangle = -1 \}$. Anyway, both connected 
components are conformal copies 
of AdS, glued along their common conformal boundary $\approx \mbox{Ein}_2$. 
This decomposition is \emph{not} the decomposition discussed in the remark~\ref{autreads}: 
indeed, $\mbox{Ein}_3 \approx {\mathbb S}^2 \times {\mathbb S}^1$ is not homeomorphic 
to $S(E) \approx {\mathbb S}^3$! See also remark~\ref{autretemps}.
\erque

\rque
\label{rk.stereo}
Similar constructions can be performed even when $v$ is not spacelike. When $v$ is 
timelike, we obtain an open domain ${\mathcal S}(v)$  conformally isometric to de $3$-dimensional de Sitter space. When $v$ is lightlike, this procedure provides a conformal idendification, called \emph{stereographic projection,\/} between every connected component of the complement of a lightcone with the three dimensional Minkowski space.
\erque

\rque
\label{scan}
In order to get a satisfying understanding of the geometry involved, it is useful 
to ``scan'' these domains in ${\mathbb S}^{2} \times {\mathbb S}^{1}$: the standard anti-de Sitter domain
$\mathcal A$ is the domain in ${\mathbb D}^{2} \times {\mathbb S}^{1}$ is 
${\mathbb S}^{2} \times {\mathbb S}^{1}$ , where
${\mathbb D}^{2}$ is the upper hemisphere $\{ x_{2} > 0 \}$.
A typical de Sitter domain is ${\mathbb S}^{2} \times ]0,\pi[$.
A typical affine domain (\S~\ref{sub.affine}) of AdS conformally embedded 
in ${\mbox{Ein}}_3$ is ${\mathbb D}^2 \times ]0, \pi[$: this is the intersection between an 
anti-de Sitter domain and some de Sitter domain.

Finally, the conformal embedding of AdS in $\mbox{Ein}_3$ lifts to a conformal embedding of
$\widetilde{\mbox{AdS}}$ in $\widehat{\mbox{Ein}}_3$: it follows that $\widetilde{\mbox{AdS}}$ 
is conformally equivalent to ${\mathbb D}^2 \times {\mathbb R}$ equipped with the metric $ds^2 - dt^2$, 
where $ds^2$ is the restriction to the hemisphere ${\mathbb D}^2$ of the spherical metric.

More precisely, if $O$ is the ``North pole'' of ${\mathbb D}^2$
(i.e. the unique point in ${\mathbb D}^{2}$ which is at $d$-distance
$\pi/2$ to every element of $\partial{\mathbb D}^{2}$), then $\widetilde{\mbox{AdS}}$ is isometric to 
${\mathbb D}^2 \times {\mathbb R}$ equipped with the metric:

\[ \frac{1}{\cos^2(d(p,O))}(ds_0^2 - dt^2) \]

\erque

\rque
\label{patchds}
Sometimes, we will consider the projection of $\mbox{Ein}_{n}$ in 
$P({\mathbb R}^{2,n}) \approx {\mathbb RP}^{n+1}$:
we denote it by $\overline{\mbox{Ein}}_{n}$. It is still time-orientable. Observe that de Sitter domains in 
$\mbox{Ein}_{n}$ projects
injectively in $\overline{\mbox{Ein}}_{n}$ as complementary parts of spacelike hyperplanes.
It follows that de Sitter domains in $\mbox{Ein}_{n}$ correspond to intersections
of $\overline{\mbox{Ein}}_{n}$ with affine patches of ${\mathbb RP}^{n+1}$ with spacelike boundaries. 
\erque

\rque
\label{noms}
When $v$ is spacelike, $\partial{\mathcal A}(v)$ is a copy of $\mbox{Ein}_{2}$: we call it a
{\em Einstein flat subspace.\/} We can also see $\partial{\mathcal A}(v)$ as the {\em conformal boundary of ${\mathcal A}(v)$.\/} The closure 
$\overline{\mathcal A}(v) = {\mathcal A}(v) \cup \partial{\mathcal A}(v)$ is canonically identified with the closure ${\mathbb A}{\mathbb D}{\mathbb S} \cup \mbox{Ein}_{2}$ 
of ${\mathbb A}{\mathbb D}{\mathbb S}$ in $S(E)$, but this identification is not analytic, 
even if its restrictions to respectively  ${\mathcal A}(v)$, 
$\partial{\mathcal A}(v)$ are individually analytic.
\erque

\rque
\label{opposite}
The projection of $\widehat{\mbox{Ein}}_{n}$ over $\overline{\mbox{Ein}}_{n}$ is a cyclic
covering. Denote by $\delta$ a generator of the group of covering transformations.
It is coherent with the convention in \S~\ref{sub.universel} since the restriction 
of this covering transformation to the image of the natural embedding $\widetilde{\mbox{AdS}} \subset \widehat{\mbox{Ein}}_3$ is indeed the generator of the Galois group of $\bar{p}$.

We say that two points in $\widehat{\mbox{Ein}}_{n}$ are {\em opposite\/} if
one of them is the image under $\delta$ of the other. It is easy to give another
equivalent definition: for every element $x$ of $\widehat{\mbox{Ein}}_{n}$, the lightlike
geodesics containing $x$ admits many other common intersection points,
that are precisely the iterates of $x$ under $\delta$. Therefore, two
points $x$, $y$ are opposite if every lightlike geodesic containing one of them
contains the other and if the lightlike segments joining $x$ to $y$ are all disjoint.
\erque

\rque
\label{paropposite}
The Einstein space $\widehat{\mbox{Ein}}_{n}$ admits of course many different parametri\-zations by 
${\mathbb S}^{n-1} \times {\mathbb R}$ for which the conformal lorentzian structure is
represented by $ds^{2}-dt^{2}$. Anyway, 
in all these parametrizations, pair of opposite points 
always have coordinates of the form $(x, \theta)$, $(-x,\theta+\pi)$.
\erque

\rque
\label{rk.einpsl}
According to \S~\ref{sub.psl},  
$\overline{\mbox{Ein}}_2$ is naturally identified with the projection in 
$P(\mbox{gl}(2, {\mathbb R}))$ of nonzero non-invertible matrices. 
For any such matrice $A$, let $K(A)$ be the kernel of $A$ and $I(A)$ 
the image of $A$. Then, $A \mapsto ([I(A)], [K(A)])$ identifies
$\overline{\mbox{Ein}}_2$ with ${\mathbb R}P^1 \times {\mathbb R}P^1$. 
Denote by ${\mathbb R}P^1_L$ (resp. ${\mathbb R}P^1_R$) the leaf space 
of the \emph{left\/} (resp. \emph{right\/}) foliation $\overline{\mathcal G}_L$ 
(resp. $\overline{\mathcal G}_R$), i.e. the foliation of 
${\mathbb R}P^1 \times {\mathbb R}P^1$ with leaves 
$\{ \ast \} \times {\mathbb R}P^1$ (resp. ${\mathbb R}P^1 \times \{ \ast \}$). 
An usual point of view is to consider every point of $\overline{\mbox{Ein}}_2$ 
as the intersection between a leaf of the left foliation and a leaf of the right foliation. 
In this spirit, we can write
$\overline{\mbox{Ein}}_2 = {\mathbb R}P^1_L \times {\mathbb R}P^1_R$. 

The extension of the isometric action of $G \times G$ 
(with $G = \mbox{PSL}(2, {\mathbb R})$) to 
$\overline{\mbox{Ein}}_2$ corresponds to 
the diagonal action of $G \times G$ on ${\mathbb R}P^1 \times {\mathbb R}P^1$, 
where the action of $G$ on ${\mathbb R}P^1$ is the usual projective action. 
Let $G_L$ (resp. $G_R$) be the group of left (resp. right) translations of 
$G$ on itself: the group of conformal isometries of $\overline{\mathcal Q}$ 
is then the product $G_L \times G_R$. Observe that the leaves of the left and
right foliations are lightlike geodesics.

$\mbox{Ein}_2$ can be considered in a similar way: 
it is bifoliated by two transverse foliations ${\mathcal G}_L$, ${\mathcal G}_R$. 
Every leaf of the left or right foliation is canonically the double covering 
$\widehat{{\mathbb R}P}^1$ of ${\mathbb R}P^1$. 
Every leaf of ${\mathcal G}_L$ intersects every leaf of ${\mathcal G}_R$ at 
\emph{two} points, one opposite to the other. Finally, the 
leaf space of the left (resp. right) foliation is canonically 
identified to ${\mathbb R}P^1_L$ (resp. ${\mathbb R}P^1_R$).

This description lifts to $\widehat{\mbox{Ein}}_2$: it admits 
a pair of transverse foliations $\widehat{\mathcal G}_L$, 
$\widehat{\mathcal G}_R$. The respective leaf spaces are still projective lines 
${\mathbb R}P^1_L$, ${\mathbb R}P^1_R$. The leaves themselves are universal 
coverings of the projective line. Finally, the intersection between a leaf 
of the right foliation and a leaf of the left foliation is a $\delta$-orbit.
\erque

\section{Causality relation}
\label{sec.causaladsein}

In this section, we discuss the notion of causality in AdS and $\mbox{Ein}_n$.
A fundamental observation is that the causality relation in AdS and $\mbox{Ein}_n$
is trivial: for any pair $(x,y)$ in AdS or $\mbox{Ein}_n$, there is a timelike 
curve joining $x$ and $y$! Actually, this notion is interesting only in the universal 
covering $\widetilde{\mbox{AdS}} \approx {\mathbb D}^2 \times {\mathbb R}$ and 
in the cyclic covering 
$\widehat{\mbox{Ein}}_{n} \approx {\mathbb S}^{n-1} \times {\mathbb R}$. 

The main purpose of this {\S} is to show that, even if the achronality notion 
is not \emph{stricto sensu} well defined in AdS or in $\mbox{Ein}_n$, projections 
in these spaces of achronal subspaces of $\widetilde{\mbox{AdS}}$ or of 
$\widehat{\mbox{Ein}}_{n}$ are nicely described.

\subsection{Achronality in $\widehat{\mbox{Ein}}_n$}
\label{sub.causalein}

We leave as an exercise to the reader the following lemma (we just stress out 
that any causal curve can intersect every ${\mathbb S}^{n-1} \times \{ \ast \}$ in at most one point):

\begin{lem}
\label{pop}
Two points $(x_{1}, \tilde{\theta}_{1})$ and $(x_{2}, \tilde{\theta}_{2})$ 
are causally related in $\widehat{\mbox{Ein}}_{n} \approx {\mathbb S}^{n-1} 
\times {\mathbb R}$ if and only if the distance in ${\mathbb S}^{n-1}$ between 
$x_{1}$ and $x_{2}$ is less or equal to $\vert \theta_{2}-\theta_{1} \vert$. These points are strictly causally related if the
distance between $x_{1}$ and $x_{2}$ is less to $\vert \theta_{2}-\theta_{1} \vert$.
In particular, they are necessarily causally related if $\vert \theta_{2}-\theta_{1} \vert$
is greater than $\pi$.\fin
\end{lem}

\begin{cor}
\label{graphein}
Achronal subsets of $\widehat{\mbox{Ein}}_n$ are graphs of $1$-lipschitz functions $f: E \rightarrow {\mathbb R}$, where $E$ is a subset of ${\mathbb S}^{n-1}$. Such a subset is strictly achronal if and only if $f$ is contracting.\fin
\end{cor}

\begin{cor}
\label{fermons}
The closure of an achronal subset is achronal.\fin
\end{cor}

\begin{cor}
\label{einstrong}
$\widehat{\mbox{Ein}}_{n}$ is strongly causal.\fin
\end{cor}

\begin{cor}
Every ${\mathbb S}^{n-1} \times \{ \ast \}$ is a Cauchy hypersurface in $\widehat{\mbox{Ein}}_{n}$. In particular, $\widehat{\mbox{Ein}}_{n}$ is globally hyperbolic.\fin
\end{cor}

\begin{lem}
\label{Edansds}
Every closed achronal subset $\Lambda$ of $\widehat{\mbox{Ein}}_{n}$ is contained in a 
de Sitter domain, except if it is contained in the past lightcone and future lightcone of two opposite elements
of itself (see remark~\ref{opposite}).
\end{lem}

\preu
Let $(x^{+},\theta^{+})$, 
$(x^{-},\theta^{-})$ be elements of $\Lambda$ where the $\theta$-coordinate attains respectively its
maximum and mimimum value: if
$\theta^{+}-\theta^{-}$ is strictly less than $\pi$, then the lemma is proved,
since $\Lambda$ is contained in some de Sitter domain of the form 
${\mathbb S}^{n-1} \times ]-\theta^{-}-\epsilon, \theta^{+}+\epsilon[$, with $2\epsilon$ less
than $\pi-\theta^{+}+\theta^{-}$.
On the other hand, $\theta^{+}-\theta^{-}$ is less than $\pi$ since $\pi$ is the diameter of the
hemisphere.
Hence, we have only to deal with the case $\theta^{+}-\theta^{-} = \pi$.
In this case, the distance between $x^{+}$ and $x^{-}$ on the sphere has to
be precisely $\pi$ and $(x^{+},\theta^{+})$, 
$(x^{-},\theta^{-})$ are opposite points (in the meaning of definition \ref{opposite}).
Moreover, for any $(x,\theta)$ in $\Lambda$, 
$x$ lies on a minimizing geodesic of ${\mathbb S}^{n-1}$ between
$x^{+}$ and $x^{-}$. It follows that $\theta$ must be equal to 
$\theta^{+}-d(x,x^{+}) = \theta^{-}+d(x,x^{-})$.
In other words, $\Lambda$ is contained in the past lightcone of $(x^{+},\theta^{+})$ 
and the future lightcone of $(x^{-},\theta^{-})$.\fin

The particular case appearing in lemma \ref{Edansds} deserves a particular appellation.

\begin{defin}
\label{defgeneric}
A subset of $\widehat{\mbox{Ein}}_{n}$ is pure lightlike if it is
contained in the past lightcone and future lightcone of two opposite
elements of itself. If not, it is generic. 
\end{defin}

We point out the obvious fact that a strictly achronal set is generic.

\rque
\label{purelight}
The proof of \ref{Edansds} actually shows that an achronal subset is pure
lightlike as soon as it contains two opposite points.
\erque

\subsection{Achronality in $\mbox{Ein}_n$}
\label{sub.einachronal}

\begin{defin}
An achronal (resp. strictly achronal) subset of ${\mbox{Ein}}_n$ is the projection of an achronal (resp. strictly) achronal subset of $\widehat{\mbox{Ein}}_n$.
\end{defin}

Our main purpose here is to provide an effective
criterion recognizing achronal subsets of ${\mbox{Ein}}_n$. 

Observe that the scalar product $\langle [x] \mid [y] \rangle$ of two elements 
of the sphere $S({\mathbb R}^{2,n})$ of half-rays is not well-defined
but has a well-defined sign.

\begin{prop}
\label{pro.accro}
A subset $\Lambda$ of $\mbox{Ein}_n$ is (strictly) achronal if and only if
it is contained in the closure of a de Sitter domain and 
for every pair $([x],[y])$ of elements of $\Lambda$ the scalar product
$\langle [x] \mid [y] \rangle$ is non-positive (resp. negative).
\end{prop}

\preu
We detail the proof only in the achronal case; the strictly achronal case being similar.

Assume first that $\Lambda$ is achronal, i.e. is the projection of an achronal subset 
$\widehat{\Lambda}$ of $\widehat{\mbox{Ein}}_n$. If $\widehat{\Lambda}$ is pure 
lightlike the conclusion holds quite easily; all the scalar products $\langle [x] \mid [y] \rangle$
are then $0$. Assume now that $\widehat{\Lambda}$ is generic.
According to lemma~\ref{Edansds}, 
$\Lambda$ has to be contained in some de Sitter domain ${\mathcal S}(v_0)$ (with $v_0$ 
a timelike vector of ${\mathbb R}^{2,n}$). In this de Sitter domain it appears clearly
that if $\langle [x] \mid [y] \rangle$ is positive then the segment in the affine domain
${\mathcal S}(v_0)$ with extremities $[x]$, $[y]$ is timelike 
(check in obvious cases and use the transitivity 
of the action of the stabilizer of $v_{0}$ in $SO(2,n)$ on the sets of timelike, lightlike and spacelike lines in ${\mathcal S}(v_0)$). The first implication of the proposition then follows.

For the reverse implication: assume that 
$\Lambda$ is contained in the closure of a de Sitter domain
and that for every pair $([x], [y])$ of elements of 
$\Lambda$ the scalar product $\langle [x] \mid [y] \rangle$ is non-positive. 
Consider any connected component $\widetilde{\mathcal S}(v_0)$
of the preimage of ${\mathcal S}(v_0)$
and let $\widehat{\Lambda}$ be the preimage of $\Lambda$
in the closure $\overline{\mathcal S}(v_0)$ of $\widetilde{\mathcal S}(v_0)$: 
the projection of $\widehat{\Lambda}$ in $S({\mathbb R}^{2,n})$
is $\Lambda$. Using the conformal identification $\overline{\mathcal S}(v_0)
\approx {\mathbb S}^{n-1} \times [-\pi/2, +\pi/2]$ it is quite straightforward 
to check that $\widehat{\Lambda}$ is achronal (hint: de Sitter domains 
in $\widehat{\mbox{Ein}}_n$ are causally convex).
\fin

The way to recognize projections of pure lightlike subsets of $\widehat{\mbox{Ein}}_n$ is obvious:

\begin{defin}
\label{defgeneric2}
A closed subset $\Lambda$ of ${\mbox{Ein}}_{n}$ is pure lightlike if:

-- it contains two opposite points $[x_0]$, $-[x_0]$,

-- for every element $[x]$ of $\Lambda$ the scalar product $\langle [x] \mid [x_0] \rangle$
is zero,

-- $\Lambda$ is contained in the closure of a de Sitter domain.

Non pure lightlike closed subsets of ${\mbox{Ein}}_n$ are 
generic.
\end{defin}

\begin{cor}
\label{kkk}
A generic subset of $\mbox{Ein}_n$ is (strictly) achronal if and only if it is
contained in some de Sitter domain and (strictly) achronal in every de Sitter domain containing it.
\fin
\end{cor}

We have a well-defined notion of convex hull in $S({\mathbb R}^{2,2}) = S(E)$.

\begin{cor}
A generic subset of $\mbox{Ein}_2$ is achronal if and only if its convex hull in
$S(E)$ is contained in the closure in $S(E)$ of the Klein model 
${\mathbb A}{\mathbb D}{\mathbb S}$.
\end{cor}

\preu
Immediate corollary of Proposition~\ref{pro.accro}.\fin

Remember that {\em extreme points\/} of a closed convex $C$ set are points
which do not belong to segments $]x,y[$ with $x$, $y$ in $C$:

\begin{lem}
\label{concau}
An achronal subset of $\mbox{Ein}_2$ is  strictly achronal precisely when the intersection points between its convex hull $C$ and $\mbox{Ein}_{2}$ are all extreme points of $C$.
\end{lem}

\preu
A non-extreme point of the convex hull belonging
to $\mbox{Ein}_{2}$ is the projection of a sum $\sum_{i=1,..,k} t_{i}u_{i}$ 
where every $u_{i}$ projects to an element of $\Lambda$, $0 < t_{i} < 1$, 
$k \geq 2$ and $[u_i] \neq [u_j]$ if $i \neq j$. Then, 
$\langle \sum t_{i}u_{i} \vert \sum t_{i}u_{i} \rangle = 
\sum t_{i}t_{j} \langle u_{i} \vert u_{j} \rangle$
can be zero only if all products $\langle u_{i} \vert u_{j} \rangle$ 
are null, which precisely means that every $u_i$ is causally related to every $u_j$.
\fin

\rque
\label{parametrons}
Let $y$, $y'$ be two non-causally related points in
$\widehat{\mbox{Ein}}_{n}$. Let $p$, $p'$ be two points in ${\mathbb R}^{2,n}$ 
such that $p(y) = [p]$, $p(y') = [p']$.
Since $y$ and $y'$ are not causally related, according to
proposition~\ref{pro.accro}, the quantity $\langle p \vert p' \rangle$ is negative: we can select $p$, $p'$ such that this quantity is
actually $-2$. Then, there is a basis of ${\mathbb R}^{2,n}$ for which the quadratic form $Q_n$ still admits the expression $-u^2 - v^2 + x_1^2 + x_2^2 + \ldots + x_n^2$ and for which
the coordinates of $p$ and $p'$ are respectively
$(1, 0, 1, 0, 0, \ldots , 0)$ and $(1, 0, -1, 0, 0, \ldots , 0)$.
For this choice of coordinates, when we select the section $\sigma: {\mbox{Ein}}_n
\rightarrow {\mathbb R}^{2,n}$ 
taking value in the sphere $u^2 + v^2 + x_1^2 + \ldots + x_n^2 = 2$, we obtain
an identification ${\mbox{Ein}}_n \approx {\mathbb S}^{n-1} \times {\mathbb R}$ where
the conformal structure is still represented by $ds^2 - dt^2$
but now with the additional requirement that 
$p$, $p'$ have coordinates $(x,0)$, $(-x, 0)$. 
\erque

\subsection{Achronality in AdS}
\label{sub.acausalads}

Keeping in mind the identification of $\widetilde{\mbox{AdS}}$ with ${\mathbb D}^2 \times {\mathbb R} \subset {\mathbb S}^2 \times {\mathbb R} \approx \widehat{\mbox{Ein}}_3$:

\begin{lem}
\label{adsachronal}
Two points in $\widetilde{\mbox{AdS}}$ are (strictly) causally related if and only if they are (strictly) causally related in $\widehat{\mbox{Ein}}_3$.
\end{lem}

\preu
Let $\tilde{x}$, $\tilde{y}$ be two points in $\widetilde{\mbox{AdS}}$. Clearly, if they are (strictly) causally related in $\widetilde{\mbox{AdS}}$, they are (strictly) causally related in
$\widehat{\mbox{Ein}}_3$. The inverse follows directly from Lemma~\ref{pop}.\fin

According to Remark~\ref{strongstable} and Corollary~\ref{einstrong}:

\begin{lem}
$\widetilde{\mbox{AdS}}$ is strongly causal.\fin
\end{lem}

Corollaries \ref{graphein} and \ref{fermons} imply:

\begin{lem}
Achronal subsets of $\widetilde{\mbox{AdS}}$ are graphs of $1$-lipschitz functions 
$f: E \rightarrow {\mathbb R}$, where $E$ is a subset of 
${\mathbb D}^{2}$. Such a subset is strictly achronal if and only 
if $f$ is contracting.\fin
\end{lem}

\begin{lem}
The intersection between the closure of an achronal subset of 
$\widetilde{\mbox{AdS}}$ and $\partial \widetilde{\mbox{AdS}}$ 
is an achronal subset of $\widehat{\mbox{Ein}}_2$.
\fin
\end{lem}

Define generic subsets of $\widetilde{\mbox{AdS}}$ as subsets for which the intersection between their closure and $\partial \widetilde{\mbox{AdS}}$ is a generic in $\widehat{\mbox{Ein}}_2$. 
Lemma~\ref{Edansds} now becomes:

\begin{lem}
\label{Edansds2}
Generic achronal subsets of $\widetilde{\mbox{AdS}}$ are contained in affine domains.
\fin
\end{lem}

Define (strictly) achronal domains of AdS as projections of (strictly) achronal domains of 
$\widetilde{\mbox{AdS}}$.
Proposition~\ref{kkk} now becomes:

\begin{prop}
\label{kkk2}
A generic subset of AdS is (strictly) achronal if and only if it is
contained in some affine domain and (strictly) achronal in every affine domain containing it.
\fin
\end{prop}

\subsection{Causality relation between AdS and $\partial\mbox{AdS}$}
\label{adsdads}
Since $\mbox{Ein}_2$ is the boundary of AdS in $\mbox{Ein}_3$, we have the notion that points in AdS can be causally related to points in $\mbox{Ein}_2 \approx \partial\mbox{AdS}$. 
This notion can be easily understood by considering the identification between AdS and ${\mathbb D}^2 \times {\mathbb S}^1$: the boundary of AdS is then $\partial{\mathbb D}^2 \times {\mathbb S}^1 \approx \mbox{Ein}_2$. 
Then, $(x, \theta)$ in ${\mathbb D}^2 \times {\mathbb S}^1$ is causally related to $(y, \theta')$ in $\partial{\mathbb D}^2 \times {\mathbb S}^1$ if and only if $d(x,y) \leq \mid \theta - \theta' \mid$. 

The conformal completion ${\mathcal A} \cup \partial{\mathcal A}$ of AdS is 
naturally identified with the Klein completion 
${\mathbb A}{\mathbb D}{\mathbb S} \cup \mbox{Ein}_2$ 
(see remark~\ref{noms}).
It is also useful to understand the causality relation between points of AdS 
and points in $\partial\mbox{AdS}$, but when the last one is considered as 
the Klein boundary, not the conformal boundary:

\begin{lem}
Let $x$ be an element of ${\mathbb A}{\mathbb D}{\mathbb S}$ 
and let $y$ be an element of the Klein boundary $\mbox{Ein}_2$.
Let $I = ]x,y[$ be the shortest segment between $x$ and $y$ 
(observe that $x$ and $y$ are never opposite in $S(E)$). 
Then, $x$ and $y$ are causally related if and only if $I$ is contained in 
${\mathbb A}{\mathbb D}{\mathbb S}$. There are strictly causally 
related if and only if the projective line $d$ containing $I$ is not lightlike, 
i.e. is transverse to $\mbox{Ein}_2$. If $d$ is tangent to 
$\mbox{Ein}_2$, then $I$ is a lightlike geodesic: considered as elements of 
$\mbox{Ein}_3$, $x$ and $y$ are joined by a lightlike geodesic. \fin
\end{lem}

\rque
\label{autretemps}
In remark~\ref{2copies}, we have observed that $\mbox{Ein}_3$ is obtained by glueing conformally along their boundaries two copies of AdS. In remark~\ref{autreads}, we have seen
that $S(E)$ can also be considered as the union of two copies of AdS, glued along their common Klein boundary. But there is a main difference here: ${\mathbb A}{\mathbb D}{\mathbb S}$ is the projection of $\{ Q=-1 \}$ equipped with the restriction of $Q$, whereas the complement in $S(E)$ of its closure is the projection of $\{ Q=1 \}$ equipped with the restriction of $-Q$. Hence, the identification between the boundaries of these copies of AdS does not preserve the causality notion: it sends causal curves to achronal topological circles!
\erque

\section{Dualities}
\label{dual}
Let $E^{\ast}$ be the dual of $E$. The quadratic form $Q$ defines a map $\flat : E \rightarrow E^{\ast}$
by $x^{\flat}(y) = Q(x,y)$. The image under $\flat$ of $Q$ is a quadratic form
$Q^{\ast}$ on $E^{\ast}$.
Let $S(E)$ and $S(E^{\ast})$ be the associated half-projective spaces. 
The map $\flat$ induces a polarity $S(E) \rightarrow S(E^{\ast})$ that we still denote
by $\flat$. We denote by $\sharp$ the inverse map of $\flat$.

We have denoted by ${\mbox{Ein}}_{2}$ the projection of the null cone of $Q$ in $S(E)$;
the nullcone of $Q^{\ast}$ is a dual copy ${\mbox{Ein}}_{2}^{\ast}$
of ${\mbox{Ein}}_{2}$. Elements of ${\mbox{Ein}}_{2}^{\ast}$ can be interpreted
as lightcones in ${\mbox{Ein}}_{2}$; more precisely, $x^{\flat}$ is the lightcone emitted
from $x$ in ${\mathbb A}{\mathbb D}{\mathbb S}$.

Observe that $\flat$ and $\sharp$ respect the causality notion.
In particular:

\begin{lem}
The image ${\Lambda}^{\flat}$ of a (strictly) achronal subset ${\Lambda}$
of ${\mbox{Ein}}_{n}$ by $\flat$ is (strictly) achronal.\fin
\end{lem}

We will use another notion of duality, more traditional, completely
independant from the notion discussed above: the duality of
convex subsets of $S(E)$. We recall basic facts (cf. \cite{metriquehilbert}):

A convex cone $J$ of $E$ is a convex subset stable by positive homotheties. It is {\em proper\/}
if it is nonempty and its closure $\bar{J}$ does not contain a complete affine line.
A convex subset $C$ of $S(E)$ is the projection of a convex cone $J(C)$ of $E$.
It is called proper if $J(C)$ is proper. 

For any convex cone $J$, we define its dual by 
$J^{\ast} = \{ \alpha \in E^{\ast} / \;\forall x \in \bar{J} \setminus \{ 0 \}\;\; \alpha(x) < 0 \}$.
This provides a construction of dual convex $C^{\ast} \subset S(E^{\ast})$ for any convex subset
of $S(E)$ (which could be empty!).

\begin{prop}
\label{convpn}
A convex subset $C$ has empty interior if and only if its dual $C^{\ast}$ is not proper.
If $C$ is open and proper, then the same is true for $C^{\ast}$,
and $C^{\ast\ast}=C$.\fin
\end{prop}

Recall that a support hyperplane to an open convex subset $C$ is a projective hyperplane meeting
the closure of $C$ but not $C$ itself.

\begin{prop}
\label{dualsupport}
Let $C$ be a proper open convex subset of $S(E)$. The support hyperplanes of
$C$ are the projections in $S(E)$ of the boundary points of $C^{\ast}$. More precisely, 
if $[\mbox{Ker} \alpha]$ is a support hyperplane of
$C$ at $[x] \in \partial C$, then $[\mbox{Ker} x^{\ast}]$ is a support hyperplane of
$C^{\ast}$ at $[\alpha] \in \partial C^{\ast}$.\fin
\end{prop}

\section{Spacelike and nontimelike surfaces}
\label{SPA}

The notion of nontimelike hypersurfaces in $\mbox{Ein}_{n} \approx {\mathbb S}^{n-1} \times {\mathbb S}^{1}$ or $\mbox{AdS} \approx {\mathbb D}^2 \times {\mathbb R}$ can be easily extended to the nonsmooth case: define it as closed subsets which are locally the graphs of $1$-Lipschitz maps from ${\mathbb S}^{n-1}$ into ${\mathbb S}^{1}$ or from
 ${\mathbb D}^{2}$ into ${\mathbb S}^{1}$. If moreover the Lipschitz maps are contracting,  i.e. are functions $f$ such that the equality $\vert f(x)-f(y)\vert = d(x,y)$ implies $x=y$, then 
the nontimelike surface is said \emph{spacelike}. The same notions apply in the coverings $\widehat{\mbox{Ein}}_n$ and $\widetilde{\mbox{AdS}}$.

Since it is Lipschitz, $f$ is differentiable almost everywhere. 
For any $C^1$-curve $c: [0,a] \rightarrow {\mathbb D}^2$ in 
the domain of definition of $f$, let $l(c)$ be the integral over $[0, a]$ 
of the square root of the AdS-norm of $(c'(t), D_{c(t)}f(c'(t)))$. 
Define then the distance between $(x, f(x))$ and $(y, f(y))$ as the 
infimum of the $l(c)$. This procedure 
endows the spacelike surface $S$ with a distance. Of course, 
this construction applies more generally to spacelike hypersurfaces 
in any lorentzian space. For more details, see \cite{beem}.

Achronal (resp.  strictly achronal) hypersurfaces are nontimelike (resp. spacelike) hypersurfaces, but the converse is not true. The main goal of this {\S} is to discuss under which additionnal hypothesis a nontimelike hypersurface is achronal.

\subsection{The redshift phenomenom}
\label{sub.redshift}

Let $M$ be a lorentzian manifold. For any timelike tangent vector $v$ at a point $x$ of $M$, the orthogonal projection in $T_xM$ on the orthogonal hyperplane $v^\perp$ increases the length of spacelike vectors: this fact is at the origin of the so-called ``redshift phenomenom''. It implies the well-known ``twins paradox''. There is another power\-full consequence, already observed in \cite{mess} (see also \cite{harris1, harris2}).

Assume the existence of one-parameter group $\Phi^{t}$
of isometries of $M$ such that the orbits of $\Phi^t$ are the fibers of a fibration $\pi: M \rightarrow Q$.
The base space $Q$ can be equipped with a riemannian metric as follows: for any point $x$ of $Q$ and any tangent
vector $v$ of $Q$ at $x$, define the norm of $v$ as the norm in $M$ of any vector orthogonal to the fiber $\pi^{-1}(x)$ and projecting on $v$ by the differential $d\pi$. This is well-defined since
the $\Phi^{t}$ are isometries.

Recall that we equipped spacelike surfaces with a distance function:

\begin{lem}
\label{langevin}
Let $\varphi: S \rightarrow M$ be an isometric immersion of a riemannian manifold such
that $\varphi(S)$ is an immersed locally acausal hypersurface.
The composition $\pi \circ f: S \rightarrow Q$ is distance increasing.
\end{lem}

\preu
When $\varphi$ is $C^1$, the lemma follows from the  observation above. 
The general case is a limit case. Details are left to the reader.
\fin

When $M$ is the anti-de Sitter space AdS, we can take as one-parameter subgroup the subgroup $SO(2)$ of $\mbox{SO}_0(2,2)$ acting in $E$
on the $(u,v)$ coordinates and fixing the coordinates $x_{1}, x_{2}$ (recall that we actually used this subgroup to define the time-orientation of AdS). The quotient space
of this timelike action equipped as above with a riemannian metric is isometric to the
hyperbolic space : therefore,
spacelike hypersurfaces in anti-de Sitter space correspond to distance increasing maps into 
${\mathbb H}^{2}$. 

\begin{prop}[Lemma $6$ in \cite{mess}]
\label{pro.graphe}
Let $S$ be a complete riemannian surface and let 
$\varphi: S \rightarrow \mbox{AdS}$ be an isometric immersion. Then, $\varphi$ is an embedding and $\varphi(S)$ is the graph of some contracting map 
${\mathbb H}^{2} \rightarrow {\mathbb S}^{1}$.\fin
\end{prop}

Every timelike geodesic in AdS is the fiber of some fibration $\pi$ as above.
It follows that under the hypothesis of Proposition~\ref{pro.graphe}, 
$\varphi(S)$ meets every timelike
geodesic in one and only one point.

\subsection{Proper nontimelike surfaces}
The conformal coordinates enable to extend the proposition~\ref{pro.graphe} to the non-timelike case. 
By \emph{proper} nontimelike surface $f: S \rightarrow \mbox{AdS}$, we mean an immersion such 
that the immersion $f$ is proper, i.e. that the preimage of a compact domain is compact.

\begin{prop}
\label{pro.graphe2}
Let $\varphi: S \rightarrow \mbox{AdS}$ be a proper nontimelike surface without boundary.
Then, $\varphi$ is an embedding and $\varphi(S)$ is the graph of some $1$-Lipschitz map 
${\mathbb D}^{2} \rightarrow {\mathbb S}^{1}$.
\end{prop}

\preu
The projection $\pi: \mbox{Ein}_{3} \rightarrow {\mathbb S}^{2}$ induces a
projection $\pi_{a}: {\mathcal A} \rightarrow {\mathbb D}^{2}$ (cf. remark \ref{scan}). We claim that
$\pi_{a} \circ \varphi$ is covering map. We will justify it by proving that it has the path lifting
property. For this purpose, it is enough to prove it for paths in the open hemisphere
${\mathbb D}^{2}$ which are
segments of geodesics of ${\mathbb D}^{2} \subset {\mathbb S}^{2}$: let $[a,b]$ be such a geodesic segment and such that
$a=\pi(\varphi(x))$ for some element $x$ of $S$.  Since $\pi \circ \varphi$ is open, there is a section $\sigma$ of $\pi \circ \varphi$ 
defined over a subinterval $[a,c[ \subset [a,b]$. The point is that 
$p \circ \varphi \circ \sigma: [a,c[ \subset [a,b] \rightarrow {\mathbb S}^{1}$ is then a $1$-Lipschitz map (where $p$
is the projection on the second factor). Therefore, it can be continuously extended over $c$. The properness
of $\varphi$ then implies that $\sigma$ can be extended over $c$. Hence, $c=b$ and $\sigma$ can be extended over all
$[a,b]$: $\pi \circ \varphi$ has the path lifting property.

We thus know that $\pi \circ \varphi$ is a covering map over ${\mathbb D}^{2}$, hence a homeomorphism. It follows
that $\varphi(S)$ is the graph of a $1$-Lipschitz map ${\mathbb D}^{2} \rightarrow {\mathbb S}^{1}$ and that
$\varphi$ is an embedding.\fin

From now, we will assume that $S$ fits {\em inside\/} AdS, i.e. that
$\varphi: S \rightarrow \mbox{AdS}$ is an inclusion map.
The additional advantage of our point of view is that the proof of lemma $7$ of \cite{mess},
which was a delicate matter in this paper, now appears as completely obvious:

\begin{prop}
\label{facile!}
A proper nontimelike hypersurface $S$ in AdS extends continuously in 
${\mathcal A} \cup \partial{\mathcal A} \subset \mbox{Ein}_{3}$ as a closed topological disk,
whose boundary $\partial S$ is a topological nontimelike circle in $\partial{\mathcal A} \approx \mbox{Ein}_{2}$.
\end{prop}

\preu
Any $1$-Lipschitz map from ${\mathbb D}^{2}$ into ${\mathbb S}^{1}$ extends continuously as
a $1$-Lipschitz on the closure of ${\mathbb D}^{2}$ in ${\mathbb S}^{2}$.\fin

Consider now the universal covering $\widetilde{\mbox{AdS}}$. Select any connected component $\tilde{S}$ of
the preimage of $S$ in $\widetilde{\mbox{AdS}}$ by the covering map. Clearly, $\tilde{S}$
is the graph of a $1$-Lipschitz maps from
${\mathbb D}^{2}$ into ${\mathbb R}$. Moreover:

\begin{cor}
\label{p1}
A proper nontimelike surface in $\widetilde{\mbox{AdS}}$ meets every timelike geodesic in one and only one 
point.\fin
\end{cor}

\begin{cor} 
Proper nontimelike surfaces in AdS are achronal subsets of AdS.
\fin
\end{cor}

Observe also:

\begin{lem}
\label{adhachronal}
Let $S$ be a proper spacelike hypersurface in $\widetilde{S}$. Then, the closure 
$\overline{S}$ of $S$ in $\widetilde{\mbox{AdS}} \cup \partial\widetilde{\mbox{AdS}}$ 
is not necessarily strictly achronal, but if two points in this closure are causally related, then they both belong to $\partial\widetilde{\mbox{AdS}}$. In particular, no point of $S$ is causally related to a point of $\overline{S} \cap \partial\widetilde{\mbox{AdS}}$.
\end{lem}

\sketch
The closure $\overline{S}$ is the graph of a $1$-Lipschitz function $f: \overline{\mathbb D}^2
\rightarrow \mathbb R$ such that the restriction of $f$ to ${\mathbb D}^2$ is contracting.
If a point $(x, f(x))$ of $S$ is causally related to a point $(y, f(y))$ of $\overline{S} \cap 
\partial\widetilde{\mbox{AdS}}$, and if $\mid f(y) - f(x) \mid = d(x,y)$, then 
for any $z$ on the geodesic $[x,y]$ in $\overline{\mathbb D}^2$ the equality 
$\mid f(z) - f(x) \mid = d(x,z)$ must hold. This is a contradiction since $S$ is assumed
acausal.\fin

Pure lightlike nontimelike surfaces in ${\mathbb A}{\mathbb D}{\mathbb S} \approx \mbox{AdS}$ are easy to 
describe:

\begin{cor}
Pure lightlike surfaces in ${\mathbb A}{\mathbb D}{\mathbb S}$ are connected components of intersections 
between ${\mathbb A}{\mathbb D}{\mathbb S}$ and projections in $S(E)$ of lightlike hyperplanes of $E$.\fin
\end{cor}

\subsection{Embeddings in ${\mathbb H}^2 \times {\mathbb H}^2$}
\label{sub.isomh2}

Recall the projectivized timelike tangent bundle $\mathcal T$ (cf. definition~\ref{def.timebundle}); 
more precisely, the connected component ${\mathcal T}^+$ corresponding to future oriented timelike vectors.

\begin{defin}
\label{def.gauss}
The \emph{Gauss flow\/} is the flow $G^t$ on ${\mathcal T}^+$ defined by:

\[ G^t(x,y ) = (x\cos{t}+ y\sin{t}, -x\sin{t} + y\cos{t}) \]

\end{defin}

This flow commutes with the $\mbox{O}_0(2,2)$-action. Moreover, it is easy to check that $G^t$  is isometric. The Killing vector field $Z$ generating $G^t$ is easy to describe (see \S~\ref{sub.timebundle} for the convention on tangent vectors to ${\mathcal T}^+$):

\[ Z(x,y) = (y, -x) \]

Let $Q_G$ be the orbit space of $G^t$ and let $\pi_G: {\mathcal T}^+ \rightarrow Q_G$ be the quotient map. We equip $Q_G$ with a riemannian metric as discussed in 
\S~\ref{sub.redshift}: the norm of a tangent vector $\zeta$ to $Q_G$ is the norm of any tangent 
$w$ to ${\mathcal T}^+$ orthogonal to $Z$ and such that $d\pi_G(w) = \zeta$. 

\rque
\label{rk.+precis}
Let's be more precise: let $w = (u,v)$ be a tangent vector at $(x,y) \in {\mathcal T}^+$. The norm of $\zeta = d\pi_G(w)$ is the norm of $w+\lambda{Z}$, where $\lambda$ is the unique real number such that $w+\lambda{Z}$ is orthogonal to $Z$. A straightforward computation shows $\lambda = \langle y \mid u \rangle = -\langle x \mid v \rangle$. Hence, the norm of $\zeta$ is $\frac{1}{4}(\mid u \mid^2 + \mid v \mid^2) + \frac{1}{2}\langle x \mid v \rangle^2$.
\erque

\begin{prop}
\label{pro.QG}
The riemannian orbit space $Q_G$ is homothetic to the riemannian product ${\mathbb H}^2 \times {\mathbb H}^2$ of two copies of the hyperbolic plane.
\end{prop}

\preu
Identify AdS with $G = \mbox{SL}(2,{\mathbb R})$ and consider the upper-half plane model
of ${\mathbb H}^2$. Let $i$ denote the point $\sqrt{-1}$ in ${\mathbb H}^2$. Let $x_0$ be the identity 
matrice and $y_0$ be the matrice representing the rotation by angle $\pi/2$ around $i$.
Observe that under the indentifications above, $(x_0,y_0)$ belongs to ${\mathcal T}^+ \subset \mbox{SL}(2,{\mathbb R}) \times \mbox{SL}(2,{\mathbb R})$. 
Define a $G \times G$-equivariant map $F: {\mathcal T}^+ \rightarrow {\mathbb H}^2 \times {\mathbb H}^2$ as follows:
if $(x, y) = (g_Lx_0g^{-1}_R, g_Ly_0g_R^{-1})$, then $F(x,y) = (g_Li,g_Ri )$.
Observe that it is well-defined: indeed, if $(g_L, g_R)$ fixes $(x_0, y_0)$, then $g_L = g_R$ commutes with $y_0$: $g_L = g_R$ preserves $i$ in ${\mathbb H}^2$. Moreover, the preimage of $(i, i)$ is precisely the $G^t$-orbit of $(x_0,y_0)$: it follows that $F$ induces a homeomorphism between $Q_G$ and ${\mathbb H}^2 \times {\mathbb H}^2$.
The only remaining point to check is that $F$ is a homothety. 
Since it is equivariant, we just have to consider the differential of $f$ at $(x_0, y_0)$. 
The computation can be performed as follows:
let $A$ be an element of the Lie algebra $\mathcal G$:

\[ A = \left(\begin{array}{cc}
    \alpha & \beta \\
    \gamma & -\alpha
\end{array}\right) \]

Using lemma~\ref{rk.+precis}, we obtain that the norm in ${\mathcal T}^+$ of the tangent vector to $x_0$ of the curve $t \mapsto (\exp(tA)x_0, \exp(tA)y_0)$ is $\frac{\alpha^2}{2} + \frac{1}{8}(\beta+\gamma)^2$, whereas the norm in ${\mathbb H}^2 \times {\mathbb H}^2$ of the tangent vector at $(i,i)$ of the image curve $t \mapsto (\exp(tA)i, i)$ is $(\beta+\gamma)^2+4\alpha^2$. It follows that the restriction of $F$ to every left $G_L$-orbit is an isometry on the image ${\mathbb H}^2 \times \{ \ast \}$ with the metric divided by $2\sqrt{2}$.

A similar calculus holds for curves $t \mapsto (x_0\exp(-tA), y_0\exp(-tA))$, proving that the restriction of $F$ to every $G_R$-orbit is an homothety of factor $8^{-1/2}$ on $\{ \ast \}\times {\mathbb H}^2$. The proposition follows from the fact that in ${\mathcal T}^+$, $G_L$-orbits are orthogonal to $G_R$-orbits and that every ${\mathbb H}^2 \times \{ \ast \}$ are orthogonal to every $\{ \ast \} \times {\mathbb H}^2$.
\fin

Consider a $C^1$ embedded spacelike surface $S \subset \mbox{AdS}$. 
For any $x$ in $S$, let $n(x)$ be the unique future-oriented unit timelike 
vector normal to $S$ at $x$: $(x, n(x))$ is an element of 
$PT^+_{-1}\mbox{AdS} \approx {\mathcal T}^+$. In other words, 
the embedding of $S$ in AdS lifts to an embedding 
$n: S \rightarrow {\mathcal T}^+$.

\begin{defin}
\label{def.gaussmap}
The Gauss map of $S$ is $n: S \rightarrow {\mathcal T}^+$.
\end{defin}

\begin{lem}
\label{lem.isom}
The image of the Gauss map is a topological spacelike surface in ${\mathcal T}^+$.
Moreover, the restriction of $\pi_G$ to the image of the Gauss map 
endowed with the induced metric is an isometry.
\end{lem}

\preu
Consider first the case where $f$ is $C^2$. Then, the Gauss map is $C^1$. 
Let $(u,v)$ be a tangent vector to the image of $n$. By definition of the 
Gauss map, the tangent vector $u$ to $S$ satisfies: 
$\langle y \mid u \rangle = 0$. Hence, 
$\langle x \mid v \rangle = -\langle y \mid u \rangle = 0$: 
$u$ and $v$ both belong to the spacelike $2$-plane $x^\perp \cap y^\perp$: the sum of their norms is positive.
These identities mean also that $(u,v)$ is orthogonal to the Killing vector field $Z(x,y)$. 
It follows that the restriction of $\pi_G$ to the image is an isometry.
The $C^1$-case is a limit case: any $C^1$ spacelike surface can be $C^1$-approximated by a $C^2$ spacelike 
surface. It follows that even in the $C^1$-case, the image of the Gauss map is locally achronal. 
Observe that a locally achronal surface which is not locally acausal must
contain a lightlike geodesic segment: we can apply the argument above, leading to a contradiction. 
It follows that the image of the Gauss map is a topological spacelike surface.

The distance between two points in the image of the Gauss map is the supremum of the length
of Lipschitz curves joining the two points. Let $c: t \mapsto (x(t), y(t))$ be such a 
Lipschitz curve.
It is differentiable almost everywhere, with derivative $(x'(t), y'(t))$. 
The tangent vector $x'(t)$, where it is defined, is orthogonal to $y(t)$. 
The derivation of the identity $\langle x(t) \mid y(t) \rangle = 0$ then implies
that $\langle x(t) \mid y'(t) \rangle = 0$. It follows that $(x'(t), y'(t))$
is orthogonal to $Z(x(t), y(t))$ almost everywhere. Hence the length of $c$
is equal to the length of its projection in $Q_G$. The restriction of $\pi_G$ 
to the image of the Gauss map is thus
an isometry.
\fin

\rque
The Gauss map $n: S \rightarrow {\mathcal T}^+$ in general is not isometric! Actually, the 
metric along the image of $n$ involves the second fundamental form of $S$. Observe that 
$n$ is isometric if and only if $S$ is totally geodesic and that $n$ is conformal if 
and only if $S$ is totally umbilic.
\erque

\rque
\label{rk.pkgauss}
Our choice of terminology is justified by the following observation: if $S_t$ is the image of 
$S$ under the Gauss flow in the usual meaning, then $n(S_t)$ is the image of $n(S)$ by the 
Gauss flow $G^t$ we have defined above.
\erque

This observation extends to a much less regular situation: the case where $S$ is maybe 
non $C^1$, but \emph{convex.\/} This notion is meaningful, due to the local real projective 
structure
of $\mbox{AdS} \approx {\mathbb A}{\mathbb D}{\mathbb S}$:

\begin{defin}
\label{def.surfconvexe}
An embedded topological spacelike surface $S$ in $AdS$ is future-convex 
if any point $x$ in $S$ admits a geodesically convex neighborhood $U$ in $AdS$ such that for any 
$y$ in $S \cap U$, the geodesic segment $[x,y]$ is contained in the causal future of 
$S \cap U$ relatively to $U$.
\end{defin}

\rque
Observe that geodesic segments in $U$ are restrictions of projective lines of 
${\mathbb A}{\mathbb D}{\mathbb S}$. 
With the local description of spacelike subsets as graphs of
functions, it follows that $S$ is future-convex if and only if it is convex in 
${\mathbb A}{\mathbb D}{\mathbb S} \subset S(E)$ in the usual meaning and that 
it is ``curved in the future direction''.
\erque

\begin{defin}
\label{def.gaussgraph}
Let $S$ be a future-convex spacelike surface in AdS. The Gauss graph of $S$ is the set of pairs 
$(x,y) \in {\mathcal T}^+$ such that:

-- $x$ belongs to $S$,

-- there is a neighborhood $U$ of $x$ in AdS such that for every for any $x'$ in $U \cap S$ 
the scalar product $\langle x' \mid y \rangle$ is nonpositive. In other words, 
the connected component $P$ of $y^\ast$ containing $x$ is a support hyperplane of $S$, 
such that $S \cap U$ is contained in the causal future of $P \cap U$.
\end{defin}

\rque
When $S$ is $C^1$, the Gauss graph is the graph of the Gauss map as defined in the definition~\ref{def.gaussmap}.
\erque
In the sequel, we assume the reader acquainted with the familiar notion of convex surfaces in $P(E)$.

\begin{prop}
The Gauss graph of a future-convex spacelike surface is a locally embedded topological spacelike surface in ${\mathcal T}^+$.
\end{prop}

\sketch
Let $(x_0,y_0)$ be an element of the Gauss graph ${\mathcal N}(S)$ of a future-convex spacelike surface $S$.
There is a neighborhood $U$ of $x_0$ such that 
$S \cap U$ is contained in the boundary of a proper compact convex subset $C$ of $S(E)$. 
Then, ${\mathcal N}(S)$ is contained in the set $\mathcal C$ of pairs $(x,y)$, where $x$ belongs to $\partial C$, $y^\flat$ belongs to the boundary of the dual convex $C^\ast \subset S(E^\ast)$ and $y^\flat$ is a support hyperplane to $C$ at $x$. The set $\mathcal C$ is notoriously a topological surface and it should be clear to the reader that ${\mathcal N}(S)$ contains a neighborhood of $(x_0, y_0)$ in $\mathcal C$. It follows that ${\mathcal N}(S)$ is a topological surface near $(x_0, y_0)$.

The local achronality of ${\mathcal N}(S)$ follows from the local achronality of $S$ and the fact that
if $y$ is a point near $y_0$ in the future of $y_0$, then $\langle x_0 \mid y \rangle$ is positive,
which is a contradiction since $\langle y \mid x \rangle$ should be nonpositive for every $x$ in $S$.
\fin

Hence, the Gauss graph admits a natural distance (recall the definition in the beginning of this {\S}). Once more, we consider the general case as a ``limit case''of the regular one, leaving to the reader the proof of the following lemma:

\begin{lem}
\label{lem.isomconvex}
The restriction of $\pi_G$ to ${\mathcal N}(S)$ is an isometry.
\fin
\end{lem}

The distance in ${\mathcal N}(S)$ is evaluated by the computation of the ``length'' of Lipschitz curves. We will also need the following fact:

\begin{lem}
\label{lem.uvpositif}
Let $c: [0,a] \rightarrow {\mathcal N}(S)$ be a Lipschitz curve. It is differentiable almost everywhere. Any tangent vector $c'(t) = (u,v)$ satisfies:
\[ \langle u \mid v \rangle \geq 0 \]
\end{lem}

\sketch
Once more, we only consider the regular case: we assume that $S$ is $C^2$ and that the curve $(x(t), y(t))$ has tangent vectors $(u(t), v(t)) = (x'(t), y'(t))$. The derivation of $\langle y \mid x' \rangle =0$ implies $\langle x' \mid y' \rangle = - \langle y \mid x'' \rangle$. Since $S$ is future oriented,
the second derivative $x''(t)$ must point towards the future of $S$, hence, in the future of
the support hyperplane at $x$ contained in $y^\ast$. The negativity of $\langle y \mid x'' \rangle$ follows.

The general case is similar, once observed that convex surfaces are $C^2$ almost everywhere.
\fin

\section{Cauchy developments}
\label{sec.cauchy}

In this {\S}, we study Cauchy development in AdS of achronal subsets
of AdS. It turns out that Cauchy developments  can be defined as \emph{invisible domains\/} 
from achronal subsets of $\partial AdS$. We start with the most familiar notion of Cauchy development of spacelike surfaces (domain of dependance in \cite{mess}) and then to extend to the most general context.

\subsection{Cauchy developments of spacelike surfaces}

We consider in this section a proper spacelike hypersurface $\tilde{S}$ in $\widetilde{\mbox{AdS}}$.  Actually,  all the results apply if $\tilde{S}$ is more generally any strictly achronal surface. 

According to proposition \ref{facile!}, the boundary $\partial\tilde{S}$ of $\tilde{S}$ 
in $\partial\widetilde{\mbox{AdS}} \approx \widehat{\mbox{Ein}}_{2}$ is an achronal topological circle. We denote by $S$, $\partial S$  the projections in AdS, $\mbox{Ein}_2$.

\begin{lem}
\label{etendre}
The past development $P(\tilde{S})$ is the set of points $x$ in ${\widetilde{\mbox{AdS}}}$ such that every lightlike geodesic containing $x$ meets $\tilde{S}$ in its future.
\end{lem}

Of course, the analogous property for $F(\tilde{S})$ is true. 

\preud{etendre}
Assume that every lightlike geodesic containing $x$ meets $\tilde{S}$ in its future.
According to corollary~\ref{p1} every timelike geodesic containing $x$ intersect
$\tilde{S}$ in one and only one point. The set of timelike geodesics containing $x$
coincide with the space of future oriented timelike elements of 
$T_{x}\widetilde{\mbox{AdS}}$, i.e. a copy of ${\mathbb H}^2$. In particular,
it is connected, and the boundary of this space is the space of lightlike geodesics containing
$x$. It follows that every timelike geodesic containing $x$ meets $\tilde{S}$ in the future
of $x$. The union of the nonspacelike geodesic segments $[x,y]$ with $y$ in
$\tilde{S}$ is a topological disk $B'$. Any properly embedded causal path starting from $x$ cannot escape from $B'$ through $L$: it must therefore intersect $B' \subset \tilde{S}$.\fin

We pursue our investigation in the Klein model ${{\mathbb A}{\mathbb D}{\mathbb S}}$.
According to Lemma~\ref{Edansds2}, since it is generic, $\tilde{S}$ is contained in some affine domain. Therefore, it projects injectively in AdS as a strictly achronal surface $S$, contained in some affine domain.

\begin{defin}
We define $T(S)$ as the set of points $x$ in ${\mathbb A}{\mathbb D}{\mathbb S}$ such that the affine domain ${A}(x)$ contains $S$.
\end{defin}

\begin{lem}
\label{Touv}
$T(S)$ is a neighborhood of $S$.
\end{lem}

\preu
By definition, $T(S)$ is the set of elements $x$ of ${\mathbb A}{\mathbb D}{\mathbb S}$ such that $\langle x \mid  y \rangle$ is negative for every $y$ in $S$. In other words, in the terminology of \S~\ref{dual},
$T(S)$ is the intersection between ${\mathbb A}{\mathbb D}{\mathbb S}$ and the image under $\sharp$ of the dual of the convex hull $\mbox{Conv}(S)$ of $S$ in $S(E)$. 

Consider an element $x_{0}$ of $S$. It is the projection of some vector $v_0$ in 
$E$. Then, $S$ is the projection of some subset $S(v_0)$ in $P_0$ (see \S~\ref{sub.einachronal}). According to proposition~\ref{pro.accro}, since $S$ is achronal, 
for every $x$, $y$ in $S(v_0)$, the scalar product 
$\langle x \mid y \rangle$ is negative.
It follows that $T(S)$ contains $S$. 

Moreover, an element $x$ of $S(v_0)$ does not project to a point in the interior of $T(S)$ if and only if
there is a sequence of points $x_n$ in $S(v_0)$ such that $\langle x \mid x_n \rangle$ tends to $0$. Up to some subsequence, the projections in $x_n$ in ${\mathbb A}{\mathbb D}{\mathbb S}$ converge to some element $\bar{x}$ of ${\mathbb A}{\mathbb D}{\mathbb S} \cup \partial{\mathbb A}{\mathbb D}{\mathbb S}$. Then, $\bar{x}$ would be a point in the closure of $S$ in 
$\mbox{AdS} \cup \partial\mbox{AdS}$ causally related to the point $x$ in $S$. It contradicts lemma~\ref{adhachronal}.\fin

Let $T_0(S)$ be the interior of $T(S)$. It contains $S$. 
Select any element $x_0$ in $T_0(S)$.
Then, $S$ is contained in the affine domain $A_0 = A(x_0)$. Actually, the fact that $x_0$ belongs to the interior of $T(S)$ means that $S$ is contained in a compact domain of
the affine patch $V_0 = V(x_0)$. Hence, the closure $\overline{S}$ in $\mbox{AdS} \cup \partial\mbox{AdS}$ is a closed topological disc in $V(x_0)$, with boundary  $\partial{S}$ contained in the one-sheet hyperboloid $\mbox{Ein}_2 \cap V_0$.

\begin{prop}
\label{T=C}
The restriction of $p$ to the Cauchy development ${\mathcal C}(\tilde{S})$ is injective, with image $T_0(S)$.
\end{prop}

\preu
First observe that $T_0(S)$ is contained in every affine domain $A = A(x)$, with $x$ in $S$. 
Let $\widetilde{A}$ be the affine domain in $\widetilde{\mbox{AdS}}$ containing $\tilde{S}$ such that $p(\widetilde{A}) = A$.

Let $\tilde{x}$ be an element of $P(\tilde{S})$. Consider a conformal parametrization
of $\widetilde{\mbox{AdS}}$ by ${\mathbb D}^2 \times {\mathbb R}$ such that 
$\tilde{x}$ has coordinates $(x_0, 0)$, where $x_0$ is the north pole, i.e. is the unique point of
${\mathbb D}^2$ at distance $\pi/2$ of $\partial{\mathbb D}^2$. In these coordinates, $\tilde{S}$ is the graph of some contracting function $f: {\mathbb D}^2 \rightarrow {\mathbb R}$. Since $\tilde{x}$ is in the past of $\tilde{S}$, we have $f(x_0) > 0$. 

Every future oriented lightlike geodesic starting from $\tilde{x}$ intersect $\tilde{S}$: it follows that there is an open topological disc $B$ in ${\mathbb D}^2$ containing $x_0$ and such that $f(x) = d(x_0, x)$ for every $x$ in $\partial B$. Since $f$ is contracting and since any point in $\overline{\mathbb D}^2$ is at distance strictly less than $\pi/2$ of some point in $\partial B$,
it follows that the extension $\bar{f}$ of $f$ over $\overline{\mathbb D}^2$ takes value in $]-\pi/2, \pi/2 [$. Since $\overline{\mathbb D}^2$ is compact, $\bar{f}$ takes actually value in some closed intervall $[-\pi/2 + \epsilon, \pi/2 - \epsilon]$. It follows that $p(\tilde{x})$ belongs to $T_0(S)$. In other words, the image $p(P(\tilde{S}))$ is contained in $T_0(S)$. Hence,
since $P(\tilde{S})$ is connected and since $A$ contains $T_0(S)$, $P(\tilde{S})$ is contained in $\widetilde{A}$. 

Applying a similar argument to $F(\tilde{S})$, we obtain that $C(\tilde{S})$ is contained in $\widetilde{A}$ and that $p(C(\tilde{S})) \subset T_0(S)$. 

Assume now that $\tilde{x}$ is an element of $\widetilde{A}$ such that $p(\tilde{x})$ belongs to $T_0(S)$ and select once more a conformal parametrization
of $\widetilde{\mbox{AdS}}$ by ${\mathbb D}^2 \times {\mathbb R}$ such that 
$\tilde{x}$ has coordinates $(x_0, 0)$, where $x_0$ is the north pole.  The affine domain 
associated to $\tilde{x}$ is then the open domain $\{ (y, \theta) / \mid \theta \mid < \pi/2 \}$.
Hence, since $\tilde{S}$ has to be contained in this affine domain, the map $f: {\mathbb D}^2 \rightarrow {\mathbb R}$ admitting $\tilde{S}$ as graph takes value in $]-\pi/2, \pi/2 [$.
More precisely, since $p(\tilde{x})$  belongs to the interior of $T(S)$, $f$ takes value in some
interval $[-\pi/2 + \epsilon, \pi/2 - \epsilon]$. 

If $f(x_0) = 0$, $\tilde{x}$ belongs to $\tilde{S} \subset C(\tilde{S})$. Assume $f(\tilde{x}_0) > 0$, i.e. assume that $\tilde{x}$ is in the past of $\tilde{S}$. Define $g(y) = f(y) - d(y, x_0)$. This function is continuous, positive on $x_0$ and negative near $\partial{\mathbb D}^2$. Hence, any geodesic ray in ${\mathbb D}^2$ starting from $x_0$ admits some point where $g$ is $0$. It means that any future oriented lightlike geodesic starting form $\tilde{x}$ intersect $\tilde{S}$: $\tilde{x}$ belongs to $P(\tilde{S})$.

A similar argument proves that in the remaining case, i.e. when $\tilde{x}$ belongs to the future of $\tilde{S}$, then it belongs to $F(\tilde{S})$. The proposition follows.\fin

\begin{cor}
$T_0(S)$ is globally hyperbolic, with Cauchy surface ${S}$.
\fin
\end{cor}

\begin{cor}
\label{dS}
The Cauchy development ${\mathcal C}(\tilde{S})$ is contained in a de Sitter domain.
\fin
\end{cor}

\subsection{The Cauchy development as the invisible domain from $\partial S$}

\begin{defin}
For any point $x$ of $\mbox{Ein}_{2}$, let $T_{x}$ 
be the projective hypersurface in $S(E)$ containing $x$ and tangent to 
$\mbox{Ein}_{2}$. For any pair of points $x$, $y$ 
of $\mbox{Ein}_{2}$, let $E_{xy}$
be the open half-space in $S(E)$ bounded by $T_{x}$ and $T_{y}$ 
and containing the segment $]x,y[$ contained in ${\mathbb A}{\mathbb D}{\mathbb S}$.
\end{defin}

\begin{defin}
\label{defabs1}
Let $E(\partial{S})$ be the intersection of all $E_{xy}$ when $(x,y)$ describes 
$\partial{S} \times \partial{S}$ minus the diagonal.
\end{defin}

\rque
\label{defabs}
The presentation above is suitable for getting some geometrical vision 
of $E(\partial{S})$. It is more relevant for the proofs below to
consider the following equivalent definition: 
recall that $\partial S$ is the projection in ${\mathbb A}{\mathbb D}{\mathbb S}$ 
of a compact subset ${\partial S}(v_0)$ in $\{ Q=0 \} \cap P(v_0)$, where $P(v_0)$ is some affine hyperplane in $E$. Then, $E(\partial S)$ is the projection of the set of elements $v$ of $E$ satisfying $\langle v \mid x \rangle < 0$ for every $x$ in ${\partial S}(v_0)$. 
In other words:
{\em $E(\partial{S})$ is the dual of the convex hull in $S(E^{\ast})$ of
$\partial{S}^{\flat}$.\/} We write:
$E(\partial{S}) = \mbox{Conv}(\partial{S}^{\flat})^{\ast}$.
Alternatively: $E(\partial{S}) = (\mbox{Conv}(\partial{S})^{\ast})^{\sharp}$.
Observe that the compactness of ${\partial S}(v_0)$ implies that $E(\partial S)$ is open.
\erque

\rque
\label{deff} 
Using the conformal model $\widetilde{\mbox{AdS}} \approx {\mathbb D}^2 \times {\mathbb R}$, we obtain another equivalent definition: let $f:\partial{\mathbb D}^{2} \rightarrow {\mathbb R}$ be the $1$-Lipschitz map admitting as graph the topological circle $\partial\tilde{S}$.
Define $f_{\pm}: {\mathbb D}^{2} \rightarrow {\mathbb R}$ by:

- $f_{-}(x)$ is the supremum of $f(y)-d(x,y)$ when $y$ describes ${\mathbb S}^{1}$,

- $f_{+}(x)$ is the infimum of $f(y)+d(x,y)$ when $y$ describes ${\mathbb S}^{1}$.

Then, it is easy to check that $f_{\pm}$ are both $1$-Lipschitz extensions of $f$. Moreover, it follows easily
from the arguments used in the proof of proposition~\ref{T=C} that 
$E(\partial{S}) \cap {\mathbb A}{\mathbb D}{\mathbb S}$ is the 
projection of the domain 
$E(\partial\tilde{S}) = \{ (x,\theta) \in {\mathbb D}^{2} \times {\mathbb R} / f_{-}(x) < \theta < f_{+}(x) \}$. The advantage of our definition \ref{defabs1} is
to explicit the {\em convex\/} character of ${E}(\partial{S})$.
\erque

\rque
\label{definvisible}
The definition in remark~\ref{deff} can be interpreted in the following way (recall \S~\ref{adsdads}): \emph{$E(\partial S)$ is the set of points of AdS which are not causally related to any element of $\partial S$.\/} Hence, it is appropriate to consider $E(\partial S)$ as
the \emph{invisible domain from $\partial S$.\/}
\erque

\rque
The extensions $f_{\pm}$ above can be defined in any metric space $X$ for any $1$-Lipschitz map $f$ 
defined over a closed subset $Y \subset X$: they are still $1$-Lipschitz and any extension $F$ of $f$
satisfy $f_{-} \leq g \leq f_{+}$. In general, $f_{-}$ and $f_{+}$ can coincide on some closed
subset of $X \setminus Y$: the points belonging to some minimizing geodesic joining two 
points $x$, $y$ of $Y$ such that $\vert f(x)-f(y)\vert = d(x,y)$. 
In the case where $X$ is the closed unit disc in the euclidean plane and $Y$ its boundary,
the set $f_{+}=f_{-}$ is a lamination. 

In the particular case we consider here, where $X$ is a hemisphere
and $Y$ its boundary, $Y$ is totally geodesic and extremities in $Y$ of minimizing geodesics
are opposite points of the boundary sphere. Moreover, minimizing geodesics joining two given
points form a foliation of the hemisphere. We recover easily from this that when $f$ correspond
to the generic topological sphere $\partial S$, then the associated $E(\partial S)\cap {\mathbb A}{\mathbb D}{\mathbb S}$ is open, i.e.,
$f_{-}<f_{+}$. 
\erque

The definition in remark~\ref{deff} implies easily:

\begin{lem}
\label{convfer}
The convex $E(\partial{S})$ is contained in ${\mathbb A}{\mathbb D}{\mathbb S}$. 
The intersection between its closure and $\mbox{Ein}_{2}$ is $\partial{S}$.
\end{lem}

\preu
Indeed, since $f_- = f_+$ over $\partial{\mathbb D}^2$, the intersection
between $\partial{\mathbb D}^2 \times {\mathbb R}$ and the closure 
$\{ (x,\theta) \in {\mathbb D}^{2} \times {\mathbb R} / f_{-}(x) \leq \theta \leq f_{+}(x) \}$ of
$E(\partial\tilde{S})$ in ${\mathbb D}^2 \times {\mathbb R}$ is simply the graph of $f$. The lemma follows easily,
since $E(\partial S)$ is convex.\fin 

In the same spirit:

\begin{lem}
Any proper nontimelike topological hypersurface contained in ${\mathbb A}{\mathbb D}{\mathbb S}$ and containing $x$, $y$
in its closure is necessarily contained in the closure of $E_{xy}$.\fin
\end{lem}

$E(\partial{S})$ is thus a convex subset of $S(E)$ containing $S$:
in particular, it is not empty! Actually, using the definition in remark~\ref{defabs} and since $\partial S$ is in the closure of $S$, it is straightforward to show the inclusion
$T_0(S) \subset E(\partial S)$. The inverse inclusion is true:

\begin{prop}
\label{pro.E=T}
The invisible domain $E(\partial S)$ is equal to the Cauchy development $T_0(S)$.
\end{prop}

\preu
Let $x$ be a point in $\partial T_0(S) \cap E(\partial S)$. Select as usual a conformal parametrization of the affine domain $A(x)$ by ${\mathbb D}^2 \times ]-\pi/2, +\pi/2[$, such that $x$ has coordinates $(x_0, 0)$, where $x_0$ is the North pole. The surface $S$ is the graph of a $1$-Lipschitz function $f: {\mathbb D}^2  \rightarrow \mathbb R$. 
Since $x$ belongs to $E(\partial S)$, the restriction of $f$ to $\partial{\mathbb D}^2$ takes value in $]-\pi/2, +\pi/2[$. On the other hand, since $x$ belongs to $\partial T_0(S)$, the map
$f$ takes value in $[-\pi/2, +\pi/2]$, but there is some element $y$ of $\overline{\mathbb D}^2$ such that $f(y) = \pm \pi/2$; let's say, $+\pi/2$. Observe that $y$ cannot belong to $\partial{\mathbb D}^2$. Hence, $f(x_0) \geq \pi/2 - d(x_0, y) > 0$: $x$ is in the past of $S$. Consider the function $g: z \mapsto f(z) - d(x_0, z)$: it is negative on $\partial{\mathbb D}^2$ and $g(x_0)$ is positive: it follows, as in the proof of proposition~\ref{T=C} that $x$ belongs to the past development $P(S)$.
Similarly, if $f(y) = -\pi/2$, we infer that $x$ belongs to the future development of $S$. 
This is a contradiction since $T_0(S)$ is the Cauchy development of $S$.\fin

According to corollary~\ref{dS}, $T_0(S) = E(\partial S)$ is contained in a de Sitter domain.
We can now say more:

\begin{lem}
\label{convpro}
If $\partial{S}$ is not a round circle, then the closure of
$E(\partial{S})$ is contained in a de Sitter domain.
\end{lem}

\preu
If $\partial S$ is not a round circle, then the interior of $\mbox{Conv}(\partial{S})$ is
not empty. Thus,
the same is true for $\mbox{Conv}(\partial{S}^{\flat})$. Points in the interior
of $\mbox{Conv}(\partial{S}^{\flat})$ correspond to an open set of flat spheres avoiding
$E(\partial{S})$.\fin

\subsection{Support hyperplanes}

\begin{lem}
\label{ECdual}
The boundary of $\mbox{Conv}(\partial{S})$ (resp. $E(\partial {S})$) in $S(E)$ is the
set of points dual to support hyperplanes to the closure of $E(\partial {S})$ 
(resp. $\mbox{Conv}(\partial {S})$) in $S(E)$.
\end{lem}

\preu
Corollary of remark \ref{defabs}, proposition \ref{dualsupport}, 
and lemma \ref{convfer}.\fin

We can be slightly more precise.
When $\partial{S}$ is not contained in a round circle, the complement of $\partial{S}$ in the boundary $\partial\mbox{Conv}(\partial{S})$ is contained in ${\mathbb A}{\mathbb D}{\mathbb S}$ and admits two connected components. 

\begin{defin} 
The future (respectively past) convex boundary $\partial^+C(\partial{S})$ 
(respectively $\partial^-C(\partial{S})$) is the connected component of $\partial\mbox{Conv}(\partial{S}) \setminus \partial{S}$ such that the interior of $\mbox{Conv}(\partial{S})$ is contained in the past (resp. future) of $\partial^+C(\partial{S})$ (resp. $\partial^-C(\partial{S})$).
\end{defin}

Observe that in the flat case, $\partial^+C(\partial{S}) = \partial^-C(\partial{S})$

Similarly, the complement of $\partial{S}$ in $\partial E(\partial{S})$ has two connected components:

\begin{defin} 
The future (respectively past) boundary $\partial^+E(\partial{S})$ (resp. $\partial^-E(\partial{S})$) is the connected component of $\partial{E}(\partial{S}) \setminus \partial{S}$ such that the interior of $E(\partial{S})$ is contained in the past (resp. future) of $\partial^+E(\partial{S})$ (resp. $\partial^-E(\partial{S})$).
\end{defin}

Then:

\begin{prop}
\label{pro.futfut}
$\partial^+C(\partial{S})$ is the set of points dual to spacelike support hyperplanes to 
$E(\partial {S})$ at a point of $\partial^-E(\partial{S})$.\fin
\end{prop}

\subsection{Cosmological time}

\begin{prop}
\label{pro.Ecosmologique}
$E(\partial{S})$ has regular cosmological time.
\end{prop}

 \preu
According to remark~\ref{deff}: 

\[ E(\partial\tilde{S}) = \{ (x,\theta) \in {\mathbb D}^{2} \times {\mathbb R} / f_{-}(x) < \theta < f_{+}(x) \} \]

Since it is contained in a de Sitter domain, we can assume that $f_\pm$ take value in $[-\pi/2, \pi/2]$.
We denote by $O$ the north pole of ${\mathbb D}^2$.

Let $p_0 = (x_0,\theta_0)$ be a point in $E(\partial\tilde{S})$. By the very
definition of $f_+$, for $x$ on the boundary 
$\partial {\mathbb D}^{2}$ we have 
$\theta_0 - d(x,x_0) < f(x)$. Hence, the map 
$x \to \theta_0 - d(x, x_0) - f_-(x) $ is negative 
on $\partial{\mathbb D}^{2}$. By compactness of 
$\partial{\mathbb D}^{2}$, there is a
compact subset $K$ of ${\mathbb D}^{2}$ such that each $x$ satisfying 
$f_-(x) < \theta_0 - d(x, x_0)$ belongs to $K$. 
But the past $I^-(p_0)$ of $p_0$ in $E(\partial\tilde{S})$ is the set of
points $(x, \theta)$ with $f_-(x) < \theta < \theta_0 - d(x, x_0)$: it follows that
$I^-(p_0)$ is contained in $K \times [-\pi/2, \pi/2]$. 
In particular, the conformal factor $\frac{1}{\cos^2(d(O,x)}$
of the AdS metric versus $ds_0^2 - dt^2$ (cf. remark~\ref{scan}) is uniformly bounded on
$I^-(p_0)$ by a factor $\mu^2$. It follows easily that the time length
of causal curves contained in $I^{-}(p_0)$ is uniformly bounded by
$\mu\pi$. Hence, the cosmological time $\tau$ on  
$E(\partial\tilde{S})$ has finite existence time.

Consider now an inextendible past oriented causal curve $c$ in
$E(\partial\tilde{S})$ starting from $p_0$. Forgetting the parametrization, 
such a causal curve
can always be expressed as the graph of a $1$-Lipschitz function $x:
[\theta_0, \theta_1[ \rightarrow {\mathbb D}^{2}$ (with $\theta_1 < \theta_0$) such that
$x(\theta_0) = x_0$ and:

\[ f_-(x(\theta)) < \theta < f_+(x(\theta)) \]

Since it is $1$-Lipschitz, $\theta \to x(\theta)$ admits a limit at $\theta_1$, that
we denote by $x_1 = x(\theta_1)$. Observe that $x_1$ must belong to the
compact $K$ defined above: in particular, $x_1$ belongs to ${\mathbb
  D}^{2}$. Moreover, since $c$ is inextendible, $p_1 = (x_1, \theta_1)$
belongs to the boundary of $E(\partial\tilde{S})$: $f_-(x_1) = \theta_1$.

Assume that $\tau$ does not converge to $0$ along $c$. Then, there is
some $\epsilon > 0$, and a sequence of $\theta_n$ converging to $\theta_1$ such
that each $\tau((x(\theta_n), \theta_n))$ is bigger than $\epsilon$. In other
words, there is a past-oriented causal curve $c_n$ starting from
$(x(\theta_n), \theta_n)$ contained in $E(\partial\tilde{S})$ and with time length equal
to  $\epsilon$. All these curves are contained in the past of
$p_0$. We have just proved that this past is contained in a
region $K \times [-\pi/2, \pi/2]$ where $K$ is compact, and where the
conformal factor between the AdS metric and the Ein metric is between
$1$ and $\mu^2$. Hence the curves $c_n$, considered as causal curves
in the Einstein universe ${\mathbb S}^{2} \times {\mathbb R}$ with
the metric $ds_0^2-dt^2$ have time length between $\epsilon$ and
$\mu\epsilon$.  Since the Einstein Universe is globally hyperbolic, it
follows that the curves $c_n$ converges in the Hausdorff topology to a
past-oriented causal curve $\bar{c}$ with time length (for the
Einstein metric) strictly positive and contained in $K \times [-\pi/2,
\pi/2]$. Obviously, the starting point of $\bar{c}$ must be $p_1$. 
Hence, the final part of $\bar{c}$ is in the (strict) past
of $p_1$. The same must be true for the $c_n$, for sufficiently big
$n$. In particular, the past of $p_1$ contains points of
$E(\partial\tilde{S})$. This is a contradiction, since the past of $p_1$ is
the domain $(x, \theta)$ with $\theta < \theta_1 - d(x, x_1)$, and every element
$(x, \theta)$ of $E(\partial\tilde{S})$ satisfy  $\theta > f_-(p) \geq f_-(x_1) - d(x,
x_1) = \theta_1 - d(x, x_1)$. 
\fin

\rque
The cosmological time is Lipschitz continuous, but it is not $C^1$ in general.
It can be proved that $\tau$ is $C^{1,1}$ on $\{ \tau < \pi/2 \}$, i.e. 
the past in $E(\partial{S})$ of $\partial^+C(\partial{S})$. See \cite{BenBon, BenBon2}.
\erque

\subsection{Generic achronal circles in $\mbox{Ein}_2$}

Most considerations above apply when $\partial S$ is any achronal topological circle in
$\mbox{Ein}_2$. Anyway:

\begin{prop}
\label{remplir1}
Every generic achronal topological circle of $\mbox{Ein}_{2}$ is the boundary of a smooth
spacelike hypersurface of AdS.
\end{prop}

\preu
A generic achronal topological circle correspond to the graph $\Lambda$ of some $1$-Lipschitz map $f: {\mathbb S}^{1} \rightarrow {\mathbb R}$. We can define the open set 
$E({\Lambda})$ as in remarks \ref{defabs} and \ref{deff} (the two definitions
still coincide). Following the second definition, it is
the open set contained between the graphs of two $1$-Lipschiz maps
$f_{-}$, $f_{+}$. The proposition is proved as soon
as we prove the existence of the smooth contracting map $g: {\mathbb D}^{2} \rightarrow {\mathbb R}$ 
with $f_{-} < g < f_{+}$. Indeed, such a $g$ will necessarily coincide with $f_{+}=f_{-}$ on 
$\partial{\mathbb D}^{2}$. 

After adding some positive contant, we can assume $f_{+}>0$.
Define then, for every integer $n$, 
$g_{n}(x)=\mbox{Sup}(f_{-}(x), (1-\frac{1}{n})f_{+}(x)-\frac{1}{n})$. 
Then, the sum $g(x)= \sum \frac{g_{n}(x)}{2^{n}}$ provides a contracting 
map $g$ with all the required property, except smoothness. 
This map can be approximated with a smooth one, still contracting 
(\cite{fathimaderna}).

There is another proof:
$E({\Lambda})$ is strongly causal since $\widehat{\mbox{Ein}}_3$ is strongly causal. It follows from the definition that if $x$, $y$ are points in $E({\Lambda})$, the intersection between the future (in $\widehat{\mbox{Ein}}_3$) of $x$ and the past (in $\widehat{\mbox{Ein}}_3$) of $y$ is contained in $E({\Lambda})$: it is compact and coincide with the intersection between the future in $E({\Lambda})$ of $x$ and the past in $E({\Lambda})$ of $y$. 
Hence, according to Theorem~\ref{gerochthm1}, $E(\Lambda)$ is globally hyperbolic. In particular, it admits a \emph{smooth} Cauchy surface $S$ (see remark~\ref{rk.smooth}).
Any inextendible causal curve intersect $S$; in particular, this is true
for the curves $\{ x \} \times ]f_{-}(x), f_{+}(x)[$. It follows that $S$
is the graph of a smooth contracting map ${\mathbb D}^{2} \rightarrow 
{\mathbb R}$ which extends on $\partial{\mathbb D}^{2}$ as $f$.

We propose now a third and last proof,  more adapted to the equivariant case to be considered later 
(see \S~\ref{subsub.cerclegh}): prove as for proposition~\ref{pro.Ecosmologique} that $E(\Lambda)$ 
has regular cosmological time and then apply Theorems~\ref{gerochthm1}, \ref{lem.cosmogood}.
\fin

\subsection{Invisible domains from achronal subsets of 
$\widehat{\mbox{Ein}}_{2}$ }
\label{sub.invisible}

In this {\S} $\widetilde{\Lambda}$ is a generic closed achronal subset of $\widehat{\mbox{Ein}}_{2}$.
We assume $\mbox{card}(\widetilde{\Lambda}) \geq 2$.
Then $\widetilde{\Lambda}$ is the graph of a $1$-Lipschitz map $f_0: \Lambda_{0} \rightarrow {\mathbb R}$, 
where $\Lambda_{0}$ is a closed subset of ${\mathbb S}^{1} = \partial{\mathbb D}^2$.
We can define as before the {\em invisible domain from $\widetilde{\Lambda}$ in $\widehat{\mbox{AdS}}$\/} 
that we denote by ${E}(\widetilde{\Lambda})$: it is the set 
$\{ (x,\theta) \in {\mathbb D}^{2} \times {\mathbb R} / f_{-}(x) < \theta < f_{+}(x) \}$ where:

- $f_{-}(x) = \mbox{Sup}_{y \in \Lambda_0} \{ f(y)-d(x,y) \}$,

- $f_{+}(x) = \mbox{Inf}_{y \in \Lambda_0} \{ f(y)+d(x,y) \}$.

Of course $f_\pm$ are actually defined on the closure $\overline{\mathbb D}^2 \times {\mathbb R}$. 
Define 
$\Omega(\widetilde{\Lambda}) = \{ (x,t) \in \partial{\mathbb D}^{2} \times {\mathbb R} / 
f_{-}(x) < t < f_{+}(x) \}$. 
It is the {\em invisible domain\/} from $\widetilde{\Lambda}$ in $\widehat{\mbox{Ein}}_{2}$.
Observe that $\Omega(\widetilde{\Lambda})$ is an open subset of $\widehat{\mbox{Ein}}_{2}$.

Moreover, if we consider $\widetilde{\Lambda}$ as a closed achronal subset of $\widehat{\mbox{Ein}}_{3}$
containing $\widehat{\mbox{Ein}}_{2}$ as the boundary of an anti-de Sitter domain,
then $E(\widetilde{\Lambda})$ is the intersection between the anti-de Sitter domain and the invisible domain 
from $\widetilde{\Lambda}$ in $\widehat{\mbox{Ein}}_{3}$.

Let $\widetilde{\Lambda}^\pm$ be the graphs of the restriction of $f_\pm$ to $\partial{\mathbb D}^2$.
These graphs are achronal topological circles  in $\widehat{\mbox{Ein}}_{2}$ containing $\widetilde{\Lambda}$. 

If we change the parametrization $\widetilde{\mbox{AdS}} \approx {\mathbb D}^2 \times {\mathbb R}$, we change 
of course the closed subset $\Lambda'_0 \subset \partial{\mathbb D}^2$, but there is a diffeomorphism from 
$\partial{\mathbb D}^2$ into itself mapping $\Lambda_0$ on $\Lambda'_0$. In particular, if the first coordinates 
of $x$, $y$ for the first parametrization are extremities of a connected component of 
$\partial{\mathbb D}^2 \setminus \Lambda_0$, then their projections in $\Lambda'_0$ are also 
extremities of a connected component of $\partial{\mathbb D}^2 \setminus \Lambda_0$. 

\begin{defin}
A \emph{gap pair\/} is a pair $(x, y)$ of points $\widetilde{\Lambda}$ corresponding to points $(x_0, y_0)$ in $\Lambda_0$ which are extremities of a connected component of $\partial{\mathbb D}^2 \setminus \Lambda_0$. 
An \emph{ordered gap pair} is the data of a gap pair with a connected component of
$\partial{\mathbb D}^2 \setminus \Lambda_0$ with extremities $x_0$, $y_0$.
A gap pair is \emph{achronal\/} if $x$, $y$ are not causally related in $\widehat{\mbox{Ein}}_2$. 
An ordered gap pair $(x,y, I)$ is \emph{lightlike\/} if $x$, $y$ are extremities of a lightlike segment in 
$\widehat{\mbox{Ein}}_2$ projecting in $\partial{\mathbb D}^2$ on $I$.
An ordered gap pair $(x,y, I)$ is \emph{extreme\/} if $x$, $y$ are extremities of a lightlike segment in 
$\widehat{\mbox{Ein}}_2$ projecting in $\partial{\mathbb D}^2$ on a segment disjoint from $I$.
\end{defin}

\rque
\label{rk.fillinggap}
It is quite clear that for any lightlike gap pair $(x,y)$, 
$\widetilde{\Lambda}$ and $\widetilde{\Lambda} \cup [x,y]$ 
define the same invisible domain. Hence, by adding all 
these lightlike segments and if the initial $\widetilde{\Lambda}$ 
did contains at least one pair of non-causally related points, we can 
reduce the study to the case where all gap pairs are achronal
or extreme.
Moreover, if  $\widetilde{\Lambda}$ contains at least three points, every 
gap pair defines an ordered gap pair. Finally, extreme gap pairs occurs only 
in the case where $\widetilde{\Lambda}$ is contained in a lightlike segment.
\erque

\begin{defin}
\label{def.gapsegment}
A gap segment is a segment $]x,y[$ where $(x,y)$ is a gap pair.
\end{defin}

Gap segments associated to achronal gap pairs are contained 
in ${\mathbb A}{\mathbb D}{\mathbb S}$. 
Observe that gap segments are contained in $\partial \mbox{Conv}(\Lambda)$.

According to lemma \ref{Edansds}, the projection of $\widetilde{\Lambda}$
in ${{\mathbb A}{\mathbb D}{\mathbb S}}$ is injective and the image $\Lambda$ is a compact subset in ${\mbox{Ein}_{2}}$. Moreover, since
$\Lambda$ is contained in a de Sitter domain, we can select some element $x_0$ of AdS
such that ${\Lambda}$ is contained in the affine patch $V(x_0)$.

\begin{lem}
\label{ddansds}
If $\widetilde{\Lambda}$ contains at least two non-causally related points, then $E(\widetilde{\Lambda})$ 
is contained in a de Sitter domain.
\end{lem}

\preu
$\widetilde{\Lambda}$ contains two elements which are not causally related, i.e. with coordinates
$(x, 0)$, $(-x,0)$ (remark \ref{parametrons}). Then, for any $(y, \theta)$ in
$E(\widetilde{\Lambda})$ we have:
$$\mbox{Sup} \{ -d(x,y), -d(-x,y) \} < \theta < \mbox{inf} \{d(x,y), d(-x,y) \}$$
Hence, $-\frac{\pi}{2} < \theta < \frac{\pi}{2}$ and the lemma follows.
\fin

Thanks to this lemma, we can project everything in AdS.
Remark \ref{deff} remains true: define $C({\Lambda}) = (\mbox{Conv}({\Lambda})^{\ast})^{\sharp}=\mbox{Conv}({\Lambda}^{\flat})^{\ast}$. The projection of ${E}(\widetilde{\Lambda})$ in this affine patch is the intersection $E({\Lambda})$ between ${{\mathbb A}{\mathbb D}{\mathbb S}}$ and all the half-spaces $E_{xy}$
where $x$, $y$ describes ${\Lambda} \times {\Lambda}$ minus the diagonal.
It follows that $E({\Lambda})$ is equal to 
$C({\Lambda}) \cap {{\mathbb A}{\mathbb D}{\mathbb S}}$.
$\Omega(\Lambda)$ projects in ${\mbox{Ein}_{2}}$
to $\Omega({\Lambda}) = {\mbox{Ein}_{2}} \cap C({\Lambda})$.

\subsection{Elementary cases}
\label{sub.elementary}

In this section, we assume that $\widetilde{\Lambda}$ is a generic achronal subset, containing 
at least two points and without lightlike gap pair (see remark~\ref{rk.fillinggap}).

\begin{defin}
The \emph{common future} of ${\widetilde{\Lambda}}$ is the set of points of 
$\widehat{\mbox{Ein}}_{2}$ containing the entire
${\widetilde{\Lambda}}$ in their past lightcones. We denote it by $X^+(\widetilde{\Lambda})$.
Similarly, the \emph{common past} $X^-(\widetilde{\Lambda})$ is the set of points in $\widehat{\mbox{Ein}}_{2}$ containing ${\widetilde{\Lambda}}$ in their future lightcones.
\end{defin}

\begin{defin}
When $X^+(\widetilde{\Lambda})$ or $X^-(\widetilde{\Lambda})$ is not empty, $\widetilde{\Lambda}$ is elementary.
If not, $\widetilde{\Lambda}$ is said nonelementary. 
\end{defin}

\begin{lem}
\label{lem.elementarycases}
The elementary case admits three subcases: 

\begin{itemize}

\item \emph{the conical case:\/} it is the case where $X^+(\widetilde{\Lambda})$ or
$X^-(\widetilde{\Lambda})$ is reduced to one point $x_0$ and $X(\widetilde{\Lambda}) = \{ x_0 \}$. 
If $X^+(\widetilde{\Lambda}) =\{ x_0 \}$:
let $L_1$, $L_2$ be the two past oriented closed lightlike 
segments in $\widehat{\mbox{Ein}}_2$ with extremities $x_0$, $x_1$, where $x_1$ is the 
point opposite to $x_0$ in the past.
Then, $\widetilde{\Lambda}$ is the union of two lightlike segments $I_1$, $I_2$, not reduced to single points, 
contained respectively in $L_1$, $L_2$.

The case $X^-(\widetilde{\Lambda}) =\{ x_0 \}$ admits a similar description: inverse the time orientation.

\item \emph{the splitting case:\/} it is the case where $X(\widetilde{\Lambda})$ is a pair of non-causally related points. $\widetilde{\Lambda}$ is then a pair of non-causally related points.

\item \emph{the extreme case:\/} it is the case where $X^+(\widetilde{\Lambda})$ and $X^-(\widetilde{\Lambda})$ are lightlike rays, contained in a lightlike geodesic $\Delta$. $\widetilde{\Lambda}$ is then a lightlike segment inside $\Delta$.

\end{itemize}

\end{lem}

\preu
Reversing the time orientation if necessary, we can assume that $X^+(\widetilde{\Lambda})$
contains a point $x_0$. Let $-x_0$ be the point opposite 
to $x_0$ in the past and $L_1$, $L_2$ the two lightlike segments with extremities 
$x_0$, $-x_0$. 
Then, $\widetilde{\Lambda}$ is contained in the pure lightlike circle 
$L_1 \cup L_2$. Since we assume that $\widetilde{\Lambda}$ is generic, 
it is not the union $L_1 \cup L_2$. Since we assume that there is no lightlike gap 
pair, the intersections 
$\widetilde{\Lambda}_i = \widetilde{\Lambda} \cap L_i \;\; (i=1,2)$ are connected, 
i.e. intervals. 

Assume that $\widetilde{\Lambda}$ contains at least two non-causally related points. 
Then, $\widetilde{\Lambda}_1$ and $\widetilde{\Lambda}_2$ have non-empty interior. 
If $X(\widetilde{\Lambda}) = \{ x_0 \}$, we are in the conical case.
If $X^+(\widetilde{\Lambda})$ contains another point $x'_0$, we can define similarly two 
lightlike segments $L'_1$, $L'_2$ and $\widetilde{\Lambda}_i$ must be contained in 
$L_i \cap L'_i$.
Since these intersections are not empty, $x_0$ and $x'_0$ are not causally related, 
and these intersections are both reduced to one point. The lemma follows in this case: 
we are in the splitting case.

The last case to consider is the case where all points in $\widetilde{\Lambda}$ are 
causally related one to the other. $\widetilde{\Lambda}$ is then equal to 
$\widetilde{\Lambda}_1$ or $\widetilde{\Lambda}_2$. We are in the extreme case.
\fin

\subsection{Description of the splitting case}
\label{subsub.split}

In \cite{ba2}, we will describe $E(\widetilde{\Lambda})$ for every elementary
$\widetilde{\Lambda}$. For the present paper, we just need to
understand the splitting case $\widetilde{\Lambda} = \{ x,y \}$, 
where $x$, $y$ are two non-causally related points in 
$\widehat{\mbox{Ein}}_{2}$. Then, $\{ x,y \}$ is a gap pair, 
and there are two associated ordered gap pairs, that we denote 
respectively by $(x,y)$ and $(y,x)$. $\widetilde{\Lambda}^+$ is 
an union ${\mathcal T}^+_{xy} \cup {\mathcal T}^+_{yx}$ of two 
nontimelike segments with extremities $x$, $y$, that we call 
\emph{upper tents\/.} Such an upper tent is the union of two 
lightlike segments, one starting from $x$, the other from $y$, 
and stopping at their first intersection point, that we call the 
\emph{upper corner.\/}

\begin{figure}[ht]
\centerline{\includegraphics[width=8cm]{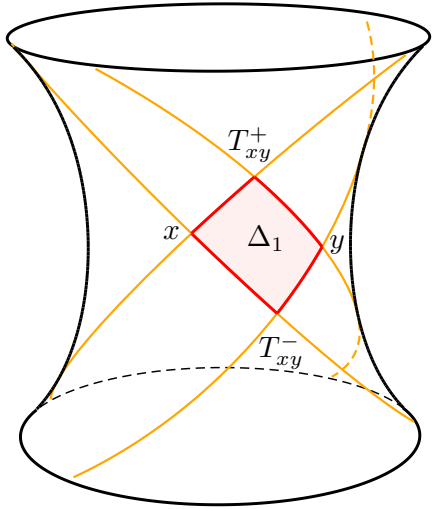}}
\caption{\label{ttents}
   \textit{Upper and lower tents}} 
\end{figure}

Similarly, $\widetilde{\Lambda}^-$ is an union ${\mathcal T}^-_{xy} \cup {\mathcal T}^-_{yx}$ of two \emph{lower tents\/} admitting a similar 
description, but where the lightlike segments starting from $x$, $y$ are now past oriented
(see Figure \ref{ttents}) and sharing a common extremity: the \emph{lower corner.\/}

The invisible domain $\Omega(\widetilde{\Lambda})$ from $\widetilde{\Lambda}$ in $\widehat{\mbox{Ein}}_{2}$ is the union of two diamond-shape regions $\widetilde{\Delta}_1$, $\widetilde{\Delta}_2$. The boundary of $\widetilde{\Delta}_1$ is the union ${\mathcal T}^+_{xy} \cup {\mathcal T}^-_{xy}$ and the boundary of $\widetilde{\Delta}_2$ is ${\mathcal T}^+_{yx} \cup {\mathcal T}^-_{yx}$. We project all the picture in some affine region $V \approx {\mathbb  R}^3$ of $S(E)$ such that:

-- $V \cap {\mathbb A}{\mathbb D}{\mathbb S}$ is the interior of the hyperboloid: $\{ x^2 + y^2 < 1 + z^2 \}$,

-- $\Lambda = \{ (1,0,0), (-1, 0, 0) \}$.

Then, $E(\Lambda)$ is region $\{ -1 < x < 1 \} \cap {\mathbb A}{\mathbb D}{\mathbb S}$.
One of the diamond-shape region $\widetilde{\Delta}_i$ projects to $\Delta_1 = \{ -1 < x < 1, y >0 , x^2+y^2 = 1+z^2 \}$, the other projects to $\Delta_2 = \{ -1 < x < 1, y <0 , x^2+y^2 = 1+z^2 \}$. The past of $\Delta_1$ in $E(\Lambda)$ is $P_1 = \{ (x,y,z) \in E(\Lambda) / z < y \}$. and the future of $\Delta_1$ in $E(\Lambda)$ is $F_1 = \{ (x,y,z) \in E(\Lambda) / z > -y \}$. We have of course a similar description for the future $F_2$ and the past $P_2$ of $\Delta_2$ in $E(\Lambda)$. Observe:

-- the intersections $F_1 \cap F_2$ and $P_1 \cap P_2$ are disjoint. They are tetraedra in $S(E)$:
$F_1 \cap F_2$ is the interior of the convex hull of $\Lambda^+$, 
and $P_1 \cap P_2$ is the interior of the 
convex hull of $\Lambda^-$.

-- the intersection $F_1 \cap P_1$ (resp. $F_2 \cap P_2$) is the intersection between ${\mathbb A}{\mathbb D}{\mathbb S}$ and the interior of a tetraedron in $S(E)$: the convex hull of $\Delta_1$ (resp. $\Delta_2$).

\begin{defin}
\label{defghsimple}
$E^+(\Lambda) = F_1 \cap F_2$ is the future globally hyperbolic convex core; $E^-(\Lambda) = P_1 \cap P_2$ 
is the past globally hyperbolic convex core.
\end{defin}

This terminology is justified by the following (easy) fact: $F_1 \cap F_2$ (resp. $P_1 \cap P_2$) is the invisible domain $E(\Lambda^+)$ (resp. $E(\Lambda^-)$). Hence, they are indeed globally hyperbolic.

The intersection between the closure of $E(\Lambda)$ in $S(E)$ and the boundary ${\mbox{Ein}_{2}}$ of ${\mathbb A}{\mathbb D}{\mathbb S}$ is the union of the closures of the diamond-shape regions. Hence, $\Delta_{1,2}$ can be thought as the conformal boundaries at infinity of $E(\Lambda)$. Starting from any point in $E(\Lambda)$, to $\Delta_i$ we have to enter in $F_i \cap P_i$, hence we can adopt the following definition:

\begin{defin}
\label{def.end}
$F_1 \cap P_1$ is an end of $E(\Lambda)$.
\end{defin}

Finally:

\begin{defin}
The future horizon is the past boundary of $F_1 \cap F_2$; the past horizon is the
future boundary of $P_1 \cap P_2$.
\end{defin}

\begin{prop}
\label{simpledecomp}
$E(\Lambda)$ is the disjoint union of the future and past globally hyperbolic cores $E^\pm(\Lambda)$, of the two ends and of the past and future horizons.\fin
\end{prop}

\begin{figure}[ht]
\centerline{\includegraphics[width=8cm]{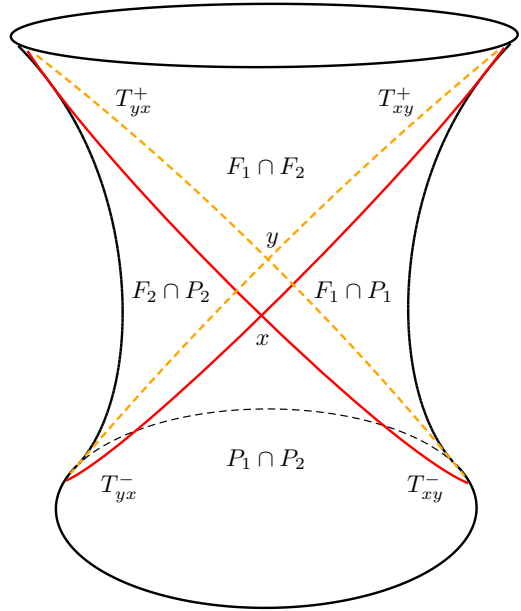}}
\caption{\label{splitfigure}
   \textit{The splitting case. The domain $E(\Lambda)$ is between 
the hyperplanes $x^\perp$ and $y^\perp$. These hyperplanes,
tangent to the hyperboloid, are not drawn, except their intersections with the hyperboloid,
which are the upper and lower tents ${\mathcal T}^\pm_{xy}$, ${\mathcal T}_{yx}^\pm$.}} 
\end{figure}

\rque
\label{rk.BTZnotations}
In the conventions of \cite{BTZ, BTZ2, brill}, the globally hyperbolic convex cores $F_1 \cap F_2$ and 
$P_1 \cap P_2$ are \emph{regions of type II,\/} also called \emph{intermediate regions.\/} 
The ends $F_1 \cap P_1$ and $F_2 \cap P_2$ are \emph{outer regions,\/} or \emph{regions of type I.\/}
\erque

\subsection{The nonelementary case}

From now, we assume that $\widetilde{\Lambda}$ is generic and 
nonelementary. Then,
it contains at least two non-causally related points: every 
$\widetilde{\Lambda}^\pm$ is generic. 
Moreover, every gap pair defines uniquely an ordered gap pair, that we can assume to be achronal 
(remark~\ref{rk.fillinggap}).

\begin{defin}
\label{def.proper}
$\widetilde{\Lambda}$ is proper if it is nonelementary and
not contained in a flat sphere.
\end{defin}

\begin{prop}
\label{2cas}
The closure of $E({\Lambda})$ in $S(E)$
is contained in an affine patch if and only if $\Lambda$ is proper.
\end{prop}

\rque
We have chosen the terminology so that $\Lambda$ is proper if and only if the convex $E(\Lambda)$ 
is proper in the meaning of \S~\ref{dual}.
\erque

\preud{2cas}
$E(\Lambda)$ is the intersection between $(\mbox{Conv}({\Lambda})^{\ast})^{\sharp}$ and 
${\mathbb A}{\mathbb D}{\mathbb S}$. It follows from Proposition~\ref{convpn} that
its closure is contained in an affine patch, except if $\mbox{Conv}({\Lambda})$ has empty interior, i.e. 
is contained in a projective hyperplane $P$. But then $P$ must be spacelike since it contains the nonelementary 
set $\Lambda$. Hence, $P \cap {\mbox{Ein}_{2}}$ is a flat sphere containing $\Lambda$.

Inversely, if $\Lambda$ is contained in a flat sphere $S(v^\perp_0 \cap \{ Q=0 \})$, then $v_0$ and $-v_0$ 
belong to the closure of $E(\Lambda)$.\fin

\rque
Observe that if $\Lambda$ is nonelementary but non-proper, the flat sphere containing it is unique.
\erque

\subsection{The decomposition in ends and globally hyperbolic cores}

We still assume that $\widetilde{\Lambda}$ is nonelementary.
It is easy to see that  $\widetilde{\Lambda}^+$ (resp. $\widetilde{\Lambda}^-$) 
is obtained from $\widetilde{\Lambda}$ by adding for any lightlike gap pair $(x,y)$ 
the lightlike segment $[x,y]$ and for any achronal gap pair $(x,y)$ the
upper (resp. lower) tent ${\mathcal T}^+_{xy}$ (resp. ${\mathcal T}^-_{xy}$).

\begin{figure}[ht]
\centerline{\includegraphics[width=8cm]{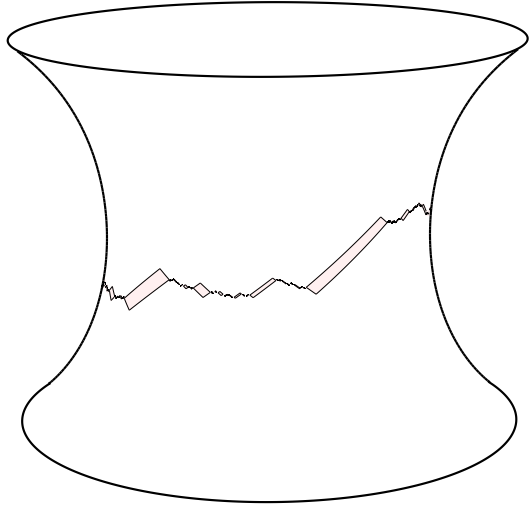}}
\caption{\label{completerlambda}
   \textit{Filling the gaps.}} 
\end{figure}

The connected components of $\Omega(\Lambda)$ are precisely the diamonds $\Delta_{xy}$. 
The convex hull in $S(E)$ of $\Delta_{xy}$ is a tetraedron (see \S~\ref{sub.elementary}), 
the intersection of this tetraedron with ${\mathbb A}{\mathbb D}{\mathbb S}$ has been described above.

\begin{defin}
For any gap pair $(x,y)$, the closed end $\overline{\mathcal E}_{xy}$ is the 
intersection between ${\mathbb A}{\mathbb D}{\mathbb S}$ and the convex hull
in $S(E)$ of $\Delta_{xy}$. 
\end{defin}

\begin{lem}
\label{boutdansE}
Every closed end $\overline{\mathcal E}_{xy}$ is contained in $E(\Lambda)$.\fin
\end{lem}

\preu
Since $E(\Lambda)$ is convex and since $\Delta_{xy}$ is contained in the closure of $E(\Lambda)$, 
$\overline{\mathcal E}_{xy}$ is contained in the closure of $E(\Lambda)$. 
Let $z^+$ (resp. $z^-$) be the upper (resp. lower) corner of $\Delta_{xy}$. 
Let $p = ax+by+cz^++dz^- \;\; (a,b,c,d \geq 0)$ be a point in $\overline{\mathcal E}_{xy}$. 
Observe that the norm of $p$ is $2ab\langle x \mid y \rangle + 2cd\langle z^- \mid z^+ \rangle$.
It has to be negative; since $\langle z^- \mid z^+ \rangle$ is positive, we have $ab > 0$.

If $p$ belongs to the boundary of $E(\Lambda)$ in ${\mathbb A}{\mathbb D}{\mathbb S}$, then there is an element 
$z$ of $\Lambda$ such that $\langle z \mid p \rangle =0$. But the scalar products of $z$ with $x$, $y$, 
$z^+$ and $z^-$ are all nonnegative.
Hence, $a \langle z \mid x \rangle$ and $b\langle z \mid y \rangle$ are both $0$. Since $a$ and $b$ are 
positive, it follows that $z$ is causally related to $x$ \emph{and\/} $y$. But $z^\pm$ are the only points in 
the affine patch in consideration which are causally related to $x$ and $y$ and they don't belong 
to $\Lambda$. We obtain a contradiction: $\overline{\mathcal E}_{xy}$ does not intersect the boundary 
of $E(\Lambda)$. The lemma follows.\fin

\begin{defin}
The end ${\mathcal E}_{xy}$ is the interior in ${\mathbb A}{\mathbb D}{\mathbb S}$ of 
$\overline{\mathcal E}_{xy}$.
\end{defin}

\begin{defin}
The future globally hyperbolic core is $E(\Lambda^+)$, the past globally hyperbolic core is $E(\Lambda^-)$.
\end{defin}

\begin{defin}
The past boundary in ${\mathbb A}{\mathbb D}{\mathbb S}$ of $E(\Lambda^+)$ is the past horizon. 
The future boundary in ${\mathbb A}{\mathbb D}{\mathbb S}$ of $E(\Lambda^-)$ is the future horizon.
\end{defin}

Future and past horizons are achronal proper hypersurfaces.

Finally, another important feature is the convex hull $\mbox{Conv}(\Lambda)$.
We denote $\partial{\mbox{Conv}(\Lambda)}$ its boundary in $S(E)$.

\begin{defin}
The edge part $\partial_{ed} \mbox{Conv}(\Lambda)$ of $\partial \mbox{Conv}(\Lambda)$ 
is the union of all gap segments (see definition~\ref{def.gapsegment}).
\end{defin}

\begin{lem}
\label{lem.C+-}
The edge part is the set of points of $\partial\mbox{Conv}(\Lambda) \cap {\mathbb A}{\mathbb D}{\mathbb S}$
admitting timelike support hyperplanes.
\end{lem}

\preu
Observe that every gap segment admits two lightlike hyperplanes: the dual planes to the middle points
$z^\pm$ of the associated upper and lower tents. The dual plane to the middle of $[z^-, z^+]$ is then
a timelike support hyperplane containing the gap segment.

Let's prove the reverse inclusion: let $q$ be an element of $S(E)$ such that the dual plane
$q^\perp$ is a timelike support hyperplane at a point $p = \sum_{i=1, ... , k}a_ip_i$ of
$\partial\mbox{Conv}(\Lambda) \cap {\mathbb A}{\mathbb D}{\mathbb S}$. Since $q^\perp$ is timelike,
$q^\perp \cap \mbox{Ein}_{2}$ is a copy of the Klein model of the anti-de Sitter space of dimension $1$.
It follows that the integer $k$ is less than $2$, and since $p$ belongs to ${\mathbb A}{\mathbb D}{\mathbb S}$,
$k = 2$. Moreover, since $q^\perp$ is a support hyperplane, $\Lambda$ is contained in one connected
component of ${\mbox{Ein}_{2}} \setminus q^\perp$. It follows that $( p_1, p_2 )$ is a gap pair, and
that $p$ belongs to the gap segment $]p_1, p_2 [$.\fin

In the proper case, i.e. when $\mbox{Conv}(\Lambda)$ has non-empty interior, 
$\Lambda \cup \partial_{ed} \mbox{Conv}(\Lambda)$ is a Jordan curve in the topological
sphere $\partial\mbox{Conv}(\Lambda)$ and 
$\partial \mbox{Conv}(\Lambda) \setminus (\Lambda \cup \partial_{ed} \mbox{Conv}(\Lambda))$ is 
the union of two (non-proper) topological discs. One of them - $C^+(\Lambda)$ - 
is in the future of the other, that we call $C^-(\Lambda)$.

In the non-proper case, i.e. when $\Lambda$ is contained in the
boundary of a totally geodesic embedding of ${\mathbb H}^2$ in ${\mathbb A}{\mathbb D}{\mathbb S}$,
$\mbox{Conv}(\Lambda)$ coincide with the convex hull of $\Lambda \subset \partial{\mathbb H}^2$
and the edge part $\partial_{ed} \mbox{Conv}(\Lambda)$ is the boundary of this convex hull
in ${\mathbb H}^2$. We simply define 
$C^+(\Lambda) = C^-(\Lambda) = \mbox{Conv}(\Lambda) \setminus (\partial_{ed} 
\mbox{Conv}(\Lambda) \cup \Lambda)$.

In both cases $C^{\pm}(\Lambda)$ are non timelike topological surfaces (lemma~\ref{lem.C+-}).

\begin{lem}
\label{lem.W}
Every timelike geodesic intersecting $C^{-}(\Lambda)$ intersects $C^{+}(\Lambda)$.
\end{lem}

\preu
Let $c$ be a timelike geodesic intersecting $C^{-}(\Lambda)$. Since $C^-(\Lambda)$ is nontimelike, 
$c$ must enter in $\mbox{Conv}(\Lambda)$. The only possible exit for $c$ is then through $C^+(\Lambda)$.\fin

\begin{cor}
\label{cor.W}
Let $W$ be the projection of $C^{-}(\Lambda) \subset {\mathbb A}{\mathbb D}{\mathbb S} 
\approx {\mathbb D}^2 \times {\mathbb S}^1$ on the first factor ${\mathbb D}^2$. Then, 
it is also the projection of $C^+(\Lambda)$ on the first factor. Every connected component 
of $\partial W$ is a curve $c$ in ${\mathbb D}^2$ joining two elements $x$, $y$ in 
$\Lambda_0 \subset \partial{\mathbb D}^2$ such that:

-- $]x,y[$ is a gap of $\Lambda_0$,

-- the curve $c$ disconnect $W$ from the gap segment $]x,y[ \subset \partial{\mathbb D}^2$.

\end{cor}

\preu
The first statement follows immediatly from Lemma~\ref{lem.W}. Then, $C^{\pm}(\Lambda)$ are 
graphs of functions $g_\pm: W \rightarrow \mathbb R$. The proofs of the topological
description of $W$, which follows from the concave-convex properties of $C^\pm(\Lambda)$,
are left to the reader.\fin

\rque
Be careful! $W$ is \emph{not} in general the convex hull in ${\mathbb D}^2$ of $\Lambda_0$.
Observe that $W \subset {\mathbb D}^2$ depends on the selected conformal parametrization
${\mathbb A}{\mathbb D}{\mathbb S} \approx {\mathbb D}^2 \times {\mathbb S}^1$.
\erque

\begin{prop}
\label{decompropre}
$E(\Lambda)$ is the union of the past and future globally hyperbolic cores with the closed ends 
associated to gap pairs.
\end{prop}

\preu
One of the inclusion follows from Lemma~\ref{boutdansE} and the obvious inclusion 
$E(\Lambda^+) \cup E(\Lambda^-) \subset E(\Lambda)$. 

For the reverse inclusion: first observe that if $\Lambda$ is a topological circle, 
it has no ends and $E(\Lambda) = E(\Lambda^+) = E(\Lambda^-)$: there is nothing to prove.
Hence, we assume that $\Lambda$ admits at least one (achronal) gap. 
Recall that, in a suitable conformal domain $\approx  {\mathbb D}^2 \times [-\pi/2, \pi/2]$, 
the invisibility domain $E(\Lambda)$ is the domain in ${\mathbb D}^2 \times [-\pi/2, \pi/2]$ 
between the graphs of $f_\pm$. Define:

- $F_{-}(x) = \mbox{Sup}_{y \in \partial{\mathbb D}^2} \{ f_+(y)-d(x,y) \}$,

- $F_{+}(x) = \mbox{Inf}_{y \in \partial{\mathbb D}^2} \{ f_-(y)+d(x,y) \}$.

Then, $E(\Lambda^+)$ is the domain between the graphs of $F_-$, $f_+$, 
and $E(\Lambda^-)$ is the domain between the graphs of $f_-$, $F_+$. 

The discs $C^\pm(\Lambda)$ are the graphs of two functions 
$g_\pm: W \rightarrow ]-\pi/2, \pi/2[$ described in corollary~\ref{cor.W}.

\emph{Claim: for every $x$ in $W$, $F^+(x) > F^-(x)$.}

Since $\mbox{Conv}(\Lambda) \subset \mbox{Conv}(\Lambda^\pm) \subset E(\Lambda^\pm)$, we have:

\[ \forall x \in W, \;\; f_-(x) \leq F_-(x) \leq g_-(x) \leq g_+(x) \leq F_+(x) \leq f_+(x) \]

In the proper case, we actually have $g_+(x) > g_-(x)$ for $x$ in $W$: the claim follows. In the non-proper case, we can select the de Sitter domain so that $f (y)= 0$ for every $y$ in $\Lambda_0$. Then, for every $x$ in $W$, we have $g_+(x) = g_-(x) = 0$. It follows that if $F_+(x_0) = F_-(x_0)$ for some $x_0$ in $W$, then this common value is $0$.
By definition of $F_\pm$, $0$ is the supremum of $f_+(y) - d(x_0, y)$ and this supremum is attained at some $y_0$. Then, $f_-(y_0) = - f_+(y_0) = - d(x_0, y_0)$.
It means that $(y_0, d(x_0, y_0))$ is the upper corner of a upper tent associated to some gap $(x,y)$ and $(y_0, -d(x_0, y_0))$ is the lower corner of the associated lower tent. It implies that $(x_0,0)$ belongs to the gap segment $]x,y[$. Contradiction.

The claim is proved. Let now $(x, \theta)$ be a point in $E(\Lambda)$: $f_-(x) < \theta < f_+(x)$ holds.

If $x$ belongs to $W$, then since $F_-(x) < F_+(x)$, we have either $f_-(x) < \theta < F_+(x)$, 
or $F_-(x) < \theta < f_+(x)$. In the former case, $(x, \theta)$ belongs to $E^-(\Lambda)$ and in the later case, 
$(x, \theta)$ belongs to $E^+(\Lambda)$.

According to Lemma~\ref{boutdansE}, the same conclusion holds if $x$ belongs to $\overline{W}$. 

Assume now $x \in {\mathbb D}^2 \setminus \overline{W}$. $W$ is a topological disc in 
${\mathbb D}^2$ with boundary components the projection of gap segments. Hence $x$ 
belongs to the connected component of ${\mathbb D}^2 \setminus l$ which does not 
contain $W$, where $l$ is the projection of a gap segment $]y, z[$. 
Let $z^+$ be the upper corner of the upper tent ${\mathcal T}^+_{yz}$ and let $z^-$ be the 
lower corner of the lower tent ${\mathcal T}^-_{yz}$. 
Since $(x,\theta)$ cannot be related to $y$ or $z$ we easily infer that
it must belong to the intersection between the past of $z^+$ and the future of $z^-$, i.e. 
to the closed end $\overline{\mathcal E}_{yz}$.\fin

\rque
Furthermore, it follows quite easily from the proof above that the intersection $E^+(\Lambda) \cap E^-(\Lambda)$ 
is not empty: actually, it can be proved that, in the notation used in the proof of \ref{decompropre}, 
$F_+ > F_-$ on $W$, $F_+ = F_-$ on $\partial W$ and $F_+ < F_-$ on ${\mathbb D}^2 \setminus \overline{W}$.
\erque

\rque
\label{rk.unionsurfaces}
There are other ways to characterize $E(\Lambda)$. For example, it is the union of every proper 
spacelike surfaces containing $\Lambda$ in their natural extensions in $\partial \mbox{AdS}$.
\erque

\section{Synchronized isometries of AdS  }
\label{sec.cyclic}

We use the identification 
$\overline{{\mathbb A}{\mathbb D}{\mathbb S}} \approx G = 
\mbox{PSL}(2, {\mathbb R})$ (cf. \S~\ref{sub.psl}). Then, 
$\widetilde{AdS}$ can be identified with the universal covering 
$\widetilde{G} = \widetilde{\mbox{SL}}(2, {\mathbb R})$. Denote 
by $\bar{p}: \widetilde{G} \rightarrow G$ the covering map, and $Z$ 
the kernel of $\bar{p}$: $Z$ is cyclic, it is the center of $\widetilde{G}$. 
Let $\delta$ be a generator of $Z$: we select it in the future of the 
neutral element $id$.

$\widetilde{G} \times \widetilde{G}$ acts by left and right translations 
on $\widetilde{G}$.
This action is not faithfull: the elements acting trivially are precisely 
the elements in $\mathcal Z$, the image of $Z$ by the diagonal embedding. 
The isometry group $\widehat{\mbox{SO}}_0(2,2)$ is then identified with 
$(\widetilde{G} \times \widetilde{G})_{/\mathcal Z}$.

Let $\mathcal G$ be the Lie algebra $\mbox{sl}(2, {\mathbb R})$ of $G$: 
the Lie algebra
of $(\widetilde{G} \times \widetilde{G})_{/\mathcal Z}$ is 
${\mathcal G} \times {\mathcal G}$. 
We assume the reader familiar with the notion of elliptic, parabolic, 
hyperbolic elements of $\mbox{PSL}(2, {\mathbb R})$. Observe that 
hyperbolic (resp. parabolic) elements of 
$\mbox{PSL}(2, {\mathbb R})$ are the exponentials $\exp(A)$ of 
\emph{hyperbolic\/} (resp. \emph{parabolic,\/} \emph{elliptic\/}) 
elements of ${\mathcal G} = \mbox{sl}(2, {\mathbb R})$, i.e., such that 
$\mbox{det}(A) < 0$ (resp. $\mbox{det}(A)=0$, $\mbox{det}(A)>0$).

\begin{defin}
An element of $\widetilde{G}$ is \emph{hyperbolic\/} 
(resp. parabolic, elliptic) if it is the exponential 
of a hyperbolic (resp. parabolic, elliptic) element of 
$\mathcal G$. 
\end{defin}

\rque
There is another possible definition through the identification 
$\widetilde{G} \approx \widetilde{\mbox{AdS}}$: hyperbolic (resp. parabolic) 
elements of $\mathcal G$ are spacelike (resp. lightlike) tangent vectors to 
$\widetilde{\mbox{AdS}}$ at the neutral element $id$.
Hyperbolic elements of $\widetilde{\mbox{AdS}}$ are the elements which are 
not causally related to $id$ of $\widetilde{G}$. Parabolic elements are 
points in the lightcone of $id$. Hence, their union is the set of points 
in $\widetilde{\mbox{AdS}}$ which are not strictly causally related to $id$. 
In particular, they 
belong to the affine domain associated to $id$.
\erque

\rque
\label{nonsync}
Elements of $\widetilde{G}$ which are not hyperbolic, parabolic 
or elliptic have the form $\delta^k \gamma'$, where $\delta^k$ is 
a non-trivial element of the center of $\widetilde{G}$, and $\gamma'$ 
a hyperbolic or parabolic element of $\widetilde{G}$.
\erque

\rque
\label{forme}
Every element of $(\widetilde{G} \times \widetilde{G})_{/\mathcal Z}$ 
can be represented by
a pair $(\gamma_L, \gamma_R)$ such that:

-- $\gamma_L$ is the exponential of an element of $\mathcal G$,

-- $\gamma_R = \gamma'_R \delta^k$, where $\gamma'_R$ is the exponential 
of an element of
$\mathcal G$, and $\delta^k$ an element of $Z$.

\erque

\begin{defin} 
\label{def.synch}
An element $\gamma = (\gamma_{L}, \gamma_{R})$ of $\widetilde{G} \times \widetilde{G}$ is synchronised if, up to a permutation of left and right components, it has one of the following form:

\begin{itemize}

\item (hyperbolic translation): $\gamma_L$ is trivial and $\gamma_R$ is hyperbolic,

\item (parabolic translation): $\gamma_L$ is trivial and $\gamma_R$ is parabolic,

\item (hyperbolic - hyperbolic) $\gamma_L$ and $\gamma_R$ are both non-trivial and hyperbolic,

\item (parabolic - hyperbolic) $\gamma_L$ is parabolic and $\gamma_R$ is hyperbolic,

\item (parabolic - parabolic) $\gamma_L$ and $\gamma_R$ are both non-trivial and 
parabolic, 

\item (elliptic) $\gamma_L$ and $\gamma_R$ are elliptic elements conjugate in $\widetilde{G}$. 

\end{itemize}

An element $\gamma$ of $(\widetilde{G} \times \widetilde{G})_{/\mathcal Z}$
is synchronized if it is represented by a synchronized element of $\widetilde{G} \times \widetilde{G}$.

\end{defin}

We will see that synchronized isometries are precisely those preserving some generic achronal subset. This statement essentially follows from the lemma:

\begin{lem}
\label{lem.syncaffine}
An isometry $\gamma$ is synchronized if and only if there is an affine domain $U$ in $\widetilde{\mbox{AdS}}$ such that $\gamma^n(U) \cap U \neq \emptyset$ for every $n$ in $\mathbb Z$.
\end{lem}

\preu
Assume that $\gamma$ is synchronized. Consider first the case where $\gamma_L$ or $\gamma_R$ is a nontrivial elliptic element. Then, after conjugacy, we can assume $\gamma_L = \gamma_R$. Then, $\gamma$ preserves the affine domain $A(id)$.

Consider now the case where $\gamma_R$ and $\gamma_R$ are not elliptic.
After conjugacy in $\widetilde{G} \times \widetilde{G}$, we can assume that $\gamma_L$ and $\gamma_R$ are exponentials of matrices of the form:

\[ \left(\begin{array}{cc}
  a & b \\
 0 & -a \end{array}\right) \]

In particular, $\gamma_L\gamma_R^{-1}$ is the exponential of a matrix $X$ of the form above.

Consider the affine domain $A(id)$. Then, for every integer $n$, $(\gamma^n_R, \gamma_R^n) A(id) = A(id)$, and $(\gamma_L^n\gamma_R^{-n}, id )A(id) = A(\gamma_L^n\gamma_R^{-n}) = A(\exp(nX))$. Since $\exp(nX)$ belongs to $A(id)$, we obtain that $\gamma A(id) \cap A(id) \neq \emptyset$.

Assume now that $\gamma$ is \emph{not\/} synchronized, but that there exists a affine domain 
$A$ such that all the $\gamma^nA \;\; (n \in {\mathbb Z})$ intersect $A$. 
There is an integer $q$ such that $\delta^qA$ intersect $A(id)$. Then, every $\gamma^n\delta^qA$ intersects $\delta^qA$. It implies that all the $\gamma^n\delta^qA$ are contained in the past of $\delta^3 A(id)$ and the future of $\delta^{-3} A(id)$. It follows that all the $\gamma^n A(id)$ are contained in the past of $\delta^6 A(id)$ and in the future of $\delta^{-6} A(id)$.

Select a representant $(\gamma_L, \gamma_R)$ of $\gamma$ as in remark~\ref{forme}.
According to remark~\ref{nonsync}, we have three cases to consider:
\begin{enumerate}
\item $\gamma_L$ and $\gamma'_R$ are parabolic or hyperbolic, but $\delta^k$ is not trivial,
\item $\gamma_L$ is elliptic and $\gamma_R$ is the exponential of an element of $\mathcal G$, but not conjugate to $\gamma_L$,
\item $\gamma_L$ is parabolic or hyperbolic, but $\gamma_R = \gamma'_R$ is elliptic.
\end{enumerate}

In the first case, $\gamma' = (\gamma_L, \gamma'_R)$ is synchronized: hence, for every integer $n$, $\gamma'_n A(id) \cap A(id) \neq \emptyset$. 
The affine domain $\gamma^n A(id) = \delta^{kn} \gamma'_{n} A(id)$ intersect $\delta^{kn} A(id)$. Since $k \neq 0$ - let's say, $k>0$ -  if $n$ is sufficiently big, affine domains intersecting $\delta^{kn} A(id)$ cannot be contained in the past of $\delta^6 A(id)$. Contradiction.

Consider now the second case. The first subcase is the case where $\gamma_R$ is elliptic too.
Moreover, after conjugacy, we can assume that $\gamma_L$ and $\gamma_R$ commute.
Then, $\gamma' = (\gamma_R, \gamma_R)$ is synchronized, and $\gamma$ is the composition of $\gamma'$ with the left translation by the nontrivial elliptic element $\gamma_L\gamma_R^{-1}$. We can identify $\widetilde{\mbox{AdS}}$ with ${\mathbb D}^2 \times \mathbb R$ such that the left translation by  $\gamma_L\gamma_R^{-1}$ is a non-trivial translation along the $\mathbb R$-factor.
It follows that for $n$ sufficiently big $\gamma^n A(id)$ is not contained in the past of $\delta^6 A(id)$ and the past of $\delta^{-6} A(id)$. Contradiction.

Assume now that we are still in the second case, but with $\gamma_R$ nonelliptic: $\gamma' = (id, \gamma_R)$ is synchronized, and $\gamma$ is the composition of $\gamma'$ with the left translation by $\gamma_L$. We obtain a contradiction as above.

The third case reduce to the second one after composing with the involution $g \mapsto g^{-1}$, which permutes $\gamma_R$ with $\gamma_L$.\fin

\section{Invariant achronal subsets}
\label{achronalinvariant}

Let $\Gamma$ be a subgroup of $\widehat{\mbox{SO}}_0(2,2)$ preserving a generic closed 
achronal subset $\widetilde{\Lambda}$ of $\widehat{\mbox{Ein}}_2$. We assume that
$\widetilde{\Lambda}$ is non-elementary, and without lightlike pairs (we recall once more 
remark~\ref{rk.fillinggap}). 

We now consider any discrete subgroup $\Gamma$ of $\widehat{\mbox{SO}}_0(2,2)$. 
According to Lemma~\ref{Edansds} and lemma~\ref{lem.syncaffine}, 
every element of $\Gamma$ is synchronized.
We assume moreover that $\Gamma$ is torsion free: it follows that $\Gamma$ does 
not contain synchronized elements $(\gamma_L, \gamma_R)$, where $\gamma_L$ and $\gamma_R$ 
are elliptic elements of $\widetilde{G}$ conjugate in $\widetilde{G}$. Indeed, the torsion 
free hypothesis prevents $\gamma_L$, $\gamma_R$ to have rationnal rotation angle, and if 
this rotation angle was irrationnal, $\Gamma$ would not be discrete.

The action of $\Gamma$ on $\widetilde{\mbox{AdS}}$ and $\widehat{\mbox{Ein}}_2$ 
preserves the invisible domains $\Omega(\widetilde{\Lambda})$ and 
$E(\widetilde{\Lambda})$.

\begin{thm}
\label{noneleOK}
Let $\widetilde{\Lambda}$ be a nonelementary generic achronal subset, 
preserved by a torsionfree discrete group $\Gamma \subset \mbox{SO}_0(2,2)$. 
Then, the action of
$\Gamma$ on $\Omega(\widetilde{\Lambda})$ and $E(\widetilde{\Lambda})$ 
are free, properly discontinuous, and the quotient spacetime 
$M_{\widetilde{\Lambda}}(\Gamma) = 
\Gamma\backslash{E}(\widetilde{\Lambda})$ is strongly causal.
\end{thm}

\subsection{The hyperbolic-hyperbolic case}

We first consider a special case, which is truely speaking an elementary one, 
but which is necessary to consider for the general case:

\begin{lem}
\label{lem.splitstronglycausal}
Assume that $\Gamma$ is cyclic, generated by some $\gamma  = (\gamma_L, \gamma_R)$, and 
that:

--  $\gamma _R = id$ or, 

-- $\gamma_L = id$ or, 

-- $x$ (resp. $y$) is an attractive (resp. repulsive) fixed point of $\gamma$.

Then, the quotient space $M_{xy}(\Gamma) = \Gamma\backslash{E}(x,y)$ is strongly causal. 
Moreover, the projections in $M_{xy}(\Gamma)$ of the open ends $P_i \cap F_i$, and 
of the globally hyperbolic convex cores $P_1 \cap P_2$, $F_1 \cap F_2$, are all 
causally convex domains.
\end{lem}

\preu
The strong causality is proved in \cite{ba2} (Proposition $4.11$, cases $(2)$ and $(6)$).
The causal convexity is quite obvious.
\fin

\subsection{The globally hyperbolic case}
\label{subsub.cerclegh}

According to
proposition~\ref{pro.Ecosmologique}, $E(\widetilde{\Lambda})$ is globally hyperbolic. 
We prove here :

\begin{prop}
\label{pro.gh}
If $\widetilde{\Lambda}$ is a topological circle, then the action of $\Gamma$ on $E(\widetilde{\Lambda})$ is free and properly discontinuous, and the quotient space $\Gamma\backslash{E}(\widetilde{\Lambda})$ is globally hyperbolic, with regular cosmological time.
\end{prop}

\proof
We first prove the properness of the action:

\emph{The flat case:} This is the case where  $\widetilde{S}$ is a round circle, i.e. the boundary $\partial \widetilde{S}_0$ of a totally geodesic isometric copy of ${\mathbb H}^2$ in $\widetilde{\mbox{AdS}}$. Then, $\Gamma$ preserves $\tilde{x}^+_0$, the point dual to $\widetilde{S}_0$ such that $\widetilde{S}_0$ is the past boundary of the affine domain $A(\tilde{x}_0)$ (see \S~\ref{sub.affine}). Project everything in ${\mathbb A}{\mathbb D}{\mathbb S}$. Select a basis on $E$ so that $x_0 = p(\tilde{x}_0)$ has coordinates $(1,0, ... ,0)$. Then, the stabilizer of $x_0$ is $\mbox{SO}_0(1,2)$, and $\Gamma$ is a torsion-free discrete subgroup of $\mbox{SO}_0(1,2)$. Hence, the action of $\Gamma$ on $S_0 = p(\widetilde{S}_0)$ is free and properly discontinuous. Our claim then follows from
the decomposition $x= \cos(\theta)y + \sin(\theta)x_0$ valid for any element $x$ of $E(S_0)$.

\emph{The non flat case:}
When $\widetilde{\Lambda}$ is not a round circle, it is proper (see definition~\ref{def.proper}).
We observe, as in the proof of Proposition~\ref{2cas} that 
$E(\Lambda) = p(E(\widetilde{\Lambda}))$ is $(\mbox{Conv}({\Lambda})^{\ast})^{\sharp}$, which is a proper convex domains in $S(E)$: the Hilbert metric on it is a well-defined metric (see \cite{metriquehilbert}). It follows that the action of $\Gamma$ on it is properly discontinuous. Observe that the action is free since $\Gamma$ has no torsion.

Hence, in any case, the quotient $M_{\widetilde{\Lambda}}(\Gamma)$ is a well-defined locally AdS spacetime. The proposition then follows immediatly  from Propositions~\ref{pro.Ecosmologique}, \ref{pro.quotientcosmook}, and Theorem~\ref{lem.cosmogood}.\fin

\subsection{The general case}
Even when $\widetilde{\Lambda}$ is not a topological circle, we can prove as in the section~\ref{subsub.cerclegh} that if $\widetilde{\Lambda}$ is non proper, $\Gamma$ acts freely and properly discontinuously on $E(\widetilde{\Lambda})$ by considering the Hilbert metric on the proper convex domains 
$(\mbox{Conv}({\Lambda})^{\ast})^{\sharp}$, since $E(\widetilde{\Lambda})$ is the intersection between this proper convex domain and ${\mathbb A}{\mathbb D}{\mathbb S}$.
We leave to the reader the proof of the properness of the action in the flat case.

According to Proposition~\ref{pro.strongconvex}, $M_{\widetilde{\Lambda}}(\Gamma)$ is strongly causal if and only if any point $x_0$ in $M_{\widetilde{\Lambda}}(\Gamma)$ admits a causally convex neighborhood. 
If $x_0$ belongs to the projection of the future or the past globally hyperbolic core, this projection is the required causally convex neighborhood. If not, $x_0$ belongs to the projection of a closed end $\overline{\mathcal E}_{xy}$ (cf. proposition~\ref{decompropre}). 

Observe that $\Gamma$ permutes the gaps, hence 
$\gamma\overline{\mathcal E}_{xy} \cap \overline{\mathcal E}_{xy} \neq \emptyset$ implies 
$\gamma\overline{\mathcal E}_{xy} = \overline{\mathcal E}_{xy}$. Moreover, 
in this situation, $\gamma$ preserves the gap segment $[x,y]$: since the action on this 
segment must be free and proper, the stabilizer 
$\Gamma_0$ of $\overline{\mathcal E}_{xy}$ is a trivial or cyclic group. 
In the last case, since it admits two non-causally related fixed points in 
$\widehat{Ein}_2$, elements in $\Gamma_0$ have the form $(\gamma_L, \gamma_R)$ 
where $\gamma_L$, $\gamma_R$ are both hyperbolic (one maybe trivial). 

We can be slightly more precise: the projection of 
$\overline{\mathcal E}_{xy}$ is closed. Indeed, let  $x_n$ be a sequence in $\overline{\mathcal E}_{xy}$, and $\gamma_n$ a sequence in $\Gamma$ such that $\gamma_n x_n$ converge to some point $\bar{x}$ in $E(\widetilde{\Lambda})$. Then, $\bar{x}$ belongs to some closed end $\overline{\mathcal E}_{x'y'}$ since the complement of the union of closed ends is open (it is the union of the globally hyperbolic convex cores). If $\bar{x}$ is in the interior of $\overline{\mathcal E}_{x'y'}$ , then $\overline{\mathcal E}_{x'y'} = \gamma_n \overline{\mathcal E}_{xy}$ for every sufficiently great $n$. The claim follows. 
If $\bar{x}$ is on the boundary of $\overline{\mathcal E}_{x'y'}$, then it is in the lightcone of 
some corner point $z'$ of a upper or lower tent. In other words, $\bar{x} = ax'+ by' +cz'$ with 
$a,b > 0, c \geq 0$. But every $x_n$ can be written: $x_n = a_nx + b_ny + c_nz^+ + d_nz^-$.
Hence, $\langle z \mid \gamma_n x_n \rangle$ is the sum of the nonpositive terms 
$a_n\langle \gamma_nx \mid z'\rangle$, $b_n\langle \gamma_ny \mid z'\rangle$, 
$c_n\langle \gamma_nz^+ \mid z'\rangle$ and $d_n\langle \gamma_nz^- \mid z'\rangle$. 
All these terms have to tend to $0$: it follows that $\gamma_n^{-1}z'$ tends to $z^+$ or $z^-$.
But corner points are isolated in $\widetilde{\Lambda}^\pm$. Hence, $\overline{\mathcal E}_{x'y'} = \gamma_n \overline{\mathcal E}_{xy}$ for every sufficiently great $n$, and the claim follows in this case too.

Hence, we can associate to every closed end $\overline{\mathcal E}_{xy}$ an open neighborhood $W_{xy}$ in $E(\widetilde{\Lambda})$ such that $\gamma W_{xy} \cap W_{xy} \neq \emptyset$ implies $\gamma W_{xy}= W_{xy}$.

Now, $\Omega(\widetilde{\Lambda})$ is contained in $\Omega(x,y) = \Delta_1 \cup \Delta_2$, where we can consider that $\Delta_1$ is the conformal boundary of the end $\overline{\mathcal E}_{xy}$. Since $(x,y)$ is a gap pair, $\widetilde{\Lambda}$ is a $\Gamma_0$-invariant closed achronal subset contained in the closure of $\Delta_2$. Moreover, since it is nonelementary, $\widetilde{\Lambda}$ has a non-trivial intersection with $\Delta_1$. Since it is achronal, it follows that the hypothesis of lemma~\ref{lem.splitstronglycausal} are fulfilled: either $\Gamma_0$ is a subgroup of $G_L$ or $G_R$, or $x$, $y$ are attractive or repulsive fixed points of every element of $\Gamma_0$.

Finally, there is a neighborhood $W'_{xy}$ of $\overline{\mathcal E}_{xy}$ in $E(x,y)$ such that $W_{xy}$ contains the intersection $W'_{xy} \cap E(\widetilde{\Lambda})$. According to lemma~\ref{lem.splitstronglycausal}, the neighborhood $W'_{xy}$ can be selected so that its projection in $\Gamma_0\backslash{E}(x,y)$ is a causally convex domain. 
Then, the projection of $W'_{xy} \cap E(\widetilde{\Lambda})$ in $M_{\widetilde{\Lambda}}(\Gamma)$ is a causally convex neighborhood of $x_0$.\fin

\rque ({\bf The case with torsion })
The results above can be extended to the case with torsion: the action of a discrete group $\Gamma$ on the invisibility domain of a nonelementary achronal subset $\widetilde{\Lambda}$ is properly discontinuous and \emph{strongly causal\/} in the following meaning:

\begin{defin}
\label{strongcausal}
The action of a group $\Gamma \subset \widehat{SO}_0(2,2)$ on an open subset $E$
of $\widetilde{\mbox{AdS}}$ is strongly causal if every element $x$ of $E$ admits an open neighborhood 
$U$ such that, for every element $\gamma$ of $\Gamma_{0}$, either ${\mathbb A}{\mathbb D}{\mathbb S}$ is a fixed point of
$\gamma$, or no element of $U$ is causally related to an element $\gamma(U)$.
\end{defin}

The quotient $\Gamma\backslash{E}(\widetilde{\Lambda})$ is a AdS-spacetimes with singularities, the singularities (``particles'') being timelike lines. Observe that when $\Gamma$ is finitely generated this quotient is finitely covered by a AdS spacetime without singularity, since according to Selberg lemma the discrete group $\Gamma$ contains then a finite index torsion-free subgroup.
\erque

\subsection{Existence of invariant achronal subsets}
\label{lambdafromgamma}
We have proved that in most cases, if a discrete group $\Gamma \subset \widehat{\mbox{SO}}_0(2,2)$ preserves a generic closed achronal subset $\widetilde{\Lambda}$ containing at least two points, then the action of $\Gamma$ on $E(\widetilde{\Lambda})$ is proper and strongly causal. We now try to answer to the question: \emph{given a torsion-free discrete  subgroup $\Gamma$ of $\widehat{\mbox{SO}}_0(2,2)$, is there a $\Gamma$-invariant generic closed achronal subset of $\widehat{\mbox{Ein}}_2$?\/}
According to lemmas~\ref{Edansds} and \ref{lem.syncaffine}, in order to preserve such an achronal subset, every element of $\Gamma$ must be synchronized, and $\Gamma$ projects injectively in $G \times G$, with $G = \mbox{PSL}(2,{\mathbb R})$. Moreover, this projection must be faithfull, with discrete image. Hence, $\Gamma$ has to be the image of some faithfull morphism $\rho: \Gamma \rightarrow G \times G$. 

We then reformulate the question above in the following way:

\begin{defin}
\label{def.admissible}
Let $\rho_L: \Gamma \rightarrow G$ and $\rho_R: \Gamma \rightarrow G$ two morphisms.
The representation $\rho = (\rho_L, \rho_R)$ is admissible if and only if it is faithfull, has discrete image, and lifts to  some representation $\tilde{\rho}: \Gamma \rightarrow (\widetilde{G} \times \widetilde{G})_{/{\mathcal Z}}$ preserving a generic closed achronal subset of $\widehat{\mbox{Ein}}_2$ containing at least two points.

A $\rho$-admissible closed subset for an admissible representation $\rho$ is the projection in $\overline{\mbox{Ein}}_2$ of $\tilde{\rho}$-invariant generic closed achronal subset of $\widehat{\mbox{Ein}}_2$ containing at least two points.
\end{defin}

\emph{Problem: characterize admissible representations.\/}

\begin{thm}
\label{thm.admissible}
Let $\Gamma$ be a torsionfree group, and 
$\rho: \Gamma \rightarrow G \times G$ a faithfull representation.
Then, $\rho$ is admissible if and only if one the following occurs:
\begin{enumerate}

\item \emph{The abelian case:\/} $\rho(\Gamma)$ is a discrete subgroup of $A_{hyp}$, $A_{ext}$ or $A_{par}$ 
where:

-- $A_{hyp} =\{ (\exp(\lambda\Delta), \exp(\mu\Delta)) / \lambda, \mu \in {\mathbb R}\}$,

-- $A_{ext} =\{ (\exp(\lambda\Delta), \exp(\eta{H})) / \lambda, \eta \in {\mathbb R}\}$,

-- $A_{par} =\{ (\exp(\lambda{H}), \exp(\lambda{H})) / \lambda \in {\mathbb R}\}$.

\item \emph{The non-abelian case:\/} The left and right morphisms $\rho_L$, $\rho_R$ are faithfull with discrete image, and the marked surfaces $\rho_L(\Gamma)\backslash{\mathbb H}^2$, $\rho_R(\Gamma)\backslash{\mathbb H}^2$ are homeomorphic, i.e. there is a $\Gamma$-equivariant homeomorphism $f: {\mathbb H}^2 \rightarrow {\mathbb H}^2$ satisfying:

\[ \forall \gamma \in \Gamma, \; f \circ \rho_L(\gamma) = \rho_R(\gamma) \circ f\]

\end{enumerate}
\end{thm}

Observe that the restriction of an admissible representation to any non-trivial subgroup 
of $\Gamma$ is still admissible. 

The abelian case needs a study of the elementary case: hence
we postpone its proof to \cite{ba2}, and assume from now that $\Gamma$ is not abelian.
The main step in the proof of Theorem~\ref{thm.admissible} is to prove:

\begin{prop}
If $\rho$ is admissible, then the representations $\rho_L$, $\rho_R$ are faithfull.
\end{prop}

\preu
Assume by contradiction that the kernel $\Gamma_L$ of $\rho_L$ is not trivial, and that 
$\rho$ is admissible. Let $\widetilde{\Lambda}$ be a generic achronal $\tilde{\rho}(\Gamma)$-invariant 
closed subset of $\widehat{\mbox{Ein}}_2$. Let $\widetilde{\Lambda}^\pm$ be the up and low 
completions of $\widetilde{\Lambda}$.
Recall the description of the left and right foliations $\widehat{\mathcal G}_L$, 
$\widehat{\mathcal G}_R$ in remark~\ref{rk.einpsl}. Then, every $\tilde{\rho}(\gamma)$, 
for every $\gamma$ in $\Gamma_L$ preserves individually every leaf of $\widehat{\mathcal G}_L$. 
Hence, for every such a leaf, the intersection $\widetilde{\Lambda}^\pm \cap l$, if non-empty, 
is a $\tilde{\rho}(\Gamma)$-invariant closed interval (this intersection is connected since 
$\widetilde{\Lambda}^\pm$ is a topological circle). The extremities of this interval - maybe 
reduced to one point - project in ${\mathbb R}P^1_R$ as fixed points for every element of 
$\rho_R(\Gamma_L)$. 
Hence, the fixed point set $F_R$ of $\rho_R(\Gamma_L)$ in ${\mathbb R}P^1_R$ is not empty. 
On the other hand, an element of $G$ with three fixed point in ${\mathbb R}P^1$ is trivial, 
and the restriction of $\rho_R$ to $\Gamma_L$ is faithfull, since $\rho$ is faithfull: $F$ 
contains at most two points. The action of $\rho_R(\Gamma)$ on $\rho_R(\Gamma_L)$ permutes 
these two points, and an instant of reflexion is enough to realize that, since 
$\tilde{\rho}(\Gamma)$ preserves the chronological orientation, every element of 
$\rho_R(\Gamma)$ must preserves every element of $F_R$. 

Assume that if $F_R$ is reduced to one point. Let $r_0$ be the corresponding leaf of the right foliation $\widehat{\mathcal G}_R$. The argument above implies that 
for every leaf $l$ of $\widehat{\mathcal G}_L$, $l \cap \widetilde{\Lambda}^\pm$ is either empty, either a point in $l \cap r_0$, or a closed interval projecting on the entire ${\mathbb R}P^1_R$. If the last case occurs, then $\widetilde{\Lambda}^\pm$ is pure lightlike: it means that $\widetilde{\Lambda}$ is elementary, more precisely, that it is conical or extreme. 
Actually, the conical case would imply that $F_R$ contains two points (the non-causally related 
extremities). Hence, $\widetilde{\Lambda}$ is a lightlike segment contained in $r_0$.
The projections in ${\mathbb R}P^1_L$ of the two extremities of this lightlike segment are distinct $\rho_L(\Gamma)$-fixed points. It follows that $\rho(\Gamma)$ is contained in a conjugate of $A_{ext}$. It is absurd since $\Gamma$ is not abelian.

Therefore, $F_R$ contains two points. After conjugacy, $\rho_R(\Gamma)$ is 
contained in the $1$-parameter group 
$\{ \exp(\lambda\Delta) / \lambda \in \mathbb R \}$. Since $\Gamma$ is not abelian, 
it means that $\rho_R$ is not injective too! Apply once more all the arguments above: 
it follows that $\rho_L(\Gamma)$ admits two distinct fixed points in ${\mathbb R}P^1_L$. 
In other words, 
$\rho(\Gamma)$ is contained in the abelian group $A_{hyp}$. Contradiction.\fin

\begin{cor}
If $\rho$ is admissible, $\rho_L$ and $\rho_R$ have discrete image in $G$.
\end{cor}

\preu
Since $\rho$ is admissible, $\rho_L(\Gamma)$ has no elliptic element. Hence, if $\rho_L(\Gamma)$ 
is not discrete, the neutral component of its closure is a the stabilizer of one or two 
points in ${\mathbb R}P^1_L$. These fixed points are permuted, and actually preserved, 
by every $\rho_L(\gamma)$. Since $\rho_L$ is faithfull and $\Gamma$ is not abelian, 
it means that there is only one fixed point, i.e. $\rho_L(\Gamma)$ is contained in a 
solvable group conjugate to Aff. 
It follows that $\rho_R(\Gamma)$ is solvable too, hence contained also up to conjugacy in 
Aff. The elements of the commutator subgroup $\rho([\Gamma,\Gamma])$, which is not trivial 
since $\Gamma$ is not abelian, are parabolic elements. 
The representation $\rho: [\Gamma, \Gamma] \rightarrow G \times G$ is then
an admissible representation of an abelian group, case which is studied in \cite{ba2}.
It follows that up to conjugacy  
$\rho([\Gamma,\Gamma])$ is contained in $A_{par}$. Since $\rho$ has discrete image, 
$\rho_L([\Gamma,\Gamma]) \approx \rho_R([\Gamma, \Gamma])$ is a cyclic group, preserving
a copy of the affine line ${\mathbb R}$ in ${\mathbb R}P^1_{L,R}$, and acting 
on this line as translations. But the action by conjugacy of $\Gamma$
on $[ \Gamma, \Gamma ]$ induces an action by homotheties on this dicrete goup
of translations; hence either $\Gamma$ is contained in $A_{par}$, or $[\Gamma, \Gamma]$
is trivial. In the later case, $\Gamma$ is abelian: contradiction. In the former case,
we can apply the arguments above: $\rho_L(\Gamma) = \rho_R(\Gamma)$ are cyclic groups
of translations.
\fin

Observe that if $\rho$ is admissible, then the future extension $\widetilde{\Lambda}^+$ 
(for example) of a $\tilde{\rho}(\Gamma)$-invariant closed generic achronal subset 
is a topological circle which defined a \emph{monotone\/} map 
$f: {\mathbb R}P^1_L \rightarrow {\mathbb R}P^1_R$ which is equivariant:

\[ \forall \gamma \in \Gamma, \; f \circ \rho_L(\gamma) = \rho_R(\gamma) \circ f\]

Here, by monotone map, we mean a relation which can send a single point in ${\mathbb R}P^1_L$ on a closed segment of ${\mathbb R}P_R^1$, and such that every $f^{-1}(x)$ is a point or a closed segment of ${\mathbb R}P^1_L$. Equivalently, a monotone map is the quotient of a monotone relation $\tilde{f}: \widetilde{\mathbb R}P^1_L \rightarrow \widetilde{\mathbb R}P^1_R$
which commutes with the Galois groups:

\[ \tilde{f}\circ \delta_L = \delta_R \circ \tilde{f} \]

The projective line ${\mathbb R}P^1_L$, ${\mathbb R}P^1_R$ has to be considered as the conformal boundary of the hyperbolic plane ${\mathbb H}^2$ on which are acting respectively $\rho_L(\Gamma)$, $\rho_R(\Gamma)$. Theorem~\ref{thm.admissible} follows from the well-known fact that the existence of a homeomorphism between the marked surfaces $\Sigma_L= \rho_L(\Gamma)\backslash{\mathbb H}^2$, $\Sigma_R=\rho_R(\Gamma)\backslash{\mathbb H}^2$ is equivalent to the existence of an equivariant monotone map as above. \fin


\rque
\label{rk.euler}
In the non-abelian case, there is a much more elegant and concise formulation of Theorem~\ref{thm.admissible}, using the notion of bounded Euler cohomology class, which is exactly the obstruction for the existence of a equivariant monotone map semi-conjugating two actions of a given group on the circle (see \cite{ghyseuler1, ghyseuler2}):

\begin{thm}
\label{thm.admissibleeuler}
Let $\Gamma$ a non-abelian group without torsion, and $\rho: \Gamma \rightarrow G \times G$ a faithfull representation. Then, $\rho$ is admissible if and only if the left and right representations $\rho_L$, $\rho_R$ are faithfull discrete representations with the same Euler bounded cohomology class.\fin
\end{thm}

\erque

\subsection{Minimal invariant achronal subsets}
In almost all this section, $\Gamma$ is a non-abelian group, and $\rho: \Gamma \rightarrow G \times G$ an admissible representation. 

\begin{defin}
$\overline{\Lambda}(\rho)$ is the closure of the set of attractive fixed points in $P(E)$.
\end{defin}

Observe that attractive fixed points in $P(E)$ of elements of $G$ belong to $\overline{\mbox{Ein}_2}$. Hence, $\overline{\Lambda}(\rho)$ is contained in $\overline{\mbox{Ein}}_2$.

\begin{thm}
\label{thm.minimal}
Let $\Gamma$ be a non-abelian torsion-free group, and $\rho: \Gamma \rightarrow G \times G$.
Then, every $\rho(\Gamma)$-invariant closed subset of $P(E)$ contains $\overline{\Lambda}(\rho)$.
\end{thm}

\begin{cor}
\label{cor.cor}
Let $(\Gamma,\rho)$ be pair satisfying the hypothesis of Theorem~\ref{thm.minimal}. Then,
$\overline{\Lambda}(\rho)$ is a $\rho(\Gamma)$-invariant generic nonelementary achronal subset of 
$\overline{\mbox{Ein}}_2$. Furthermore, for every $\rho(\Gamma)$-invariant closed achronal subset 
$\Lambda$ in $\mbox{Ein}_2$, the invisibility domain $E(\Lambda)$ projects injectively in 
$\overline{{\mathbb A}{\mathbb D}{\mathbb S}}$ inside $E(\overline{\Lambda}(\rho))$.
\fin
\end{cor}

The essential step for the proof of Theorem~\ref{thm.minimal} is:

\begin{lem}
\label{le.strongirreducible}
If $\rho(\Gamma)$ is admissible, and does not preserve a point in ${\mathbb A}{\mathbb D}{\mathbb S}$, 
then $\rho(\Gamma)$ is strongly irreducible, i.e. for every finite index subgroup $\Gamma' \subset \Gamma$, there is not $\rho(\Gamma')$-invariant proper projective subspace in $S(E)$.
\end{lem}

\preu
Assume by contradiction that $\rho$ is admissible, and that some finite index subgroup $\Gamma' \subset \Gamma$ preserves a projective subspace $S(F) \subset S(F)$, where $F \neq E$ is a non-trivial linear subspace of $E$. Observe that $F^\perp$ is also preserved by $\rho(\Gamma')$. 

If $F$ and $F^\perp$ does not contain $Q$-isotropic vectors, 
then $\Gamma'$ is contained in $O(2) \times O(2)$. 
It is impossible since elements of $\Gamma$ are non-elliptic synchronized. 

Hence, $S(F) \cap \mbox{Ein}_2$ or $S(F^\perp) \cap \mbox{Ein}_2$ is not empty. Such an intersection is either a $\rho(\Gamma')$-fixed point in $\overline{\mbox{Ein}}_2 \approx {\mathbb R}P^1_L \times {\mathbb R}P^1_R$,  a lightlike geodesic, or an invariant round circle $\partial x_0^\ast$. The last case is impossible, since $x_0$ would be a $\rho(\Gamma)$ fixed point, situation that we have excluded by hypothesis. In the two other remaining cases,  after switching if necessary the left and right factors, we obtain that $\rho_L(\Gamma')$ admits a global fixed point in ${\mathbb R}P^1_L$. It is impossible since $\rho_L(\Gamma)$ is a non-abelian discrete subgroup of $G$.\fin

\preud{thm.minimal}
We first observe that, according to Theorem~\ref{thm.admissible}, and since any non abelian discrete subgroup of $G$ admits at least one hyperbolic element, $\rho(\Gamma)$ is proximal, i.e. $\overline{\Lambda}(\rho)$ is non-empty. Then, Theorem~\ref{thm.minimal} is an immediate corollary of lemma~\ref{le.strongirreducible} and Lemma $2.5-(2)$ of \cite{benoistconvex}.\fin.

\subsection{Convexity and causality}
\label{subsec.benoistconvexe}
In this {\S}, we consider a representation $\rho: \Gamma \rightarrow G \times G$, where $\Gamma$ is not abelian.

\begin{defin}
\label{def.proximal}
An element of $\mbox{GL}(E)$ is proximal if its action on $P(E)$ admits an attractive fixed point. It is positively proximal if its action on $S(E)$ has two attractive fixed points (one opposite to the other). 

A faithfull representation $\rho: \Gamma \rightarrow \mbox{GL}(E)$ is positively proximal if 
$\rho(\Gamma)$ contains at least a proximal element, and that every proximal element of $\rho(\Gamma)$ is positively proximal.
\end{defin}

The main result of \cite{benoistconvex} (Propositions $1.1$, $1.2$) is:

\begin{thm}
A strongly irreducible representation is positively proximal if and only if it preserves a proper convex domains in $P(E)$.\fin
\end{thm}

It follows that in our case, if $\rho$ is admissible, then it is positively proximal (observe that this statement is true, even if $\rho$ is not strongly irreducible, i.e. preserves a point in AdS).
We wonder here about the inverse statement: is any positively proximal $\rho$ admissible?

In the following, we are using results in \cite{benoistconvex} which are established for strongly irreducible representations, but these results are easily checked when $\rho(\Gamma) \subset\mbox{SO}_0(2,2)$ admits a fixed point in AdS.

Let ${\mathcal F}(E)$ be the flag variety, i.e. the space of pairs $([u],[u^\ast])$ in $P(E) \times P(E^\ast)$ such that $u^\ast(u) = 0$. Here, we can of course define ${\mathcal F}(E)$ as the space of pairs of $Q$-orthogonal elements of $P(E)$: ${\mathcal F}(E) =\{ ([u], [v]) \in P(E) \times P(E) / \langle u \mid v \rangle = 0 \}$. The group $\mbox{SO}_0(2,2)$ acts naturally on it, by the diagonal action. The closure in ${\mathcal F}(E)$  of the set of attractive fixed points of elements of $\rho(\Gamma)$ is: $\Lambda^{\mathcal F} = \{ ([u],[u]) \in {\mathcal F}(E) / [u] \in \Lambda^{\mathbb P} \}$. 

\begin{prop}[Lemma $2.5-(3)$ of \cite{benoistconvex}]
\label{pro.Fminimal}
Any $\rho(\Gamma)$-invariant subset of ${\mathcal F}(E) $ contains $\Lambda^{\mathcal F}$.\fin
\end{prop}

Moreover, the statement $(3)-d$ expresses here:

\begin{prop}
\label{pro.transverse}
$\Lambda^{\mathbb P} \times \Lambda^{\mathbb P}$ contains a dense subset $Y$, which is \emph{transverse}, i.e. for every $([u], [v])$ in $Y$, the scalar product $\langle u \mid v \rangle$ is nonzero.\fin
\end{prop}

Observe that in the statement above, we cannot define the sign of $\langle u \mid v \rangle$, since $[u]$, $[v]$ are only elements of $P(E)$. But there is a sign if we lift all the picture in $S(E)$. Define $\Lambda^{S}$ as the preimage in $S(E)$ of $P(E)$. It can  also be defined as the closure of the set of attractive fixed points of elements of $\rho(\Gamma)$.

\begin{lem}[Proposition $3.15$ in \cite{benoistconvex}]
\label{le.dansS}
A strongly irreducible representation $\rho$ is positively proximal 
if and only if the action of $\rho(\Gamma)$ on $\Lambda^S$ is not minimal. 
If it is the case, $\Lambda^S$ is the union of two disjoint minimal closed invariant subset $\Lambda^S_1$ and $\Lambda^S_2 = -\Lambda^S_1$.\fin
\end{lem}

\begin{lem}
\label{le.positive}
The closed minimal subset $\Lambda_1^S$ is \emph{positive,\/} i.e. we have the following alternative:

\begin{enumerate}
\item for every element $([u],[v])$ in $\Lambda_1^S \times \Lambda^S_1$, we have $\langle u \mid v \rangle \leq 0$, or
\item for every element $([u],[v])$ in $\Lambda_1^S \times \Lambda^S_1$, we have $\langle u \mid v \rangle \geq 0$.
\end{enumerate}
\end{lem}

\preu
It follows from the proof of Proposition $3.11$ in \cite{benoistconvex}.\fin

In the first case of the alternative of lemma~\ref{le.positive}, $\Lambda_1^S$ is an achronal closed subset of $\mbox{Ein}_2$. According to Proposition~\ref{pro.transverse}, it contains at least a pair of non-causally related points: $\Lambda_1^S$ is generic. It follows that $\rho$ is admissible.

But, in the second case, $\rho$ is not admissible! Every pair of points in $\Lambda_1^S$ is causally related, and most of these pairs are strictly causally related.
The situation can be entirely understood in the light of remarks~\ref{autreads}, \ref{2copies}: we have to consider $\overline{\mbox{Ein}}_2$ not as the Klein boundary of $\overline{{\mathbb A}{\mathbb D}{\mathbb S}}$, but as the Klein boundary of the complementary AdS copy, which is the projection of $\{ Q = +1 \}$ equipped with the restriction of $-Q$. We obtain the notion of $-$-admissible representations $\rho: \Gamma \rightarrow \mbox{SO}_0(2,2)$: the representations conjugate to admissible representations in the previous meaning by an anti-isometry of $E$ permuting $\{ Q > 0 \}$ and $\{ Q < 0 \}$.

\begin{prop}
\label{pro.-}
Let $\Gamma$ be a torsionfree nonabelian group, and $\rho: \Gamma \rightarrow \mbox{SO}_0(2,2)$ a faithfull representation with discrete image. Then, $\rho$ is positively proximal if and only if it is admissible or $-$-admissible.\fin
\end{prop}

\section{Cauchy complete globally hyperbolic AdS spacetimes}
\label{sec.cauchycomplet}
Let $M$ be a $3$-dimensional manifold equipped with a lorentzian metric of constant 
curvature $-1$, i.e. locally
modeled on $\mbox{AdS}$. Let $p: \widetilde{M} \rightarrow M$ be the 
universal covering, and $\Gamma$
the fundamental group of $M$, considered as the group of desk transformations 
of $p$. We recall some
basic facts of $(G,X)$-structure's theory, applied in our context:

- There is a developing map ${\mathcal D}: \widetilde{M} \rightarrow \widetilde{\mbox{AdS}}$, which is a 
local homeomorphism;

- There is a holonomy morphism $\widetilde{\rho}: \Gamma \rightarrow \widehat{SO}(2,2)$ 
for which $\mathcal D$ is 
equivariant (here, $\widehat{SO}(2,2)$ is the isometry group of $\widetilde{\mbox{AdS}}$). 
More precisely, for every $\gamma$
in $\Gamma$, $\rho(\gamma) \circ {\mathcal D} = {\mathcal D} \circ \gamma$.

Assume that $M$ is globally hyperbolic, and \emph{Cauchy complete}, i.e. admits a Cauchy surface $S$ such that the restriction of the ambient lorentzian metric on $S$ is a complete riemannian metric. Then, the preimage $\widetilde{S}$ in $\widetilde{M}$ is a Cauchy surface for $\widetilde{M}$. According to Proposition~\ref{pro.graphe}, the restriction to $\widetilde{S}$ of the developping map $\mathcal D$ is an embedding, and the image ${\mathcal D}(\widetilde{S})$ is the graph of a contracting map from ${\mathbb D}^2$ into $\mathbb R$. Since $\widetilde{S}$ is a Cauchy surface of $\widetilde{M}$, the image of $\mathcal D$ must be contained in the Cauchy development $T({\mathcal D}(\widetilde{S})) = {\mathcal C}({\mathcal D}(\widetilde{S}))$. Since local homeomorphisms between subintervals of $\mathbb R$ are always injective, the restriction of $\mathcal D$ to any timelike geodesic is injective: it follows that $\mathcal D$ is injective. 

Now, observe that the holonomy group $\widetilde{\rho}$ preserves 
${\mathcal D}(\widetilde{S})$ and $E(\partial {\mathcal D}(\widetilde{S}))$. 
Hence, according to Theorem~\ref{noneleOK} and Proposition~\ref{pro.gh}:

\begin{prop}
\label{pro.maximalgh}
Any globally hyperbolic AdS spacetime $M$,with complete Cauchy surface $S$, embeds isometrically in the spacetime $M_{\partial \widetilde{S}}(\widetilde{\rho}(\Gamma))$, where $\widetilde{\rho}: \Gamma \rightarrow \widehat{\mbox{SO}}(2,2)$ is the holonomy morphism, and $\widetilde{S}$ the image by the developping map of a Cauchy surface in $\widetilde{M}$.\fin
\end{prop}

Recall the notion of maximal global hyperbolicity (definition~\ref{def.maxigh}).

\begin{cor}
\label{cor.maxigh}
The maximal Cauchy complete globally hyperbolic spacetimes are the quotient spacetimes $M_{\widetilde{\Lambda}}(\Gamma)$, where $\widetilde{\Lambda}$ is a generic achronal topological circle in $\widehat{\mbox{Ein}}_2$, and $\Gamma$ a discrete torsion-free subgroup of $\widehat{\mbox{SO}}(2,2)$ preserving $\widetilde{\Lambda}$.
\end{cor}

\preu
In the light of Propositions~\ref{pro.gh} and \ref{pro.maximalgh}, the only remaining point to check is the fact that every spacetime $M_{\widetilde{\Lambda}}(\Gamma)$ is indeed Cauchy complete. 

Let ${\mathcal N}(\Lambda)$ be the subset of ${\mathcal T}^+$ formed by pairs $(x,y)$ such that $y$ belongs to the future connected component $C^+(\Lambda)$ of $\partial\mbox{Conv}(\Lambda) \setminus \Lambda$. Then, $x$ must belong the past connected component $\partial^-E(\Lambda)$ of ${\partial}E(\Lambda) \setminus \Lambda$ (see Proposition~\ref{pro.futfut}).
Moreover, $y$ must belong to the spacelike part $C_{spa}^+(\Lambda)$ of $C^+(\Lambda)$. 

Equivalently, ${\mathcal N}(\Lambda)$ is the Gauss graph of $C^+_{spa}(\Lambda)$ (see Definition~\ref{def.gaussgraph}; the surface $C^+_{spa}(\Lambda)$ is actually past-convex, i.e. the roles of $x$, $y$ has been permuted, but without incidence on the results of \S~\ref{sub.isomh2}).

Consider now the cosmological time $\tau$ on $E(\Lambda)$; more precisely, the level set $\Sigma = \tau^{-1}(\pi/4)$. For any $p$ in $\Sigma$, the surface formed by points $x$ in the past of $p$ and at AdS-distance $\pi/4$ is concave. Since $\partial^-E(\Lambda)$ is convex,
it follows that there is one and only one point $x(p)$ in $\partial^-E(\Lambda)$ at distance $\pi/4$ from $p$. Let $y(p)$ be the first intersection in the future of $p$ of the future oriented geodesic $(p, x(p))$ with $\partial\mbox{Conv}(\Lambda)$. Then, $y(p)$ belongs to $C^+_{spa}(\Lambda)$, and $(x(p), y(p))$ belongs to ${\mathcal N}(\Lambda)$. Hence, $p$
is equal to $\frac{1}{2}(x(p) + y(p))$. Inversely, for every $(x,y)$ in ${\mathcal N}(\Lambda)$,
$\frac{1}{2}(x+y)$ belongs to $\Sigma$. 

We have thus defined a homeomorphism $(x,y): \Sigma \rightarrow {\mathcal N}(S)$, differentiable almost everywhere\footnote{This map is actually $C^1$. See \cite{BenBon2}.} such that $p = \frac{1}{2}(x(p) + y(p))$. The norm in AdS of a tangent vector $\frac{1}{2}(u+v)$ of $\Sigma$ is $\frac{1}{4}(\mid u \mid^2 + \mid v \mid^2) + \frac{1}{2}\langle u \mid v \rangle$, whereas the norm of the image of this tangent vector by the differential of $(x,y)$ is
$\frac{1}{4}(\mid u \mid^2 + \mid v \mid^2)$. We therefore obtain as a corollary\footnote{The reader can now understand the choice of the factor $1/4$ introduced in the definition of the lorentzian metric on $\mathcal T$.} of Lemma~\ref{lem.uvpositif} the following key point:

\emph{The map $(x,y): \Sigma \rightarrow {\mathcal N}(S)$ decreases the distance.\/}

Consider a Cauchy sequence $p_n$ in $\Sigma$. Denote $(x_n, y_n) = (x(p_n), y(p_n))$. Then, according to the key point we have just established, $(x_n, y_n)$ is a Cauchy sequence in ${\mathcal N}(S)$.
According to lemma~\ref{lem.isomconvex}, the $\pi_G$-projection of $(x_n, y_n)$ is also a Cauchy sequence in $Q_G$. But, according to Proposition~\ref{pro.QG}, $Q_G$ is complete.
Hence, if $\Delta_n$ denotes the timelike geodesic containing $x_n$ and $y_n$, the $\Delta_n$ converges to a timelike geodesic $\Delta_\infty$. This geodesic intersects $\partial^-E(\Lambda)$ at an unique point $x_\infty$, and intersects $C^+(\Lambda)$ at an unique point $y_\infty$. Then, the $x_n$ converge to $x_\infty$, and the $y_n$ converge to $y_\infty$.
It follows that the $p_n$ converge to $p_\infty = \frac{1}{2}(x_\infty + y_\infty)$. This point $p_\infty$ belongs to $\Sigma$.

Therefore, the spacelike surface $\Sigma$ is Cauchy complete. 
Since it is $\Gamma$-invariant, its projection in $M_\Lambda(\Gamma)$ is a spacelike Cauchy complete surface. According to Proposition~\ref{pro.Ecosmologique}, $\Sigma$ is also a Cauchy surface in $M_\Lambda(\Gamma)$.
The proposition is proved.
\fin

\rque
An interesting - and difficult - problem is to find a \emph{smooth\/} 
Cauchy complete Cauchy surface in $M_\Lambda(\Gamma)$.
\erque

\end{document}